\documentclass[a4paper,10pt]{article}
\everymath{\displaystyle}
\usepackage{amssymb,amsmath,amsthm}
\usepackage{graphicx}

\def\N{\mathbb N}
\def\R{\mathbb R}
\def\T{\mathbb T}
\def\Zi{{\mathbb Z}\setminus {\cal I}}
\usepackage{amsfonts}
\usepackage{hyperref}
\usepackage{latexsym}
\usepackage{amsthm}
\usepackage{amsmath}
\usepackage{amssymb}
\usepackage{amsopn}
\usepackage{geometry}
\usepackage{amscd}
\usepackage{mathrsfs}
\usepackage{dsfont}
\usepackage[english]{babel}
\usepackage{caption}
\font\teneufm=eufm10 \font\seveneufm=eufm7 \font\fiveeufm=eufm5
\newfam\eufmfam
\textfont\eufmfam=\teneufm \scriptfont\eufmfam=\seveneufm
\scriptscriptfont\eufmfam=\fiveeufm
\def\eufm#1{{\fam\eufmfam\relax#1}}

\def\Ngot{{\eufm N}}
\def\ugot{{\eufm u}}
\def\vgot{{\eufm v}}

\newcommand{\Om}{\Omega}

\newcommand{\aaa}{{\mathtt{a}}}

\newenvironment{pf}{\noindent{\sc Proof}.\enspace}{\rule{2mm}{2mm}\smallskip}
\newenvironment{pfn}{\noindent{\sc Proof} \enspace}{\rule{2mm}{2mm}\smallskip}
\newtheorem{theorem}{Theorem}[section]
\newtheorem{proposition}{Proposition}[section]
\newtheorem{lemma}{Lemma}[section]
\newtheorem{corollary}{Corollary}[section]

\newtheorem{remark}{Remark}[section]
\newtheorem{remarks}{Remark}[section]
\newtheorem{definition}{Definition}[section]
\newcommand{\be}{\begin{equation}}
\newcommand{\ee}{\end{equation}}
\newcommand{\teta}{\theta}
\newcommand{\om}{\omega}
\newcommand{\e}{\varepsilon}

\renewcommand{\a }{\alpha }
\renewcommand{\b }{\beta }
\newcommand{\s }{\sigma }
\newcommand{\ii }{{\rm i} }
\renewcommand{\d }{\delta }

\def\Ta{\mathcal T}
\newcommand{\g }{\gamma}

\renewcommand{\l }{\lambda }

\renewcommand{\t }{\tau }
\renewcommand{\o }{\omega }
\renewcommand{\O }{\Omega }
\newcommand{\C}{\mathbb{C}}
\newcommand{\Z}{\mathbb{Z}}

\newcommand{\mm}{{\rm m}}

\newcommand{\pluto}{\mathtt j}
\newcommand{\norma}{\|}

\newcommand{\rro}{{\rho}}

\newcommand{\cc}{{1}}
\newcommand{\CC}{{6}}
\newcommand{\dueCC}{{12}}
\newcommand{\duecc}{{2}}
\newcommand{\trecc}{{3}}
\newcommand{\quattrocc}{{4}}

\newcommand{\CCterzi}{{2}}
\newcommand{\CCquarti}{{3/2}}

\newcommand{\OO}{\mathtt{\O}}
\newcommand{\oo}{\mathtt{\o}}
\newcommand{\ot}{\vec{\mathtt{\o}}}
\newcommand{\XX}{\mathcal{X}}
\newcommand{\PP}{\mathcal{P}}

\usepackage{color}

\title{\bf KAM  for reversible
 derivative wave equations}

\author{Massimiliano Berti, Luca Biasco, Michela Procesi}

\date{}

\begin{document}

\maketitle

\noindent {\bf Abstract:} We prove the existence of Cantor
families of small amplitude, analytic, linearly stable quasi-periodic solutions
of reversible derivative wave equations.
\\[1.5mm]
{\sc 2000AMS subject classification}: 37K55, 35L05.


\section{Introduction}

 An important question in KAM theory for PDEs concerns equations with
derivatives  in the nonlinearity. Only few results are known, mainly restricted to dispersive equations.
For Hamiltonian perturbations of KdV, existence and stability of quasi-periodic solutions  was first proved  by Kuksin  \cite{K2}-\cite{k1}
in the  late '90, see also  Kappeler-P\"oschel \cite{KaP}. This approach has been
recently extended by  Liu-Yuan
\cite{LY} for Hamiltonian DNLS and by Zhang-Gao-Yuan \cite{ZGY}
for the reversible DNLS equation $ \ii u_t + u_{{\mathtt x}{\mathtt x}} + |u_{\mathtt x}|^2 u = 0 $.

The {\it derivative} nonlinear wave equation (DNLW), which is {\it not} dispersive,
is {\it excluded}  by  these approaches (for semilinear wave equations see
\cite{k1},  \cite{W1},  \cite{B1}, \cite{Po3}, \cite{CY}, \cite{BB12}).
Existence of periodic solutions (without stability) for  the
derivative Klein-Gordon equation
\be\label{yttyxx} {\mathtt y}_{tt} - {\mathtt
y}_{{\mathtt x}{\mathtt x}} + \mm {\mathtt y} + {\mathtt y}_t^2  =
0 \, , \quad \mm > 0 \, , \quad {\mathtt x} \in \T:=
\R/2\pi\Z \, ,
\ee
was first proved by Bourgain in \cite{B2}, extending
the approach of  Craig-Wayne in \cite{CW}.
Then Craig  \cite{C} focused  on the natural question of establishing similar results
for more general  derivative wave equations
\be\label{DNLW}
{\mathtt y}_{tt} - {\mathtt y}_{{\mathtt x}{\mathtt x}} + \mm {\mathtt y} =
g({\mathtt x}, {\mathtt y}, {\mathtt y}_{\mathtt x}, {\mathtt
y}_t) \, , \quad {\mathtt x} \in \T  \, ,
\ee
asking, for example, if $ {\mathtt y}_{tt} - {\mathtt y}_{{\mathtt x}{\mathtt x}} = {\mathtt y}_{\mathtt x}^3 $
possesses periodic solutions,
see  \cite{C}, section 7.3.

In \cite{BBP1} we recently extended KAM theory for the Hamiltonian  model
$$
{\mathtt y}_{tt} - {\mathtt y}_{{\mathtt x}{\mathtt x}} + \mm
{\mathtt y} + f(D {\mathtt y})=0 \, , \quad \mm >  0 \, , \quad D
:=\sqrt{-\partial_{{\mathtt x}{\mathtt x}}+\mm} \, , \quad
{\mathtt x} \in \T \, .
$$
This kind of  pseudo-differential equations were
introduced by Bourgain \cite{B1} and Craig  \cite{C} as models to
study the effect of derivatives versus dispersive phenomena. Clearly
 \cite{BBP1} does not apply to the derivative wave equations \eqref{DNLW}, which are not  Hamiltonian.

\smallskip

In order to prove existence of periodic/quasi-periodic solutions for
\eqref{DNLW}, conditions on the nonlinearity $ g $ have to be necessarily imposed.
For
example,  \eqref{DNLW} with the nonlinear friction term $ g =
{\mathtt y}_t^3 $ has no nontrivial smooth periodic/quasi-periodic solutions,
see Proposition \ref{nonex}. This case may be
ruled out by assuming the reversibility condition
\be\label{geven}
g({\mathtt x}, {\mathtt y}, {\mathtt y}_{\mathtt x}, - {\mathtt v}
) = g({\mathtt x}, {\mathtt y}, {\mathtt y}_{\mathtt x},  {\mathtt
v} )
\ee
satisfied for example by  \eqref{yttyxx}. Under condition
\eqref{geven} the equation \eqref{DNLW} is time-{\it reversible},
namely the associated  first order system
\be\label{firstor}
{\mathtt y}_t =  {\mathtt v}  \, , \quad {\mathtt v} _{t} =
{\mathtt y}_{{\mathtt x}{\mathtt x}} - \mm {\mathtt y} +
g({\mathtt x}, {\mathtt y}, {\mathtt y}_{\mathtt x}, {\mathtt v})
\ee
is reversible with respect to the involution
\be\label{invol1}
S ( {\mathtt y},  {\mathtt v}  ) := ( {\mathtt
y}, -  {\mathtt v}  ) \, \, , \quad S^2 = I \, .
\ee
For finite-dimensional systems it is  known (since Moser \cite{Mo67})
  that reversibility may replace the Hamiltonian structure
 in order to allow the existence of quasi-periodic solutions, see also
 Arnold \cite{Arn} and Sevryuk \cite{Se1}.
However, for \eqref{DNLW}, it 
 is  {\it not} sufficient. For
example
$ {\mathtt y}_{tt} - {\mathtt y}_{{\mathtt x}{\mathtt x}}  =
{\mathtt y}_{\mathtt x}^3 $
is time reversible
but it has {\it no} smooth periodic/quasi-periodic solutions except the
constants (in Proposition \ref{nonex} we exhibit more general time-reversible nonlinearities for which DNLW has only trivial quasi-periodic solutions).
In order to find
quasi-periodic solutions we also require the ``space-reversibility''
 assumption 
\be\label{REALITY}
g( - {\mathtt x}, {\mathtt y}, - {\mathtt y}_{\mathtt x} ,  {\mathtt v}  ) =
 g( {\mathtt x}, {\mathtt y},  {\mathtt y}_{\mathtt x} ,  {\mathtt v}  )
\ee
which  rules out  nonlinearities like
$ {\mathtt y}_{\mathtt x}^3, {\mathtt y}_{\mathtt x}^5, \ldots $. 
Actually,  condition \eqref{REALITY} is as natural as \eqref{geven}.
Indeed, for the wave equation  \eqref{DNLW},  the role of
 time and space variables $ (t, {\mathtt x}) $ is highly
symmetric, and, considering $ {\mathtt x} $ ``as time" (spatial
dynamics idea)
\eqref{REALITY} is nothing but the corresponding reversibility condition
and terms like $ {\mathtt y}_{\mathtt x}^3, {\mathtt y}_{\mathtt x}^5, \ldots  $ are frictions.

\smallskip

In this paper we prove  {\it existence} and {\it stability} of {\it analytic}
quasi-periodic solutions for
derivative wave equations \eqref{DNLW} satisfying \eqref{geven}, \eqref{REALITY},
see Theorem \ref{thm:DNLW}.
By the above considerations,  this is a
very natural class of DNLW equations which may admit  quasi-periodic solutions.
After  Theorem \ref{thm:DNLW} we shall further comment on the assumptions.
These results were presented in the note \cite{BBP2}.


\smallskip

Before describing our main results, we mention the classical
bifurcation theorems of Rabinowitz \cite{R67} about periodic
solutions  (with period $T\in \pi\mathbb Q$) of dissipative forced
derivative wave equations
$$
{\mathtt y}_{tt} - {\mathtt y}_{{\mathtt x}{\mathtt x}}  + \a
{\mathtt y}_t + \e F({\mathtt x},t, {\mathtt y}, {\mathtt y}_{\mathtt x},
{\mathtt y}_t ) = 0 \, ,  \quad {\mathtt x} \in [0,\pi] 
$$
with Dirichlet boundary conditions, and in \cite{R68} with a
fully-non-linear forcing  term
$F=F({\mathtt x},t, {\mathtt y}, {\mathtt y}_{\mathtt x},
{\mathtt y}_t ,$ $ {\mathtt y}_{tt}, {\mathtt y}_{t{\mathtt x}},
{\mathtt y}_{{\mathtt x}{\mathtt x}}).$
Note that for forced PDEs the nonlinearity does not need to
be reversible or Hamiltonian.

\subsection{Main results}

We consider derivative wave equations \eqref{DNLW} where $ \mm > 0
$, the nonlinearity
$ g : \T \times {\cal U} \to \R $, $  {\cal U} \subset \R^3 $ open  neighborhood  of $ 0 $,
is real analytic and satisfies the
assumptions \eqref{geven}, \eqref{REALITY}. We require $ g $
to vanish at least quadratically at $ ({\mathtt y},{\mathtt
y}_{\mathtt x},  \mathtt v ) = (0,0,0 ) $, namely
$$
g({\mathtt x}, 0,0,0)  = (\partial_{\mathtt
y} g)({\mathtt x}, 0,0,0 ) = (\partial_{{\mathtt y}_{\mathtt x}}
g)({\mathtt x}, 0,0,0 ) = (\partial_{\mathtt v} g)({\mathtt x},
0,0,0 ) = 0 \, .
$$
Because of \eqref{geven}, it is natural to
look for ``reversible"  solutions,  namely such that
$ {\mathtt y}(t, {\mathtt x} ) $ is even and $  {\mathtt v}
(t,{\mathtt x}) $ is odd in time, and, because of
\eqref{REALITY}, it is natural to restrict to solutions which are
even is $ {\mathtt x} $ (standing waves).
Hence we  look for quasi-periodic
solutions  of \eqref{DNLW}
satisfying
\be\label{simmet} 
{\mathtt y}(t,{\mathtt x}) = {\mathtt
y}(t,-{\mathtt x}) \, , \ \forall t \, , \quad {\mathtt
y}(-t,{\mathtt x}) = {\mathtt y}(t,{\mathtt x}) \, , \  \forall
{\mathtt x \in \T } \, . 
\ee
For every finite choice of the {\it
tangential} sites $ {\cal I}^+ \subset \N \setminus \{0 \} $,  the
linear Klein-Gordon equation
\be\label{LW} {\mathtt y}_{tt} -
{\mathtt y}_{{\mathtt x}{\mathtt x}} + \mm {\mathtt y} = 0  \, ,
\quad {\mathtt x} \in \T  \, , \ee
possesses  the family of
quasi-periodic standing wave solutions
\be\label{ylineare}
{\mathtt y} = {\mathop\sum}_{j \in {\cal I}^+} \sqrt{8 \xi_j} \, \l_j^{-1}
\cos ( \l_j  \, t ) \cos (j {\mathtt x}) \, , \quad \l_j :=
\sqrt{j^2 + \mm} \, ,
\ee
parametrized by the ``actions'' $ \xi_j \in \R_+ $ and with linear frequencies of oscillations 
$ \bar\om:=(\l_j)_{j\in{\cal I}^+}$. 

In order to continue such solutions for the nonlinear equation \eqref{DNLW}
--as it is well known in KAM theory--
the leading term of the  nonlinearity $ g $ has to satisfy some  non-degeneracy condition
so that the ``action-to-frequency" map ``twists".
For definiteness, we have focused on nonlinearities
\be\label{gleading}
g =g^{(=3)}(\mathtt y, \mathtt y_{\mathtt x}, \mathtt y_t)  +
g^{(\geq 5)}
(\mathtt x, \mathtt y, \mathtt y_{\mathtt x}, \mathtt y_t)
\ee
with cubic leading term
\be\label{3order}
g^{(=3)}=\kappa_1 {\mathtt y}^3+
\kappa_2 {\mathtt y} {\mathtt y}_{\mathtt x}^2
+\kappa_3 {\mathtt y} {\mathtt y}_{\mathtt t}^2\,,
\qquad \kappa_1, \kappa_2,  \kappa_3 \in\mathbb R\,,
\ee
and $g^{(\geq 5)}$ collects terms of order at least five in $(\mathtt y, \mathtt y_{\mathtt x}, \mathtt y_t) $.
We   assume 
the  {\it non-degeneracy} condition
\begin{equation}\label{NRF}
\kappa_1 +  (\kappa_2  + \kappa_3) i^2 + \kappa_3 \mm  \neq 0 \, ,  \quad \forall i \in \mathcal I^+ \, . 
\end{equation}
Note that, for each $ \mm > 0 $, condition \eqref{NRF}  is verified for all the
$ (\kappa_1, \kappa_2, \kappa_3) \in \R^3 $ outside finitely many hyperplanes, for example for {\it each}
$ (\kappa_1, \kappa_2, \kappa_3) \neq 0 $ with non negative  components  $ \kappa_j  \geq  0 $, $ j =1, 2, 3 $.

Fix a compact interval  $[\mm_1,\mm_2]\subset (0,\infty) $, and  
assume that the mass $\mm\in[\mm_1,\mm_2]$ satisfies 
the  finitely many {\it non-resonance} conditions
\begin{equation}\label{potapota}
 (\l_i^{-1} \pm \l_j^{-1}) 4 {\bar \om}     \,, \
 \l_j^{-1}  4 {\bar \om}   
 \,   
  \notin \,  (2n-1)\Z^{n/2} \setminus \{ 0 \} \, ,
 \ \forall i,j \in \N \setminus {\cal I}^+ \, , i, j \leq C_0 \, , 
\end{equation}
where $ \bar\om:=(\l_h)_{h\in{\cal I}^+}$,
$\l_h=\sqrt{h^2+\mm},$ $n$ is twice the cardinality of ${\cal I}^+$,
and $C_0$ is a suitably large constant depending  on $\mm_1,\mm_2,{\cal I}^+$. 
Note that, for a given set ${\cal I}^+ $ of  tangential sites, condition \eqref{potapota}  is verified, by analiticity, for all the masses
$\mm \in [\mm_1, \mm_2 ] $ except {\it finitely} many (and independently of $\kappa_1,\kappa_2,\kappa_3$).

\begin{theorem}\label{thm:DNLW}
Assume that the tangential sites
$ {\cal I}^+  \subset \N \setminus \{ 0 \} ,$
the mass $ \mm \in [\mm_1, \mm_2 ] $ and  $\kappa_1,\kappa_2,\kappa_3 \in \R $ satisfy
\eqref{NRF}, \eqref{potapota}. Then
the DNLW equation
\eqref{DNLW} with a real analytic nonlinearity satisfying
\eqref{geven}, \eqref{REALITY},
 \eqref{gleading}-\eqref{3order}
admits small-amplitude, analytic (both in $ t $ and $ \mathtt x
$), quasi-periodic solutions \be\label{solutions} {\mathtt y}  =
{\mathop\sum}_{j \in {\cal I}^+}  \sqrt{ 8 \xi_j} \, \l_j^{-1}   \cos (
\om_j^\infty (\xi) \, t ) \cos ( j {\mathtt x} ) + o( \sqrt{\xi}
), \quad
 \om^\infty_j (\xi) \stackrel{\xi \to 0}\longrightarrow \sqrt{j^2 + \mm}
 \ee
satisfying \eqref{simmet},  for a Cantor-like set of parameters
with  density $ 1$ at $ \xi = 0 $. The
quasi-periodic solutions have {\it zero Lyapunov exponents} and
the linearized equations can be {\it reduced to constant
coefficients} (in a phase space of functions even in $ \mathtt x$).
The term $ o( \sqrt{\xi} ) $ in \eqref{solutions}
is small in some analytic norm.
\end{theorem}

Let us comment on  the hypothesis of Theorem \ref{thm:DNLW}.

\begin{enumerate}
\item {\bf Reversibility in time and space.} The 
 assumptions \eqref{geven}, \eqref{REALITY},
are  natural conditions for the existence of quasi-periodic solutions of \eqref{DNLW}, 
because they imply the reversibility assumption of Moser \cite{Mo67} 
on the subspace of functions even in $ \mathtt x $ (which does {\it not} follow by requiring only
one of them), and so they allow to solve the homological equations along the KAM proof. 
Terms like $ {\mathtt y}_{\mathtt x}^p $, $ {\mathtt y}_{t}^p $ with $ p $ odd,  destroy the oscillations 
of the Birkhoff normal form and produce drifts of the actions incompatible 
with the existence of quasi-periodic solutions. 
Proposition \ref{nonex} 
proves rigorously  these  non existence results 
using suitable Lyapunov functions, for which 
 terms like $ {\mathtt y}_{\mathtt x}^p $ and $ {\mathtt y}_t^p $ act as friction terms. 
This  shows the role of condition \eqref{REALITY}. 
As an example, the nonlinearity $ g = {\mathtt y}^3 + {\mathtt y}_{\mathtt x}^5 $ satisfies all the conditions
 \eqref{geven}, \eqref{gleading}, \eqref{3order} (and \eqref{NRF} holds for each $ {\cal I}^+ $), but 
 {\it not}  \eqref{REALITY}, and non trivial quasi periodic solutions of \eqref{DNLW} do not exist.  

Thanks to \eqref{REALITY} we can restrict to solutions which are even in $ \mathtt x $ 
and this  simplifies the KAM proof because the normal form \eqref{90210} is diagonal.
However, as said above, the main reason to assume \eqref{geven} + \eqref{REALITY} is that they 
imply the  reversibility with respect to the involution used in Moser \cite{Mo67}  (see \eqref{defS}, \eqref{pandistelle}). 
This does not 
follow, for example,   by \eqref{geven} and 
the condition $ g( - {\mathtt x}, -  {\mathtt y},  {\mathtt
y}_{\mathtt x},  {\mathtt v} ) = - g( {\mathtt x},   {\mathtt y},
{\mathtt y}_{\mathtt x}, {\mathtt v}  ) $ for which 
the subspace of functions $ ({\mathtt y}, {\mathtt v} )({\mathtt x}) $ 
odd in $ {\mathtt x} $ is invariant 
(Dirichlet boundary conditions). 
One could possibly deal also with other nonlinearities 
using the involution 
 $ ({\mathtt y}({\mathtt x}), {\mathtt v}({\mathtt x})) \mapsto ( {\mathtt y}( - {\mathtt x}), - {\mathtt v}( - {\mathtt x})) $
  which implies the Moser reversibility as well.
\item
 {\bf Mass $ \mm > 0 $.}
Also the assumption on the mass $ \mm \neq  0 $ is natural. When $ \mm = 0 $, 
Proposition \ref{FritzJ} proves that
\eqref{yttyxx} has no smooth solutions for all times except the
constants. In Proposition \ref{zeromass} we prove other
non-existence results of quasi-periodic solutions for  DNLW equations
satisfying both \eqref{geven}, \eqref{REALITY}, but  with mass $
\mm = 0 $.
\item
{\bf Twist.}
The term  $g^{(=3)} ({\mathtt y}, {\mathtt y}_{\mathtt x}, {\mathtt v}) $
in \eqref{3order} is the most general cubic nonlinearity which
satisfies
\eqref{geven}, \eqref{REALITY} and which is $\mathtt x$-independent.
Proposition \ref{nonex} proves that for
 $ {\mathtt y}_{\mathtt x}^3  $, $  {\mathtt y}^2 {\mathtt y}_{\mathtt x}  $, $ {\mathtt v}^3 $,
there exist no non-trivial quasi-periodic solutions of \eqref{DNLW}. In \eqref{gleading}  the
 leading term  $g^{(=3)}$ could also depend explicitly on $\mathtt x $
 and the higher order  nonlinearities  have order four, see Remark \ref{carmelo}. 
\item {\bf $ \bf \mathtt x $-dependence.} The nonlinearity $ g $
in \eqref{DNLW} may explicitly depend on the space variable $
\mathtt x $. 
This is a novelty  with respect to
\cite{BBP1} which  used the
conservation of momentum,
see comments below.  
\item {\bf Derivative vs quasi-linear NLW.}
Klainermann-Majda \cite{KM} exhibithed
a class of quasi-linear wave equations which do not have smooth
periodic (a fortiori quasi-periodic)  solutions except the
constants.  In this respect \cite{KM} may suggest  that Theorem \ref{thm:DNLW} is
 optimal regarding the order of (integer) derivatives in the nonlinearity.
\end{enumerate}

The proof of Theorem \ref{thm:DNLW} is based on a KAM theorem (see Theorem \ref{thm:IBKAM})
whose key step is,
like in  \cite{BBP1}, to prove
 the  first order asymptotic expansion of the perturbed normal frequencies of the linearized
equations  along the iteration, see \eqref{freco}.
This enables to verify the well known second order Melnikov conditions
which allow to reduce the KAM normal form  to  constant coefficients.
Unlike the case where $ g $ does not depend on the derivatives
$ {\mathtt y}_t $,  $ {\mathtt y}_{\mathtt x} $, this expansion requires hard work.
This is  achieved by the notion of {\it quasi-T\"oplitz
vector field} introduced in section \ref{sec:QT}.
This class is
closed with
respect to Lie brackets and
Lie transform (Propositions \ref{festa2}-\ref{main2}).
This concept is  clearly modelled on the Hamiltonian case in \cite{BBP1}, \cite{PX},
and it is related to the T\"oplitz-Lipschitz functions in Eliasson-Kuksin \cite{EK}-\cite{EK1}
(see also \cite{GT}), but there are differences. Actually this notion appears natural for vector fields.
We  underline two main novelties.

\begin{enumerate}
\item  As already said, here we consider the general case of  $ {\mathtt x} $-dependent nonlinearities which break the
translation invariance.
In \cite{BBP1}, and \cite{PX}, \cite{PP}, the theory of quasi-T\"oplitz functions was developed for
 $ {\mathtt x} $-independent
nonlinearities, namely it relied on the conservation of momentum.
This property was used in essential ways, for example in order to prove that the class of quasi-T\"oplitz functions
is closed under Poisson-bracket. 
A point of conceptual interest in this paper is that we show how to use efficiently
the notion of momentum  also when this is not a conserved quantity. Monomial vector fields with
a large momentum should be less and less relevant for dynamics.
This is efficiently implemented by the introduction of the $ \mathtt a $-momentum norm  (Definition \ref{MNV})
which penalizes
the high momentum monomials, see \eqref{torquato}.
This allows
to neglect in Proposition  \ref{festa2}  the high momentum monomial vector fields, by slightly decreasing
 the parameter $ \mathtt a $.
 With this new idea
the theory of quasi-T\"oplitz vector fields is obtained
similarly to \cite{BBP1}. 
\item Another point of conceptual interest is
to use the notion of momentum working in a subspace (here of even functions).
Until now it was not clear how to proceed, see  the end of section \ref{sec:ideaproof}.
In this paper this is achieved by the symmetrization procedure described in section \ref{sec:symme}.
The key observation is that the quasi-Toplitz norm does not increase under symmetrization, Lemma \ref{topotip}.
\end{enumerate}

We will add some more technical comments about the proof in section \ref{sec:ideaproof}.
\\[1mm]
Now we complement  Theorem \ref{thm:DNLW} with some {\it non-existence} results.
\begin{proposition}\label{nonex}
Let $ p \in \N $ be odd. The DNLW equations \eqref{DNLW}
with
\be\label{e3}
i)\ g={\mathtt
y}_{\mathtt x}^p + f({\mathtt y})\,,\qquad
ii)\ g=\partial_{\mathtt x} ( {\mathtt y}^p ) +
f({\mathtt y})\,,\qquad
iii)\ g=  {\mathtt y}^p_t +
f({\mathtt y})
 \ee have no
smooth quasi-periodic solutions except trivial periodic
solutions $\mathtt y(t,\mathtt x)=c(t)$ for $i),ii)$
and
$\mathtt y(t,\mathtt x)=c(\mathtt x)$ for $iii)$, respectively. If  $f\equiv 0$ then
$c(\cdot)\equiv const.$
\end{proposition}

\begin{pf}
The function
$ M := \int_\T {\mathtt y}_{\mathtt x}
\, {\mathtt y}_t \, d{\mathtt x} $
is a Lyapunov function of \eqref{e3}-$i$) since
$$
\frac{d}{dt} M  =
\int_\T  {\mathtt y}_{\mathtt x}^{p+1}   \, d{\mathtt x} \geq 0 \, .
$$
Hence $ M $ strictly
increases along the solutions unless $ {\mathtt
y}_{\mathtt x} (t,{\mathtt x}) = 0 $, $ \forall t $, namely  $
{\mathtt y}(t,{\mathtt x}) = c(t) $.  Case $ii)$ is similar.
A Lyapunov function of \eqref{e3}-$iii$) is
$ H := \int_\T \frac{ {\mathtt v}^2}{2}
+   \frac{{\mathtt y}_{\mathtt x}^2}{2}  - F({\mathtt y}) \,
d{\mathtt x}$ where $F'=f.$
\end{pf}

The mass term $ \mm {\mathtt y} $ could be necessary to have
existence of quasi-periodic solutions.

\begin{proposition}\label{FritzJ}
The DNLW equation \be\label{FJ} {\mathtt y}_{tt} -
{\mathtt y}_{{\mathtt x}{\mathtt x}} = {\mathtt y}_t^2 \, , \quad
{\mathtt x} \in \T \, , \ee has no smooth solutions defined for
all times except the constants.
\end{proposition}

\begin{pf}
We decompose the solution
$ {\mathtt y}(t, {\mathtt x}) = {\mathtt y}_0 (t) + \tilde {\mathtt
y}(t, {\mathtt x}) $ where $ {\mathtt y}_0 := \int_\T  {\mathtt y} (t, {\mathtt x}) d   {\mathtt x}
$ and $ \tilde {\mathtt y}  := {\mathtt y} - {\mathtt y}_0 $
has zero average in $ \mathtt x $. Then,  projecting \eqref{FJ} on
the constants, we get \be\label{zeromode} \ddot {\mathtt y}_0 =
\int_\T  {\mathtt y}_t^2 d  {\mathtt x} = \int_\T (  \dot {\mathtt
y}_0 +  {\tilde {\mathtt y}}_t)^2 d  {\mathtt x} =
 {\dot {\mathtt y}_0}^2
+  2 \dot {\mathtt y}_0 \int_\T  {\tilde {\mathtt y}}_t d {\mathtt
x} +  \int_\T  {\tilde {\mathtt y}}_t^2 d  {\mathtt x} = {\dot
{\mathtt y}_0}^2 +  \int_\T  {\tilde {\mathtt y}}_t^2 d {\mathtt
x} \geq  {\dot {\mathtt y}_0}^2 \, . \ee Hence $ {\mathtt v}_0 :=
\dot {\mathtt y}_0 $ satisfies $ \dot {\mathtt v}_0 \geq  {\mathtt
v}_0^2 $ which blows up unless $ {\mathtt v}_0 \equiv 0 $. But, in
this case, \eqref{zeromode}  implies that $ {\mathtt y}_t (t,
{\mathtt x}) \equiv 0 $, $ \forall  {\mathtt x} $. Hence $
{\mathtt y}(t,  {\mathtt x}) = {\mathtt y}({\mathtt x}) $ and
\eqref{FJ} (and $ {\mathtt x} \in \T $) imply that $ {\mathtt
y}(t,  {\mathtt x} ) = $ const.
\end{pf}

The above non-existence result may be generalized as follows:

\begin{proposition}\label{zeromass}
Let $ p , q \in \N $ be even. Then the derivative NLW equations
\be\label{e9} {\mathtt y}_{tt} - {\mathtt y}_{{\mathtt x}{\mathtt
x}} = {\mathtt y}_{\mathtt x}^{p} \, , \qquad
{\mathtt y}_{tt} - {\mathtt y}_{{\mathtt x}{\mathtt x}} = {\mathtt
y}_t^{p}  \, , \qquad
 {\mathtt y}_{tt} - {\mathtt y}_{{\mathtt
x}{\mathtt x}} = {\mathtt y}_{\mathtt x}^{p} + {\mathtt y}_t^{q}
\, , \quad {\mathtt x} \in \T \, ,
\ee
have no smooth periodic/quasi-periodic solutions except the
constants.
\end{proposition}

\begin{pf}
If there exists a
periodic solution $ ({\mathtt y}(t,{\mathtt x}), {\mathtt
v}(t,{\mathtt x})) $ of the first equation, with  period $ T $, then
$$
\int_0^T \int_{\T } ({\mathtt y}_{tt} - {\mathtt y}_{{\mathtt
x}{\mathtt x}}) dt d{\mathtt x}  = 0 = \int_0^T  \int_\T  {\mathtt
y}_{\mathtt x}^{p}(t,{\mathtt x})
 d{\mathtt x} dt \, .
$$
Hence, $ \forall t \in [0,T]  $,
$ {\mathtt y}_{\mathtt x} (t,{\mathtt x}) = 0 $,  $ \forall {\mathtt x} \in {\T } $, that is
$ {\mathtt y}(t,{\mathtt x}) = c(t) $. Inserting in
\eqref{e9} we get $ c_{tt} (t ) = 0 $
and its only periodic solutions are $  c(t) = const $. For
quasi-periodic solutions the argument is the same. The other 
equations can be treated analogously.
\end{pf}

\subsection{About the proof of Theorem \ref{thm:DNLW}}\label{sec:ideaproof}

\noindent {\bf Complex formulation.}
In  the unknowns
$$
u^+ :=\frac{1}{\sqrt 2} (D{\mathtt y}-\ii  {\mathtt v} )\,,\quad
u^- := \frac{1}{\sqrt 2} (D{\mathtt y}+\ii  {\mathtt v} )\,, \quad
D :=\sqrt{-\partial_{{\mathtt x}{\mathtt x}}+\mm} \, , \ \
\ii:=\sqrt{-1} \, ,
$$
systems \eqref{firstor}  becomes the first order system
\be\label{brace}
u^+_t = \ii D u^+ + \ii \mathtt g( u^+, u^- ) \, , \quad u^-_t =  \ \
-\ii D u^- -  \ii \mathtt g(u^+, u^-)
\ee where
\begin{equation}\label{padella}
\mathtt g(u^+, u^-) = -\frac{1}{\sqrt{2}} \, g\Big( {\mathtt x},
D^{-1} \Big(\frac{ u^+ + u^- }{ \sqrt{2}}\Big) , D^{-1}
\Big(\frac{ u^+_{\mathtt x} +  u^-_{\mathtt x}}{  \sqrt{2}}\Big) ,
\frac{u^- - u^+}{ \ii \sqrt{2}} \Big) \, .
\end{equation}
Since $ g $ is 
real on real, the subspace $ {\mathtt R} := \{ \overline{
u^+} = u^- \}  $
 is invariant under the flow evolution of \eqref{brace}.
Clearly, this corresponds to real valued solutions
$ ({\mathtt y}, {\mathtt v})  $ of 
\eqref{firstor}.
By \eqref{REALITY} the subspace of even functions
\be\label{Esub}
{\mathtt E} := \big\{ u^+ ({\mathtt x}) = u^+ (-{\mathtt x}) \, , \ u^-
({\mathtt x}) = u^- (- {\mathtt x}) \big\}
\ee
is invariant. 
Moreover \eqref{brace} is reversible with respect
to the involution
\be\label{invco} S (u^+ , u^- ) = (u^-, u^+)
\ee
which is nothing but \eqref{invol1} in the variables $ (u^+ ,  u^- ) $.
\smallskip

\noindent
 {\bf Dynamical systems formulation.} We introduce 
 coordinates by Fourier transform
\be\label{Fourier} u^+ = {\mathop\sum}_{j \in \Z} u_j^+ e^{\ii j {\mathtt
x}} \, , \quad u^- = {\mathop\sum}_{j \in \Z}  u^-_j e^{- \ii j {\mathtt x}
}\,. \ee Then \eqref{brace} becomes the infinite dimensional
dynamical system
\begin{equation}\label{tuc}
{\dot u}^+_j =  \ii  \l_j u^+_j + \ii {\mathtt g}_j^+ ( \ldots,
u^+_h, u^-_h, \ldots )\, , \quad  {\dot u}^-_j
=
- \ii  \l_j u^-_j -
\ii {\mathtt g}_{j}^-  ( \ldots, u^+_h, u^-_h, \ldots ) \, ,
\end{equation}
$ \forall j \in \Z $,  where $ \l_j := \sqrt{j^2+\mm} $  are the eigenvalues of  $ D $ and
\begin{equation}\label{piggione}
\mathtt g_j^+ = \frac{1}{2 \pi} \int_{\T}\mathtt g \Big( {\mathop\sum}_{h
\in \Z} u_h^+ e^{\ii h {\mathtt x}} , {\mathop\sum}_{h \in \Z} u^-_h e^{-
\ii h {\mathtt x}} \Big) e^{- \ii j {\mathtt x}} d{\mathtt x} \, ,
\quad \quad {\mathtt g}_j^-  := {\mathtt g}_{-j}^{+} \, .
\end{equation}
By \eqref{Fourier}, the ``real" subset $ {\mathtt R} $  reads
$ \overline{u_j^+} = u_j^- $ (this is the motivation for the choice of the signs in
\eqref{Fourier}).
The  invariant subspace $ \mathtt E $ of even functions in \eqref{Esub}
reads, under Fourier transform,
\be\label{parit} E := \big\{
u^+_{j} = u^+_{-j} \,  , \ u^-_{j} =  u^-_{-j} \, , \ \forall j
\in \Z \big\} \, .
\ee
By \eqref{Fourier} the
involution \eqref{invco} reads \be\label{Seven} S : (u^+_j, u^-_j)
\to (u^-_{-j}, u^+_{-j}) \, , \quad \forall j \in \Z \, . \ee
Finally, since $ g $ is real analytic, the assumptions
\eqref{geven} and \eqref{REALITY} imply the important property
\be\label{greco}
\mathtt g^\pm_j ( \ldots, u_i^+ , u^-_i, \ldots )
\ {\it has \ real \ Taylor \ coefficients \ in } \
(u_i^+ , u^-_i ) \, .
\ee
This property is compatible with an oscillatory behavior for
\eqref{tuc},  excluding friction phenomena.
\smallskip

\noindent {\bf Abstract KAM theorem.} For every choice of the
symmetric   \textsl{tangential sites} \be\label{Isym} \mathcal I =
\mathcal I^+\cup (-\mathcal I^+) \quad {\rm with} \quad \mathcal
I^+ \subset \N \setminus \{0\}  \, , \ \sharp \mathcal I = n \, ,
\ee
we introduce (after the Birkhoff normal form of section
\ref{sec:BNF}),
 action-angle variables
 \be\label{AAV} u_j^+ = \sqrt{\xi_{|j|} + y_j} e^{\ii
x_j} \, , \, u^-_j = \sqrt{\xi_{|j|} + y_j} e^{-\ii x_j},  \  j
\in {\cal I} \, , \  (u_j^+, u^-_j)  = (z_j^+,z_j^-)\equiv(z_j,
\bar z_j) \, , \, j \notin {\cal I} \, , \ee where $ | y_j | <
\xi_{|j|}  $.  Then \eqref{tuc} is conjugated to a
parameter dependent family of vector fields (as in section \ref{sec:4})
\begin{equation}\label{HKAM}
\XX := {\cal N} + {\cal P }
\end{equation}
with a normal form  ${\cal N} $ as in \eqref{90210}, \eqref{sio}, and a perturbation $ {\cal P } $ as in
\eqref{perturbVF} which satisfies (A1)-(A4).
In particular the vector field \eqref{HKAM} is
\begin{enumerate}
\item  {\sc reversible} (Definition \ref{reversVF}) with
respect to the involution \be\label{defS} S : (x_j, y_j, z_j, \bar
z_j) \mapsto (-x_{-j}, y_{-j} , \bar z_{-j}, z_{-j}) \, , \
\forall j \in \Z \, , \quad S^2 = I \, ,
\ee
 which is nothing but \eqref{Seven} in the variables \eqref{AAV}.
\item {\sc real-coefficients} (Definition \ref{rean}), by  \eqref{greco}.
\item {\sc even.}  The vector field $ {\cal P} : E \to E $ and so the subspace
\begin{equation}\label{pandistelle}
 E := \big\{ x_j = x_{-j}\,,\  y_j = y_{-j}\,, \ j \in \mathcal I\,,
 \  \ z_j = z_{-j}\,,\  \bar z_j = \bar z_{-j}\, , \ j \in \Z\setminus \mathcal I \big\}
\end{equation}
is invariant under the flow evolution of \eqref{HKAM}.
\item  {\sc Quasi-T\"oplitz.}  The perturbation $ \cal P $ is a {\it quasi-T\"oplitz vector field},
Definition \ref{topbis_aa}.
\end{enumerate}
The reversibility  property \eqref{defS} on the subspace  $ E $ in \eqref{pandistelle} 
implies that   the average of the term $ {\cal P}^{(y)}(x,0,0,0) $ is zero (because $ {\cal P}^{(y)}(x,0,0,0) $ is an 
odd function in $ x $)
 along  the whole iteration, otherwise  quasi-periodic solutions would not exist. Note that we use 
 both \eqref{geven} and  \eqref{REALITY} for the solvability of the homological equations in Lemma \ref{lem:HOMO}. 
Then the ``real-coefficients" property implies that the corrections to the normal form are purely imaginary
(elliptic).
\\[1mm]
\noindent
{\bf 
Quasi-T\"oplitz property.}
The second order Melnikov non resonance conditions are verified
proving
that the elliptic frequencies
(after the application of the KAM Theorem \ref{thm:IBKAM})
satisfy an asymptotic expansion like $\O_j^{\infty}(\xi)=|j|+c(\xi)+O(1/|j|)$, see \eqref{freco} and \eqref{carbonara}.
Indeed, since $c(\xi)$ is {\sl independent} of $j$,
it cancels in the difference $\O_j^{\infty}(\xi)-\O_i^{\infty}(\xi)$
and the measure estimates
follow as in the semilinear case (see \cite{Po3}), where
$c(\xi)\equiv 0 $. We only state them in Theorem \ref{thm:measure}, whose proof is like in \cite{BBP1}.

The KAM corrections to the frequencies are  the coefficients of
 the 
 linear monomial vector fields
 $ z_j \partial_{z_j}$, $ {\bar z}_j \partial_{ {\bar z}_j}$ of the perturbation $P$, and
we want to show 
 that, for  $ |j| > N $, they  assume a constant value up to an error of $ O(N^{-1}) $.
Since we need to work with a class of vector fields fulfilling
the Lie algebra property,
  we can not clearly impose conditions only on these diagonal terms,
but we have to consider a larger set of  vector fields
which are only {\sl approximately} $x$-independent,
  linear and diagonal (and $y$-independent).
The quasi-T\"oplitz vector fields introduced in section \ref{sec:QT} fulfill
quantitively these requirements, see
comments above Definition  \ref{def:vfm}.

\smallskip

\noindent
{\bf The symmetrization procedure.}
In the subspace of functions even in $ \mathtt x $ the notion of {\sc momentum} of a monomial is not well defined.
For example the vector fields  $ z_{-j} \partial_{z_i} $ and $ z_j \partial_{z_i} $,  that
have {\sc different} momentum, 
 are identified.
In other words we can {\sc not} 
work directly in the cosine basis   $ \{ \cos (j \mathtt x) \}_{j \geq 0} $, which would be 
natural looking for solutions even in $ \mathtt x $ (avoiding the double eigenvalues). 

Then we proceed as follows. In  system \eqref{HKAM} we think $ x_j, y_j,  z_{j}^\pm $,
as independent variables.
In this case, since the linear frequencies $ \o_{-j} = \o_j  $ , $ \OO_{-j} = \OO_j $ are resonant, along the KAM iteration,
the monomial vector fields of the perturbation
$$
e^{\ii k\cdot x}\partial_{x_j},  e^{\ii k\cdot x}y^i\partial_{y_j},   k \in \Z^n_{\rm odd},  |i| = 0,1, \,  j \in \mathcal I \, ,
\ e^{\ii k\cdot x} z_{\pm j} \partial_{z_j},  e^{\ii k\cdot x}  \bar z_{\pm j} \partial_{\bar z_j},
\forall  k \in \Z^n_{\rm odd}\,,\  j\in\Z\setminus \mathcal I\, ,
$$
where
\be\label{Zodd}
k \in \Z^n_{\rm odd}  := \big\{ k \in \Z^n \ : \  k_{-j} = - k_j \, , \ \forall j \in\mathcal I  \big\}
\ee
can not be averaged out. On the other hand, on the
invariant subspace $ E $, where we look for the quasi-periodic solutions, the above terms can be replaced
by the constant coefficients monomial vector fields, obtained setting $ x_{-j} = x_j $, $ z^\pm_{-j} = z^\pm_j $.
Replacing
the vector field $ \cal P $ with its symmetrized  $ {\cal S} {\cal P} $ (Definition \ref{SVF})
the  $ \aaa $-momentum and quasi-T\"oplitz norms do not increase (Proposition \ref{topotip}).
Both $ \cal P $ and $ {\cal S} {\cal P} $ determine the
same dynamics on the subspace $ E $ (Proposition \ref{carbonara1}).
The vector field $ {\cal S}{\cal P } $ is symmetric and reversible as well (see \eqref{zenone})
and the homological equations  \eqref{tebe}
can  be solved, see Lemma 
\ref{lem:HOMO}.
This procedure allows the KAM iteration to be carried out.
Remark \ref{simme} shows that
the symmetrization procedure is required at each KAM step.

\smallskip

In section \ref{sec:BNF} we finally apply the abstract KAM Theorem \ref{thm:IBKAM} to prove
Theorem \ref{thm:DNLW}. The main steps are the proof that the vector field of $ g $
is quasi-T\"oplitz (Lemma \ref{aragorn}),
that the Birkhoff normal form transformation preserves the quasi-T\"oplitz property  (Proposition \ref{BNF})
and that the frequency-to-action map is twist, see \eqref{platea}, \eqref{plateabis}. 
\\[1mm]
\noindent{\bf  Acknowledgments:}  
{\it We thank L. Corsi and an anonymous referee for many useful suggestions}.

\section{Vector fields formalism}\setcounter{equation}{0}

We introduce the main properties of
the vector fields used along the paper (commutators, momentum, norms,  reversibility, degree, ...).
We shall refer often to section 2 of \cite{BBP1}. 
The first difference with respect to  \cite{BBP1} 
is that we have to work  at the level of vector fields and not of functions (Hamiltonians).

\smallskip

For a finite set $ \cal I \subset \Z $
(possibly empty) and  $ a \geq 0, p > 1/2 $,
we define the Hilbert space
\be\label{ellcal}
\ell^{a,p}_{\cal I} := \Big\{  z=\{ z_j\}_{j\in\Z\setminus {\cal I}} \, , \ z_j \in \C   
\ : \
\| z \|_{a,p}^2 := {\mathop\sum}_{j\in\Z\setminus {\cal I}} |z_j|^2 e^{2a|j|} \langle j
\rangle^{2p} < \infty  \Big\}
\ee
that, when $ \cal I = \emptyset $, we denote more simply by $ \ell^{a,p} $.
Let $ n $ be the cardinality of $ \cal I $. We consider
$ V := \C^n \times \C^n \times \ell^{a, p}_{\cal I} \times \ell^{a, p}_{\cal I} $
(denoted  by $ E $ in \cite{BBP1}) with  $(s,r)$-weighted norm
\begin{equation}\label{normaEsr}
v =  (x,y,z,\bar z) \in V \, , \quad
\|v\|_{s,r} = \frac{|x|_\infty}{s} + \frac{|y|_1}{r^2}
 +\frac{\|z\|_{a,p}}{r}+\frac{\|\bar z\|_{a,p}}{r}
\end{equation}
where $ 0 < s, r < 1 $, and
$ |x|_\infty := \max_{h =1, \ldots, n} |x_h| $,
$ |y|_1 := {\mathop\sum}_{h=1}^n |y_h| $.

Note that $z$ and $\bar z$ are \textsl{independent} variables.
We shall also use the notation
$ z_j^+ = z_j $, $ z_j^- = \bar z_j $,
and
\be\label{setV}
\mathtt{V} := \big\{ x_1,\ldots, x_n, y_1,\ldots, y_n, \ldots, z_j, \ldots, \bar z_j,\ldots  \big\} \,,
\quad j\in \Z\setminus \mathcal{I} \, .
\ee
As phase space, we  consider the toroidal domain
\be\label{Dsr}
D (s,r) := \T^n_s \times D(r) := \T^n_s \times B_{r^2} \times B_r
\times  B_r \subset V
\ee
where
$ \T^n_s := \big\{ x \in \C^n \, : \, {\rm Re}(x_h) \in \T^n := 2\pi \R^n/\Z^n \, , \,
\max_{h=1, \ldots, n} |{\rm Im} \, x_h | < s
 \big\} $, $ B_{r^2}  := \big\{ y \in \C^n \, : \, |y |_1 < r^2 \big \} $
and $ B_{r} \subset \ell^{a,p}_{\cal I} $ is the open ball of radius $ r $
centered at zero.
If $ n = 0 $  then $ D(s,r)
 \equiv B_r\times B_r \subset
\ell^{a,p}\times \ell^{a,p}$.

We also introduce the ``real" phase space
\begin{equation}\label{realps}
    \mathtt{R}(s,r):= \big\{ \, v = (x,y,z^+ , z^-)\in D(s,r)\ \ {\rm :} \ \
    x \in \T^n\,,\ y\in \R^n\,, \  \overline{z^+} = z^-  \big\}
\end{equation}
where $ \overline{z^+}  $ is the  complex conjugate of $ z^+ $.

We consider vector fields of the form
\be\label{4comp}
X(v) = (X^{(x)}(v),X^{(y)}(v),X^{(z)}(v),X^{(\bar z)}(v)) \in V
\ee
where $ v \in D(s,r) $ and
$ X^{(x)}(v),X^{(y)}(v)\in \C^n $, $ X^{(z)}(v),X^{(\bar z)}(v)\in \ell^{a,p}_{\mathcal I} $.
We also use the differential geometry notation
\be\label{XV}
X (v)= X^{(x)} \partial_x + X^{(y)} \partial_y +
X^{(z)} \partial_z  + X^{(\bar z)} \partial_{\bar z}
= {\mathop\sum}_{\vgot\in \mathtt{V}} X^{(\vgot)}\partial_{\vgot} \, ,
\ee
recall \eqref{setV}.
Equivalently we write $ X(v) =  \big( X^{(\vgot)}(v) \big)_{\vgot \in \mathtt{V}} $ where each
component is a {\it formal}  scalar power series
\begin{equation}\label{formalpowerA}
X^{(\vgot)}(v) = \sum_{(k,i,\a,\b) \in {\mathbb I}}
X^{(\vgot)}_{k, i, \alpha, \b} \, e^{\ii k \cdot x} y^i  z^\alpha \bar z^{\b}
\end{equation}
with coefficients $ X^{(\vgot)}_{k, i, \alpha, \b} \in \C $ and multi-indices in
\be\label{I}
\mathbb{I}:=\Z^n \times \N^n \times \N^{(\Zi )} \times \N^{(\Zi )}
\ee
where
$ \N^{(\Z \setminus {\cal I})}  := \big\{ \a := (\a_j)_{j \in \Z \setminus {\cal I}} \in \N^\Z \  {\rm with}  \
|\a| := {\mathop\sum}_{j \in \Z \setminus {\cal I}} \a_j < + \infty \big\} $.
In \eqref{formalpowerA} we use the standard multi-indices notation
$ z^\a {\bar z}^\b := \Pi_{j \in \Z \setminus {\cal I}} \, z_j^{\a_j}  {\bar z}_j^{\b_j} $.

The formal vector field $ X $  is absolutely convergent
in $ V  $ (with norm \eqref{normaEsr}) at $ v \in D(s,r) $
if every component $X^{(\vgot)}(v) $, $ \vgot \in \mathtt{V}$, is
absolutely convergent
and
$ \big\| \big( X^{(\vgot)}(v)\big)_{\vgot \in \mathtt{V}} \big\|_{s,r} < +
\infty $.

\begin{definition}\label{mononucleosi} {\bf (monomial vector field)}
A monomial vector field is
\be\label{commodo}
\mathfrak m_{k,i,\a,\b;\vgot'} (v) =
 {\mathfrak m}_{k,i,\a,\b}(v)  \partial_{\vgot'} \qquad
 { where} \qquad
{\mathfrak m}_{k,i,\a,\b}(v) :=
e^{\ii k\cdot x} y^i z^\a {\bar z}^{\b}
\ee
is a scalar monomial.
\end{definition}

A vector field $ X $ 
 may be decomposed
as a formal  series of vector field monomials 
\begin{equation}\label{cesare}
X(v)=
 \sum_{\vgot \in \mathtt{V}}
\sum_{(k,i,\a,\b)\in\mathbb{I}} X^{(\vgot)}_{k,i,\a,\b}
e^{\ii k\cdot x} y^i z^\a \bar z^\b \partial_{\vgot} \, .
\end{equation}
For a subset of
indices $ I \subset \mathbb{I} \times {\mathtt V} $ we define the projection
\be\label{PROJE}
(\Pi_I X)(v) := \sum_{ (k,i,\a, \b, \vgot) \in I } X_{k,i,\a,
\b}^{(\vgot)} \, e^{\ii k \cdot x} y^i z^\a {\bar z}^{\b} \partial_\vgot \, .
\ee
The commutator (or Lie bracket) of two vector fields  is
$ [X,Y](v):=d X(v) [Y(v)]- d Y(v)[X(v)]  $
namely, its   $\vgot $-component is
\be\label{commco}
[X,Y]^{(\vgot)} = {\mathop\sum}_{\vgot'\in \mathtt{V}}
\partial_{\vgot'} X^{(\vgot)} Y^{(\vgot')} -
 \partial_{\vgot'} Y^{(\vgot)} X^{(\vgot')}\,.
\ee
Given a vector field $ X $, its transformed field
under the time $1$ flow generated by  $ Y $ is
\be\label{trasfoXY}
e^{{\rm ad}_Y} X = {\mathop\sum}_{k \geq 0}\frac{1}{k!} {\rm ad}_{Y}^k X \, ,
\qquad
{\rm ad}_{Y} X:=[X,Y]\,,
\ee
where ${\rm ad}_{Y}^k:={\rm ad}_{Y}^{k-1}{\rm ad}_{Y}$ and
${\rm ad}_{Y}^0:= {\rm Id} $.

\subsection{Momentum majorant norm}\label{sec:wmn}

Fix a set of indices
\begin{equation}\label{cC}
 {\cal I} := \{\pluto_1,\ldots, \pluto_n \} \subset \Z \,.
 \end{equation}

\begin{definition}\label{scalarmonomial}
The  {\sc momentum} of the vector field monomial $  \mathfrak m_{k,i,\a,\b; \vgot} $  is
\begin{equation}\label{momp}
\pi(k,\a,\b;\vgot) :=
\left\{ \begin{array}{ll} \pi(k,\a,\b) &
\ \ \ \ {\rm if} \quad \vgot \in \{x_1,\ldots,x_n,y_1,\ldots,y_n\} \\
 \pi(k,\a,\b) -\s j & \ \  \ \ {\rm if} \quad \vgot = z_j^\s \, , \ \s = \pm \, , 
 \end{array}\right.
\end{equation}
where
\begin{equation}\label{hurra}
\pi(k,\a,\b):=
{\mathop\sum}_{i=1}^n \mathtt j_i k_i + {\mathop\sum}_{j \in \Z \setminus {\cal I}} (\a_j-\b_j)j
\end{equation}
is the momentum of the scalar monomial  $ {\mathfrak m}_{k,i,\a,\b}(v)$.
\end{definition}
We say that a vector field $ X $  satisfies momentum conservation  if and only if
it is a linear combination of monomial vector fields with zero momentum.

\smallskip

Let $ \aaa\geq 0 $. 
Given a vector field $X$ as in \eqref{cesare}
we  define its ``$ \aaa $-momentum majorant" vector field
\begin{equation}\label{claudio}
(M_\aaa X)(v) :=
\sum_{\vgot \in \mathtt{V}}
\sum_{(k,i,\a,\b)\in\mathbb{I}}
e^{\aaa|\pi(k,\a,\b;\vgot)|}
|X^{(\vgot)}_{k,i,\a,\b}|
e^{\ii k\cdot x} y^i z^\a \bar z^\b \partial_{\vgot}
\end{equation}
where $\pi(k,\a,\b;\vgot)$ is the momentum of the monomial $ \mathfrak m_{k,i,\a,\b; \vgot} $ defined in \eqref{momp}.
When $ \aaa = 0 $ we simply write $ M X $ instead of $ M_0 X $,
which coincides with the majorant vector field
 in \cite{BBP1}-section 2.1.2.

\begin{definition} {\bf ($\aaa$-momentum majorant-norm)}\label{MNV}
The $ \aaa $-momentum majorant norm of a
 formal  vector field $ X $ as in \eqref{cesare} is
\begin{eqnarray}\label{normadueA}
\norma X \norma_{s,r,\aaa} & := &
\sup_{(y,z, \bar z) \in D(r)}
\Big\| \Big(
 \sum_{k,i,\a, \b}  e^{\aaa|\pi(k,\a,\b;\vgot)|} | X_{k,i,\a, \b}^{(\vgot)}| e^{|k|s}
 |y^i| |z^\a| |{\bar z}^{\b}|
\Big)_{\vgot \in \mathtt{V}} \Big\|_{s,r} 
\end{eqnarray}
where $ |k| := |k|_1 = |k_1| + \ldots + | k_n | $.
For a function $ f : D(s,r) \to \C $ it reduces to
$ \| f \|_{s,r, \aaa} :=
\sup_{D(r)} \sum_{k,i,\a,\b}e^{\aaa |\pi(\a,\b,k)|}
| f_{k,i,\a,\b}| e^{s |k|} |y^i|  |z^\a|  |\bar z^\b|  $.
\end{definition}

When $ \aaa=0 $ the norm $\|\cdot\|_{s,r,0}$ coincides with the ``majorant norm''
 introduced in \cite{BBP1}-Definition 2.6 (where it was simply denoted by
 $\|\cdot\|_{s,r}$).
By   \eqref{normadueA} and  \eqref{claudio} we get
$  \norma X \norma_{s,r,\aaa}=\norma M_\aaa X \norma_{s,r,0} $.

\begin{remark}\label{maniavanti}
By  the above relation, the norm $\| \cdot  \|_{s,r,\mathtt a} $
satisfies the same properties of the majorant norm $\| \cdot  \|_{s,r,0}  $ and
the next lemmas for the norm $\| \cdot  \|_{s,r,\mathtt a} $
follow by the analogous lemmas in
\cite{BBP1} for  $\| \cdot  \|_{s,r,0}  $.
\end{remark}

Let  $ |X|_{s,r}:=\sup_{v \in D(s,r)} \| X(v) \|_{s,r} $. Arguing as for Lemma 2.11 in \cite{BBP1} we get

\begin{lemma}\label{XMvectA}
Assume that for some $ s, r > 0 $, $ \aaa \geq 0 $, the $ \mathtt a$-momentum majorant-norm
$ \norma  X \norma_{s,r,\aaa} < + \infty $.
Then the series in \eqref{cesare}, resp.
\eqref{claudio},   absolutely converge to the {\rm analytic}
vector field $X(v)$, resp. $M_\aaa X(v)$, for every $v\in D(s,r).$
Moreover $ |X|_{s,r} $, $ |M_\aaa X|_{s,r}  \leq  \norma X \norma_{s,r,\aaa} $.
\end{lemma}

For a  vector field $ X: D(s,r) \times \mathcal O\to V $
depending on  parameters  $ \xi \in \mathcal O \subset \R^n $,
we define the $\l$-Lipschitz (momentum majorant) norm ($ \l \geq  0 $)
 \begin{eqnarray}\label{filadelfia}
 \| X\|^{\l}_{s, r,\aaa, {\cal O}} :=
 \| X \|^{\l}_{s, r,\aaa} & := &
 \| X \|_{s,r,\aaa,\mathcal O} +
 \l \| X \|^{{\rm lip}}_{s,r,\aaa,\mathcal O} \\
& := &
 \sup_{\xi \in {\cal O}} \| X (\xi) \|_{s,r,\aaa} +
 \l \sup_{\xi,\eta\in {\cal O},\ \xi\neq \eta}
\frac{\| X(\xi) - X (\eta)\|_{s,r,\aaa}}{|\xi-\eta|}
\nonumber
\end{eqnarray}
and we set
$$
 {\cal V}_{s,r,\aaa}^\l
 := {\cal V}_{s,r,\aaa,\mathcal O}^\l
 := \big\{ X : D(s,r) \times {\mathcal O} \to V   \  : \  \| X \|^{\l}_{s, r,\aaa}<\infty\ \big\} \,.
$$
Similarly, we denote by 
$ {\cal V}_{s,r,\aaa} $ 
the linear  space of vector fields  with $ \| X \|_{s, r,\aaa} < \infty  $. 
Note that, if $X$ is independent of $ \xi $, then $ \|X\|_{s,r,\aaa}^\l = \|X\|_{s,r,\aaa} $, $ \forall \l $.

It is immediate to check that the $\| \cdot \|^{\l}_{s, r,\aaa}$ norm behaves well under projections
 \eqref{PROJE}:

\begin{lemma}  {\bf (Projection)} \label{lem:pro}  $ \forall I \subset  \mathbb I \times \mathtt V $ we have
$ \norma \Pi_I X \norma_{s,r,\aaa}
 \leq \norma X
\norma_{s,r,\aaa} $
and
$ \norma \Pi_I X
\norma_{s,r,\aaa}^{{\rm lip}} \leq \norma X
\norma_{s,r,\aaa}^{{\rm lip}} . $
\end{lemma}

Important particular cases are the  ``ultraviolet" projection
\begin{equation}\label{minas}
   (\Pi_{|k| \geq K} X)(v) := \sum_{ |k| \geq K,i,\a, \b}
X_{k,i,\a, \b}^{(\vgot)} \,
e^{\ii k \cdot x} y^i z^\a {\bar z}^{\b}  \partial_\vgot\,,
\qquad
\Pi_{|k| < K}:= {\rm Id} - \Pi_{|k| \geq K}
\end{equation}
and the ``high momentum'' projection
\begin{equation}\label{tirith}
   (\Pi_{|\pi| \geq K} X)(v):=\sum_{|\pi(k,\a,\b;\vgot)| \geq K}
X_{k,i,\a, \b}^{(\vgot)} \, e^{\ii k \cdot x} y^i z^\a {\bar z}^{\b}
\partial_\vgot\,, \qquad \Pi_{|\pi| < K}:= {\rm Id} - \Pi_{|\pi| \geq K}\,.
\end{equation}

By \eqref{normadueA} the  following smoothing estimates follow:
\begin{lemma} \label{smoothL}
{\bf (Smoothing)} $ \forall K \geq 1 $
and $\l \geq 0$
\begin{eqnarray}\label{smoothl}
\norma\Pi_{|k| \geq K}X \norma_{s',r,\aaa}^\l &\leq&
 \frac{s}{s'} \,
e^{-K(s-s')}\norma X \norma_{s,r,\aaa}^\l  \,, \qquad \forall\, 0
< s' < s
\\
\label{torquato}
 \norma \Pi_{|\pi| \geq K} X \norma_{s,r,\aaa'}^\l
&\leq&
  e^{-K(\aaa-\aaa')} \norma X\norma_{s,r,\aaa}^\l  \,,\qquad\quad\,
\forall\, 0\leq \aaa'\leq \aaa\,.
\end{eqnarray}
\end{lemma}

The space of analytic vector fields
with finite $\aaa$-momentum majorant norm form a  Lie algebra.

\begin{proposition}\label{settimiosevero} {\bf (Commutator)}
Let $ X, Y \in {\cal V}_{s,r, \mathtt a}^\l $.
Then, for $\l \geq 0,$ $ r/2 \leq r' < r $, $ s/2 \leq s' < s $,
\begin{equation}\label{diffusivumsui}
\norma [ X, Y] \norma_{s',r',\aaa}^\l \leq 2^{2n+3} \delta^{-1} \norma
X\norma_{s,r,\aaa}^\l \norma Y\norma_{s,r,\aaa}^\l \, \quad
where \quad
  \d :=  \min\Big\{ 1- \frac{s'}{s},   1-  \frac{r'}{r}  \Big\}\, .
\ee
\end{proposition}

\begin{pf}
We say that a vector field $ X $ has momentum $ \pi (X) = h $
if it is an absolutely convergent series of monomial vector fields of  momentum $  h $.
It results that, if $ X , Y $ have momentum $ \pi(X) $, $ \pi (Y)  $, respectively, then
$ \pi ([X,Y]) = \pi (X) + \pi (Y) $. Then the proof of
$ \norma [ X, Y] \norma_{s',r',\aaa} \leq 2^{2n+3} \delta^{-1} \norma
X\norma_{s,r,\aaa} \norma Y\norma_{s,r,\aaa} $
follows as in \cite{BBP1}, Lemma 2.15.
The Lipschitz estimate follows as usual.
\end{pf}

%

\subsection{Degree decomposition}\label{sec:rem}

The degree of  the monomial vector field
 $ \mathfrak m_{k,i,\a,\b;\vgot} $ is defined as
$$
d(\mathfrak m_{k,i,\a,\b;\vgot}) := |i| +|\a|+|\b| - d(\vgot) \quad
{\rm where} \quad
d(\vgot) :=
\left\{ \begin{array}{ll} 0 &
\ \ {\rm if} \quad \vgot \in \{x_1,\ldots,x_n\} \\
1 & \ \ {\rm otherwise,}
 \end{array}\right.
$$
in particular $ d (\partial_x) = 0 $, $ d (\partial_y) = $ $ d (\partial_{z_j} ) = $ $ d (\partial_{ {\bar z}_j}) = -1 $.
This notion naturally extends to any  vector field 
 by monomial decomposition:
we say that a vector field has degree $ h $ if it is an absolutely convergent series of monomial vector fields of degree $h$.

The degree $d$ gives to the vector fields the structure of a graded Lie algebra:
given two vector fields $X,Y$ of degree respectively  $ d(X) $ and $ d(Y ) $, then
\be\label{sumd}
d([X,Y])= d(X)+d(Y) \, .
\ee
For a vector field $ X $ as in \eqref{cesare} we define the homogeneous component
of degree $  l \in \N $,
\be\label{phom}
X^{(l)} := \Pi^{(l)} X :=  \sum_{  |i|+|\a|+|\b|- d(\mathtt v)=l}
 X_{k,i,\a,\b}^{(\vgot)} \,
e^{\ii k\cdot x} y^i z^\a \bar z^\b \partial_{\mathtt v}
\ee
and we set
\be\label{Pleq2}
 X^{\leq 0} :=   X^{(-1)}+ X^{(0)}
 \, .
\ee
\begin{definition}\label{defVF}
We denote by $ {\mathcal R}^{\leq 0} $ the vector fields
with
degree $ \leq 0 $. Using the compact notation
$ \ugot := (y,z,\bar z) = (y, z^+, z^- ) $,
a vector field in $ {\cal R}^{\leq 0} $ writes
\be\label{ndg}
R =R^{\leq 0}=
R^{(-1)}+R^{(0)} \, , \ \quad
R^{(-1)}=R^{\ugot}(x)\partial_\ugot  \,, \quad
R^{(0)}=R^x(x)\partial_x + R^{\ugot, \ugot}(x) \ugot \,\partial_\ugot\,,
\ee
where
$ R^x (x)\in \C^n $,
$ R^\ugot \in \C^n\times \ell_\mathcal I^{a,p} \times \ell_\mathcal I^{a,p} $,
$ R^{\ugot, \ugot}(x) \in {\cal L}(\C^n\times \ell_\mathcal I^{a,p} \times \ell_\mathcal I^{a,p}) $.
In more extended notation
\begin{eqnarray}\label{baronerosso}
R^\ugot(x)\partial_\ugot
& = &
R^y(x)\partial_y + R^z (x)\partial_z+R^{\bar z} (x)\partial_{\bar z}
\nonumber
\\
R^{\ugot,\ugot}(x)\ugot\partial_\ugot
&=&
\big(R^{y,y}(x)y + R^{y,z} (x)z+R^{y,\bar z} (x)\bar z\big)\partial_y
+
\big(R^{z,y}(x)y + R^{z,z} (x)z+R^{z,\bar z} (x)\bar z\big)\partial_z
\nonumber
\\
& &+ \big(R^{\bar z,y}(x)y + R^{\bar z,z} (x)z+R^{\bar z,\bar z} (x)\bar z\big)
\partial_{\bar z} \,.
\end{eqnarray}
\end{definition}
The terms of the vector field that we want to eliminate (or normalize) along the KAM iteration are
those in $ {\cal R}^{\leq 0} $.
The graded Lie algebra property \eqref{sumd} implies that $ { \mathcal R}^{\leq 0}$ is closed by Lie bracket:
\begin{lemma}\label{lem:closed}
If $ X, Y \in {\mathcal R}^{\leq 0} $ then $ [X,Y] \in {\mathcal R}^{\leq 0} $.
\end{lemma}

\subsection{Reversible, real-coefficients, real-on-real, even, vector fields}\label{sec:rev}

We first define the class of reversible/anti-reversible vector
fields (this concept was efficiently 
 used in \cite{BDG} for finding Birkhoff-Lewis periodic solutions of NLW). 

\begin{definition}\label{reversVF} {\bf (Reversibility)}
A vector field $ X $ 
as in \eqref{4comp} is {\sc reversible}
with respect to an involution $ S $ (namely $ S^2 = I $) if
$ X \circ S = - S \circ X $.
A vector field $Y$ is {\sc anti-reversible} if
$ Y \circ S =  S \circ Y $.
\end{definition}

When the set $ {\cal I} $ is symmetric as in \eqref{Isym} and
$ S $ is the involution in \eqref{defS}, a  vector field $ X $ is reversible  
if its
coefficients (see \eqref{formalpowerA}) satisfy
\begin{equation}\label{tango}
   X^{(\mathtt v)}_{k,i,\a,\b}=
    \left\{ \begin{array}{ll}
    X^{(\hat{\mathtt v})}_{-\hat k,\hat \imath,\hat \b,\hat \a} &
\ \ \ \ {\rm if} \quad \vgot =x_j\,, \ \ j\in\mathcal I \, ,
\\
- X^{(\hat{\mathtt v})}_{-\hat k,\hat \imath,\hat \b,\hat \a} &
\ \ \ \ {\rm if} \quad \vgot =y_j\,, \ \ j\in\mathcal I \, ,
\\
  - X^{(z^{-\s}_{-j})}_{-\hat k,\hat \imath,\hat \b,\hat \a}
   & \ \  \ \ {\rm if} \quad {\mathtt v} = z_j^\s \,,  \  \ j \in \Z\setminus \mathcal I
 \end{array}\right.
\end{equation}
where
\be\label{def:hat}
\hat k := ( k_{-j} )_{j \in {\cal I}} \, , \ \hat \imath := ( i_{-j} )_{j \in {\cal I}} \, , \
\hat \b := ( \b_{-j} )_{j \in \Z \setminus {\cal I}} \, , \ \hat \a := ( \a_{-j} )_{j \in \Z \setminus {\cal I}} \, , \
\hat {\mathtt v} := ({\mathtt v}_{-j})_{j \in \Z} \, .
\ee
\begin{definition} 
 \label{rean}
A vector field
$ X = X^{(x)} \partial_x + X^{(y)} \partial_y + X^{(z^+)} \partial_{z^+}  + X^{(z^-)} \partial_{z^-} $
is
\begin{itemize}
\item
 ``{\sc real-coefficients}"
if the Taylor-Fourier coefficients of
$ X^{(x)},\ii X^{(y)}, \ii X^{(z^+)}, \ii X^{(z^-)}  $
are {\it real},
\item
 ``{\sc anti-real-coefficients}" if $ \ii X $ is real-coefficients,
\item ``{\sc real-on-real}" if
$$
X^{(x)}(v)= \overline{X^{(x)}(v)} \, , \
 X^{(y)}(v)=  \overline{X^{(y)}(v)} \, ,
\  X^{(z^-)}(v) =  \overline{X^{(z^+)}(v)} \, ,
\quad \forall v\in \mathtt R(s,r)\,,
$$
where $\mathtt R(s,r)$ is defined in \eqref{realps},
\item
``{\sc even}" if $ X : E \to E  $ (see \eqref{pandistelle}).
 \end{itemize}
\end{definition}
On the coefficients in \eqref{formalpowerA} the  {\sc real-on-real} condition amounts to
\begin{equation}\label{reality2}
    \overline{X^{(\mathtt v)}_{k,i,\a,\b}}=
    \left\{ \begin{array}{ll}
    X^{(\mathtt v)}_{-k,i,\b,\a} &
\ \ \ \ {\rm if} \quad \vgot \in \{x_1,\ldots,x_n,y_1,\ldots,y_n\} \\
  X^{(z^{-\s}_j)}_{-k,i,\b,\a} & \ \  \ \ {\rm if} \quad \vgot = z_j^\s \, ,
 \end{array}\right.
\end{equation}
and 
the  {\sc reversibility in space} condition to
\be\label{parity}
X^{(\mathtt v)}_{k,i,\a,\b} = X^{(\hat{\mathtt v})}_{\hat k, \hat \imath, \hat \a, \hat \b}
\quad  {\rm (see \ \eqref{def:hat})} \, .
\ee
\begin{definition}\label{pajata}
We denote by
\begin{itemize}
\item
 $ \mathcal R_{rev} $
the   vector fields  which are reversible, 
 real-coefficients,
real-on-real and even. 
\item
$ \mathcal R_{a\mbox{-}rev}  $
the  vector fields which are {\it anti}-reversible,
{\it anti}-real-coefficients, real-on-real and even.
\item  $ {\mathcal R}_{rev}^{\leq 0} :=
{\mathcal R}_{rev} \cap {\cal R}^{\leq 0}  $ and
 $ {\mathcal R}_{a-rev}^{\leq 0} :=
{\mathcal R}_{a-rev} \cap {\cal R}^{\leq 0}  $.
\end{itemize}
\end{definition}
If the vector field $ X $ is reversible and  $ Y $ is anti-reversible
then $  [X,Y]  $ and $ e^{{\rm ad}_Y} X  $ (recall \eqref{trasfoXY})
are reversible. 
If $ X $, resp. $ Y $, is real-coefficients, resp. anti-real-coefficients,
 then $ [X,Y]  $, $ e^{{\rm ad}_Y} X $ are real-coefficients.
If $ X,  Y $ are real-on-real, then $ [X,Y]  $, $ e^{{\rm ad}_Y} X $ are real-on-real.
If $ X,  Y $ are even then $ [X,Y]  $, $ e^{{\rm ad}_Y} X $ are even.
Therefore we get
\begin{lemma}\label{embeh2}
If $ X \in \mathcal R_{rev}$ and
 $ Y \in \mathcal R_{a-rev}$
 then $ [X,Y]  $, $ e^{{\rm ad}_Y} X \in \mathcal R_{rev}. $
\end{lemma}

By \eqref{phom}, \eqref{Pleq2} and \eqref{parity}
we immediately get (the space $E$ was defined in \eqref{pandistelle})
\begin{equation}\label{dicastro}
X_{|E}\equiv 0 \qquad
\Longrightarrow
\qquad
(X^{\leq 0})_{|E}\equiv 0\, .
\end{equation}
\begin{lemma}\label{montalto}
If  $X_{|E}\equiv 0$ and $Y$ is even 
 then $ \big([X,Y]\big)_{|E}\equiv 0 $,  $ (e^{{\rm ad}_Y}X)_{|E}\equiv 0 $.
\end{lemma}

\section{Quasi-T\"oplitz vector fields}\label{sec:QT}\setcounter{equation}{0}


Let  $ N_0 \in \N $, $ \teta, \mu \in \R $ be parameters  such
that
\begin{equation}\label{caracalla}
\cc < \teta, \mu <  \CC \, , \quad  \dueCC N_0^{L-1}+2 \kappa
N_0^{b-1} < \cc\,, \quad \kappa:=\max_{1\leq l \leq n} |\pluto_l |
\quad  (\kappa := 0 \ {\rm if} \ {\cal I} := \emptyset ) \, ,
\end{equation}
where $ {\cal I} := \{\pluto_1,\ldots, \pluto_n \} $, see \eqref{cC},  and  with the three scales
\begin{equation}\label{figaro}
   0< b< L<1  \, ,
\end{equation}
see comments before Definition \ref{def:vfm}.
In the following we will always take
$ N\geq N_0 $.

\begin{definition}\label{LM}
A scalar monomial  $ \mathfrak m(k,i,\a,\b)=e^{\ii k\cdot x} y^i z^\alpha {\bar z}^{\b} $ is
{\bf $(N,\mu)$-low momentum} if
\begin{equation}\label{zerobis}
|k| < N^b \, ,  \quad
\a+\b=\g \quad {\rm with}\quad
{\mathop\sum}_{l \in\Z\setminus{\cal I}} |l| \g_l < \mu N^L
  \, .
\end{equation}
An $(N,\mu)$-low momentum scalar monomial is {\bf $(N,\mu,h)$-low} if
 \be\label{vincula}
 |\pi(k,\alpha,\beta) - h| <N^b \,.
 \ee
We denote by ${\cal A}^L_{s,r,\aaa}(N,\mu) $, respectively ${\cal
A}^L_{s,r,\aaa}(N,\mu,h) $, the closure of the vector space generated by  $(N,\mu)$--low,
resp. $(N,\mu,h) $--low, scalar monomials  in the norm
$\| \ \|_{s,r,\aaa} $ in Definition \ref{MNV}.
\\The projection on ${\cal
A}^L_{s,r,\aaa}(N,\mu,h) $ will be denoted by $\Pi^{L,h}_{N,\mu}$.
Note that it is a projection (see
\eqref{PROJE}) on the subset of indexes $I\subset
\mathbb{I} $ satisfying \eqref{zerobis} and
\eqref{vincula}.
\end{definition}

Clearly, the momentum \eqref{hurra} of a scalar monomial  $\mathfrak m(k,i,\a,\b)$, which is
$(N,\mu)$-low momentum,  satisfies $ | \pi(k,\a,\b) |  \leq \kappa N^b + \mu N^L  $, by \eqref{caracalla}, \eqref{zerobis}.
Hence a scalar monomial $\mathfrak m(k,i,\a,\b)$
 may be   $ (N,\mu,h)$--low  only if
\begin{equation}\label{notabene}
|h| < | \pi(k,\a,\b) |+ N^b<  \mu N^L +(\kappa+1) N^b \stackrel{  \eqref{caracalla}} < N \, .
\end{equation}
In particular
\begin{equation}\label{notabenebis}
{\cal A}^L_{s,r,\aaa}(N,\mu,h)=\emptyset\,,\qquad \forall\, |h|\geq N\,.
\end{equation}

We now define the class of $(N,\theta,\mu)$-linear vector fields.
They are linear combinations of monomial vector fields
supported only on the {\em high components} $ \partial_{{z_m^\pm}} $, $|m| > \theta N $,  which are
linear in the {\em high variables} $ z_n $, $|n| > \theta N $,
and with polynomial coefficients  in the {\em low variables} of degree bounded by $\mu N^L $, $L<1$.
We allow a mild dependence of the coefficients  on the low variables because it is
naturally  generated by commutators.
Finally the momentum and the frequency of  each $(N,\theta,\mu)$-linear monomial vector field is bounded by $ N^b $ with $b<L$.
Since $ b < 1 $ 
these vector fields are   approximately  $x$-independent ($ |k| < N^b $) and diagonal
($ |\pi| < N^b $).  The three scales $ 0 < b < L < 1$ are  `low-high" frequency decomposition
which almost decouples
the interaction between the low variables and the high modes, and it is
used in essential way in the commutator Proposition \ref{festa2}.
We denote by $ e_n $ the multi-index with 
the $n$-th component equal to $ 1 $ and with all the others equal to zero.  

\begin{definition}\label{def:vfm} 
A vector field monomial $\mathfrak m(k,i, \a,\b; \vgot)$ is
\begin{itemize}
\item
{\bf $(N,\mu)$-low} if
\begin{equation}\label{esimia}
 |\pi(k,\a,\b;\vgot)|, |k|<N^b\,,\
\a+  \b=  \g\ \  {\rm with}\ \
{\mathop\sum}_{l\in\Z\setminus \mathcal I} |l| \g_l< \mu N^L\,.
\end{equation}
\item
{\bf $(N,\theta,\mu)$--linear} if
\begin{equation}\label{perla}
\vgot  = z_m^\s , \,
 |\pi(k,\a,\b;\vgot )|,\,|k|<N^b,  \,
\a+  \b=  e_n+\g\  {\rm with} \ |m|, |n|> \theta N, \,
\sum_{l\in\Z\setminus \mathcal I} |l| \g_l < \mu N^L \,.
\end{equation}
\end{itemize}
We denote by  $ \mathcal V_{s,r,\aaa}^L(N,\mu)$, respectively
${\cal L}_{s,r,\aaa}( N, \teta, \mu )$,  the closure  in the norm $\|\  \|_{s,r, \aaa}$ of the vector space generated by
the $(N, \mu)$--low, respectively $(N,\theta,\mu)$--linear,
monomial vector fields. The elements of
$ \mathcal V_{s,r,\aaa}^L(N,\mu)$, resp.
${\cal L}_{s,r,\aaa}( N, \teta, \mu )$,
are called
$(N,\mu)$-low, resp. $(N,\theta,\mu)$--linear, vector fields.

The projections on
$ \mathcal V_{s,r,\aaa}^L(N,\mu)$, resp.
${\cal L}_{s,r,\aaa}( N, \teta, \mu )$,
are denoted by $\Pi^L_{N,\mu}$, resp. $\Pi_{N,\theta,\mu}.$
Explicitely $\Pi^L_{N,\mu}$ and $\Pi_{N,\theta,\mu} $, are the projections (see
\eqref{PROJE}) on the subsets of
indexes $I\subset  \mathbb{I}\times \mathtt V$ satisfying
\eqref{esimia} and \eqref{perla} respectively.
\end{definition}

By \eqref{perla} and \eqref{zerobis}, a $ ( N, \theta, \mu ) $-linear vector
field $X $ has the form
\begin{equation}\label{bibi}
X(v)= \sum_{|m|,|n|> \theta N, \s,\s' =\pm } \!\!  \!\!
X^{\s, m}_{\s',n}(v) z_n^{\s'}\partial_{z_m^\s}
\quad
\mbox{where} \quad   X^{\s,m}_{\s',n} \in
 \mathcal A_{s,r,\aaa}^L(N,\mu,\s m- \s'n) \,.
\end{equation}
By  Definition \ref{LM} and  \eqref{caracalla},
the coefficients
$X^{\s, m}_{\s',n}(v)$ in \eqref{bibi} do not depend on
$z_j, \bar z_j$ with $|j|\geq 6 N^L$.

\begin{lemma} Let $ X \in  \mathcal L_{s, r, \aaa}(N,\teta, \mu) $. Then the coefficients in
\eqref{bibi} satisfy
\be\label{giglio}
X^{\s, m}_{\s',n}=0\qquad \mbox{if}\quad \s \mathtt s(m)
=-\s' \mathtt s(n)
\ee
where $  \mathtt s(m) := {\rm sign}(m) $.
\end{lemma}

\begin{pf}
By \eqref{notabenebis} and
$ | \s m- \s'n | \stackrel{\eqref{giglio}} = $ $ |m| + |n | \stackrel{\eqref{perla}} \geq $ $ 2 \teta N \stackrel{\eqref{caracalla}} >  N $
we get $ \mathcal A_{s,r,\aaa}^L(N,\mu,\s m- \s'n) = \emptyset $.
\end{pf}

\begin{lemma}\label{bruco}
Let $\frak m_{k,i,\a,\b}$ be a scalar monomial (see \eqref{commodo})
such that
\begin{equation}\label{tertulliano}
    \a+\b=:\g\qquad {\rm with}\quad
{\mathop\sum}_{l\in\Z\setminus \mathcal I} |l| \g_l<12  N^L\,.
\end{equation}
Then
$$
\Pi_{N,\theta,\mu} \big(\frak m_{k,i,\a,\b}\, z_n^{\s'}\,
\partial_{z_m^\s}\big)=
   \left\{ \begin{array}{ll}
   \big(\Pi^{L,\s m - \s' n}_{N,\mu}(\frak m_{k,i,\a,\b})\big)\, z_n^{\s'}\,
  \partial_{z_m^\s}
    &
\ \ \ {\rm if} \quad |m|,|n|>\theta N \\
0 & \  \ \ {\rm otherwise.}
 \end{array}\right.
$$
\end{lemma}
\begin{pf}
It directly follows
by \eqref{caracalla}, \eqref{vincula}
and \eqref{perla}.
\end{pf}

\subsection{T\"oplitz vector fields}

We define the subclass of  $(N,\theta,\mu)$-linear  vector fields which are {\em T\"oplitz}.

\begin{definition}\label{matteo_aa} {\bf (T\"oplitz vector field)}
A $(N,\theta,\mu)$-linear vector field $ X \in {\cal L}_{s,r,\aaa}( N, \teta, \mu )  $ is  $(N,\teta,\mu)$-{\em T\"oplitz} if
the coefficients in \eqref{bibi} have the form
\begin{equation}\label{marco}
X^{\s,m}_{\s' ,n}=
X^\s_{\s'} \big(\mathtt s(m), \s m - \s' n \big)
\quad
{\rm for\ some\ \ }  X^\s_{\s'}(\varsigma,h)
\in{\cal A}^L_{s,r,\aaa}(N,\mu,h)
\end{equation}
and $\varsigma\in\{+,-\}$, $h\in\Z$. We denote by
$  \Ta_{s,r,\aaa} ( N, \teta , \mu ) $ 
the space  of  the $ (N,\teta ,\mu)$-{\it T\"oplitz} vector fields.
\end{definition}

The next lemma is used in the proof of Proposition \ref{festa2}.

\begin{lemma}\label{poisb}
Let $ X, Y \in {\mathcal T}_{s,r,\aaa}(N, \teta, \mu ) $ and
$W\in { \mathcal V}^L_{s,r,\aaa}(N,\mu_1)$ with
$\cc<\mu,\mu_1<\CC$. For all
$0<s'<s\,, \, 0<r'<r$
and  $ \teta'  \geq \teta , \mu' \leq \mu $ one has
\be\label{pro11}
\Pi_{N,\teta',\mu'} [X, W ]   \in {\mathcal T}_{s',r', \aaa}(N, \teta', \mu' ) \, .
\ee
If moreover
\begin{equation}\label{ugo}
\mu N^L+(\kappa+1) N^b <( \teta' - \teta) N
\end{equation}
then
\be\label{pr2}
\Pi_{N,\teta',\mu'} [X,Y] \in {\mathcal T}_{s',r',\aaa}(N, \teta', \mu' ) \, .
\ee
\end{lemma}

\begin{pfn}{\sc of \eqref{pro11}.}
By definition (recall \eqref{perla})
we have that $X^{(x)}, X^{(y)}$ and
 $X^{(z_m^\s)}$ vanish if $|m|\leq \theta N$.
Arguing as in  \eqref{notabene} we have that
  $W^{z_j^\s}=0$ if $|j|\geq \mu_1 N^L+(\kappa +1) N^b$.
  Note that only the components
  $[X,W]^{(\mathtt v)}$ with $\mathtt v=z_m^\s$ and $|m|>\theta N$
   contribute
  to $\Pi_{N,\teta',\mu'} [X, W ] $.
Noting that $\theta N>\mu_1 N^L+(\kappa +1) N^b$
(by \eqref{caracalla} and $ N \geq N_0 $)
we have
\begin{equation}\label{prurito}
 [X,W]^{(z_m^\s)}= \partial_x X^{(z_m^\s)} W^{(x)}
 +\partial_y X^{(z_m^\s)} W^{(y)} +
 \sum_{\s_1,|j| < \mu_1 N^L+\kappa N^b}
 \partial_{z_j^{\s_1}} X^{(z_m^\s)} W^{(z_j^{\s_1})} \,.
\end{equation}
 By \eqref{bibi} and \eqref{marco} we get
 $ X^{(z_m^\s)}=\sum_{\s', |n|>\theta N}
X^\s_{\s'} (\mathtt s(m), \s m - \s' n ) z_n^{\s'} $.
Let us consider the first term of the right hand side of
 \eqref{prurito}.
Since
$X^\s_{\s'} (\mathtt s(m), \s m - \s' n ),$
$W^{(x_l)}$
$\in{\cal A}^L_{s,r,\aaa}(N,\mu)$ (recall \eqref{marco}),
 all the monomials in
 $\partial_x X^{\s}_{\s'}(  \mathtt s(m),  \s m - \s' n ) W^{(x)}$
 satisfy \eqref{tertulliano}.
By Lemma \ref{bruco}
 we have
 \begin{eqnarray*}
    &&\Pi_{N,\teta',\mu'} \big(\partial_x X^{(z_m^\s)}
   W^{(x)}\partial_{z_m^\s}\big)=
   \left\{ \begin{array}{ll}
   \sum_{\s', |n|> \theta' N}  U^{\s,m}_{\s',n}\, z_n^{\s'}\,
  \partial_{z_m^\s},
    &
\ \ {\rm if} \quad |m|>\theta' N \\
0 & \ \ {\rm otherwise,}
 \end{array}\right.
 \\
 &&{\rm where}\qquad
 U^{\s,m}_{\s',n}:=
 \Pi^{L,\s m - \s' n}_{N,\mu'}\big(
  \partial_x X^{\s}_{\s'}(  \mathtt s(m),  \s m - \s' n ) W^{(x)}\big)\,.
 \end{eqnarray*}
It is immediate to see that $U^{\s,m}_{\s',n}$ satisfy \eqref{marco}.
The other terms in \eqref{prurito} are analogous. \eqref{pro11} follows.
\\[1mm]
{\sc Proof of \eqref{pr2}.}
We have by \eqref{commco}
\begin{equation}\label{fiodena}
[X,Y]
=:Z-Z'\,, \qquad \mbox{where}\quad
Z:=
\sum_{\s,|m| >\theta N }\Big(
\sum_{\s_1,|j| >\theta N }
\partial_{z_j^{\s_1}} X^{(z_m^\s)} Y^{(z_j^{\s_1})}
\Big) \partial_{z_m^\s}
\end{equation}
and $Z'$ is analogous exchanging
the role of $X$ and $Y$.
We have to prove that
$\Pi_{N,\teta',\mu'} Z \in {\mathcal T}_{s',r',\aaa}(N, \teta', \mu' )$.
By \eqref{bibi} and \eqref{marco} we get
$$
Z^{(z_m^\s)}
=
\sum_{\s_1,|j| >\theta N } \sum_{\s',|n|> \theta N}
 X^{\s}_{\s_1}( \mathtt s(m),\s m - \s_1 j)
 Y^{\s_1}_{\s'}(
\mathtt s(j),\s_1 j-\s' n) z_n^{\s'}\,.
$$
Since both
$X^{\s}_{\s_1}( \mathtt s(m),\s m - \s_1 j)$ and
$Y^{\s_1}_{\s'}(
\mathtt s(j),\s_1 j-\s' n)$
belong to
${\cal A}^L_{s,r,\aaa}(N,\mu)$ (recall \eqref{marco}),
 all the monomials in their product
 satisfy \eqref{tertulliano}.
By Lemma \ref{bruco} we get
$$
 \Pi_{N,\teta',\mu'} Z
 =
\sum_{\s,\s', |m|,|n|> \theta' N}
Z^{\s,m}_{\s',n}\,
 z_n^{\s'}\,\partial_{z_m^\s}
$$
where
\begin{equation}\label{meloperoso}
Z^{\s,m}_{\s',n}:=
\Pi^{L,\s m - \s' n}_{N,\mu'}\Big(
\sum_{\s_1,|j| >\theta N }
 X^{\s}_{\s_1}( \mathtt s(m),\s m - \s_1 j)
 Y^{\s_1}_{\s'}(
\mathtt s(j),\s_1 j-\s' n)\Big)\,.
\end{equation}
Note that
 $X^{\s,\s_1}(  \mathtt s(m),\s m- \s_1 j )
 \in \mathcal A^L(N,\mu,\s m- \s_1 j)$, formula \eqref{notabene}
 and condition \eqref{ugo}  imply that if $|m|>\teta' N$
 then automatically $|j|>  |m| -|\s m- \s_1 j |> \theta' N -  \mu N^L -
 (\kappa+1) N^b > \theta N$ or
 $X^{\s,\s_1}(  \mathtt s(m),\s m- \s_1 j )=0$.
Then the summation in \eqref{meloperoso} runs over $j\in \Z.$
By \eqref{giglio}
we have  $ \mathtt s(j) = \s \s_1 \mathtt s(m)  $. Therefore
$$
Z^{\s,m}_{\s',n}:=
\Pi^{L,\s m - \s' n}_{N,\mu'}\Big(
{\mathop\sum}_{\s_1,h }
 X^{\s}_{\s_1}( \mathtt s(m),h)
 Y^{\s_1}_{\s'}(\s \s_1
\mathtt s(m),\s m-\s' n- h)\Big)
$$
satisfying \eqref{marco}.
\end{pfn}

\subsection{Quasi-T\"oplitz vector fields}

Given a  vector field $ X $ and a T\"oplitz vector
$ {\tilde X} \in \mathcal T_{s,r,\aaa}(N,  \teta,\mu) $ we define 
\begin{equation}\label{limi_aa}
\hat X := N (\Pi_{N, \teta,\mu} X- \tilde X)  \, .
\end{equation}

\begin{definition} \label{topbis_aa} {\bf (Quasi-T\"oplitz)}
A vector field $ X \in \mathcal V_{s,r, \aaa} $  is called
$ (N_0,  \teta, \mu) $-quasi-T\"oplitz  if
the quasi-T\"oplitz norm
\begin{equation}\label{unobisbis_aa}
\| X   \|_{s,r, \aaa}^T := \| X  \|_{s,r, \aaa, N_0,  \teta, \mu}^T:= \sup_{N \geq N_0}
\Big[ \inf_{ {\tilde X} \in \mathcal T_{s,r,\aaa}(N,  \teta,\mu)} \Big( \max \{\| {X}\|_{s,r, \aaa},
\| {\tilde X}\|_{s,r, \aaa},\| {\hat X}\|_{s, r, \aaa} \} \Big) \Big]
\end{equation}
is finite. We define
$$
{\mathcal Q}^T_{s,r, \aaa} ( N_0,  \teta, \mu) :=
\big\{  X : D(s,r) \to V  \, : \,  \| X \|_{s,r, \aaa,N_0,  \teta, \mu}^T < \infty \big\} \, .
$$
\end{definition}

In other words, a vector field $ X $ is
$ (N_0,  \teta, \mu) $-quasi-T\"oplitz with norm $ \| X \|_{s,r, \aaa}^T$
 if,
for all $ N \geq N_0 $,  $ \forall \e > 0 $, there is $ \tilde X \in \Ta_{s,r, \aaa} (N,  \teta, \mu) $ such that
\be\label{defto}
\Pi_{N, \teta,\mu} X = \tilde X + N^{-1} \hat X   \quad {\rm and } \quad
\| X \|_{s,r,  \aaa} \, , \ \|   \tilde X \|_{s,r,  \aaa} \, , \ \|{\hat X} \|_{s,r,  \aaa} \leq \| X \|_{s,r,  \aaa}^T  + \e \, .
\ee
We call $ \tilde X \in \Ta_{s,r, \aaa} (N,  \teta, \mu) $  a
``{\it T\"oplitz approximation}" of $ X $ and $ \hat X $ the ``{\it T\"oplitz-defect}".

If $ s'\leq s, r'\leq r, \aaa'\leq \aaa, N_0'\geq N_0, \theta'\geq \theta, \mu'\leq \mu$
then
 \begin{equation}\label{dino}
 \| \cdot \|_{s',r', \aaa', N_0',  \teta', \mu'}^T
 \leq \max\{ s/ s' , (r/ r')^2 \} \|\cdot \|_{s,r, \aaa,N_0,  \teta, \mu}^T\,.
\end{equation}
\begin{lemma}{\bf (Projections 1)}\label{proJ1}
Consider a subset of indices $ I \subset \mathbb{I}\times \mathtt V $ (see \eqref{I}, \eqref{setV})
such that the projection (see \eqref{PROJE})
\be\label{tototo}
\Pi_I :  \Ta_{s,r, \aaa} (N,  \teta, \mu) \to  \Ta_{s,r, \aaa} (N,  \teta, \mu) \, , \quad \forall N \geq N_0 \, .
\ee
Then $ \Pi_I :  {\cal Q}^T_{s,r,\aaa} (N_0, \teta,\mu) \to  {\cal Q}^T_{s,r,\aaa} (N_0, \teta,\mu) $ and
\be\label{lessore}
\| \Pi_I X \|^{T}_{s,r,\aaa} \leq \| X \|^{T}_{s,r,\aaa} \, . 
\ee
Moreover, if $  X\in {\cal Q}^T_{s,r,\aaa} (N_0, \teta,\mu) $ satisfies
$\Pi_I X =X $, then, $ \forall N \geq N_0 $, $ \forall \e > 0 $,
there exists a  decomposition
$\Pi_{N,\theta,\mu} X=
\tilde X+ N^{-1}\hat X$
with a
 T\"oplitz approximation $\tilde X\in \Ta_{s,r, \aaa} (N,  \teta, \mu) $
satisfying  $\Pi_I \tilde X = \tilde X $,
 $\Pi_I \hat X = \hat X $
 and $ \| \tilde X\|_{s,r,\aaa},  \| \hat X\|_{s,r,\aaa} <   \|X\|^T_{s,r,\aaa}+\e $.
 \end{lemma}
 \begin{pf}
 By \eqref{defto} 
 (recall that $\Pi_{N,\theta,\mu}$ is a projection on an index subset, see Definition \ref{def:vfm})
\be\label{propro}
 \Pi_{N,\theta,\mu} \Pi_I X= \Pi_I \Pi_{N,\theta,\mu} X = \Pi_I \tilde X + N^{-1}\Pi_I \hat X \, .
\ee
Assumption \eqref{tototo} implies that
 $\Pi_I \tilde X \in \Ta_{s,r, \aaa} (N,  \teta, \mu)$ and so $\Pi_I \tilde X$ is a T\"oplitz approximation for $\Pi_I X $.
Hence \eqref{lessore} follows by
$ \| \Pi_I  X\|_{s,r,\aaa},  \| \Pi_I\tilde X\|_{s,r,\aaa},  \| \Pi_I\hat X\|_{s,r,\aaa}
 <     \|X\|^T_{s,r,\aaa} + \e
$
using Lemma \ref{lem:pro} and \eqref{defto}.
Now, if $\Pi_I X = X $, then \eqref{propro} shows that  $ \Pi_I \tilde X $ (which satisfies $ \Pi_I (\Pi_I \tilde X) =  \Pi_I \tilde X $),
is a T\"oplitz approximation for $ X $.
\end{pf}

For  a vector field $ X : D(s,r) \to V  $ depending on parameters $  \xi \in {\cal O } $,  we define the norm
\begin{equation}\label{fanfulla}
    \|X\|^T_{\vec p}:=\max \big\{
\sup_{\xi \in \mathcal O}
 \| X (\cdot;\xi)\|^T_{s,r,\aaa,N_0,\theta,\mu}\,,
\,
\|X\|_{s,r,\aaa,\mathcal O}^\l   \big\}
\end{equation}
where, for brevity,
\begin{equation}\label{pipino}
  \vec p:=(s,r,\aaa,N_0,\theta,\mu,\l,\mathcal O) \, .
\end{equation}
We denote
\begin{equation}\label{priscilla}
    {\mathcal Q}^T_{\vec p}:= \big\{
X\in {\cal V}_{s,r,\aaa,\mathcal O}^\l \ \, : \ \,
X(\cdot;\xi) \in {\mathcal Q}^T_{s,r, \aaa} ( N_0,  \teta, \mu)\, , \
\forall\, \xi \in \mathcal O \ {\rm and}\  \|X\|^T_{\vec p}<\infty \big\} \, .
\end{equation}
In view of the KAM step we prove that the quasi-T\"oplitz
norm does not increase under suitable
projections and that it satisfies smoothing estimates.
We denote by $  \Pi_{\rm diag} $  the projection on the space generated by the
 monomial vector fields $ z_j \partial_{z_j}$, $ \bar z_j \partial_{ \bar z_j} $.

\begin{lemma} {\bf (Projections 2)}\label{proJ}
For all $ l \in \N $,  $ K \in \N $, $ N \geq N_0 $, the projections (see \eqref{phom}, \eqref{minas}, \eqref{tirith}) 
map
\begin{equation}\label{margherita}
 \Pi^{(l)},  \Pi_{|k| < K}, \Pi_{|\pi| < K}, \Pi_{\rm diag} : \Ta_{s,r, \aaa} (N,  \teta, \mu) \to \Ta_{s,r, \aaa} (N,  \teta, \mu) \, .
\end{equation}
If $ X\in {\cal Q}^T_{\vec p} $
 then
\be\label{fhT}
\| \Pi^{(l)} X \|^{T}_{\vec p} \, ,  \ \| \Pi_{|\pi| < K} X \|^T_{\vec p}  \, , \| \Pi_{\rm diag} X \|^T_{\vec p},
\| X^{\leq 0} \|^{T}_{\vec p} \, , \ \| X - X^{\leq 0}_{|k|<K} \|^{T}_{\vec p} \leq \|  X \|^T_{\vec p} \, .
\ee
Moreover, $ \forall\, 0 < s' < s $ and $ \forall 0<\aaa'<\aaa$, setting $\vec p'=(s',r,\aaa',N_0,\theta,\mu,\lambda,\mathcal O)$ :
\begin{eqnarray}\label{smoothT}
\| \Pi_{|k|\geq K} X \|_{\vec p'}^T \leq
e^{-K(s-s')} (s \slash s' ) \|  X \|_{\vec p}^T \, , \quad
 \| \Pi_{|\pi|\geq K} X \|_{\vec p'}^T \leq
e^{-K(\aaa-\aaa')} \|  X \|_{\vec p}^T \, .
\end{eqnarray}
\end{lemma}
\begin{pf}
We prove \eqref{margherita} for $ \Pi_{|\pi|<K}$, the others are analogous.
Since $\tilde X \in \Ta_{s,r, \aaa} (N,  \teta, \mu)$ then
$ \tilde X(v) = \sum_{\s,\s', |m|,|n|> \theta N}
\tilde X^{\s, m}_{\s',n}(v) z_n^{\s'}\partial_{z_m^\s} $
for some $\tilde X^{\s, m}_{\s',n}$ satisfying \eqref{marco}.
Then
$$
\big( \Pi_{|\pi|<K} \tilde X\big) (v)=
\sum_{\s,\s', |m|,|n|> \theta N}
Y^{\s, m}_{\s',n}
(v) z_n^{\s'}\partial_{z_m^\s}\quad
{\rm where}\quad
Y^{\s, m}_{\s',n}
:=\Pi_{|\pi  + \s' n - \s m |<K} 
\tilde X^{\s, m}_{\s',n}
$$
(recall Definition \ref{LM}).
Therefore $Y^{\s, m}_{\s',n}$ satisfy \eqref{marco}
and
$\Pi_{|\pi|<K}\tilde X \in \Ta_{s,r, \aaa} (N,  \teta, \mu)$.
The estimates \eqref{fhT} follow from \eqref{margherita} and Lemma \ref{proJ1}
(in particular \eqref{lessore}). The bounds \eqref{smoothT}
follow by \eqref{smoothl}, \eqref{torquato} and similar arguments.
\end{pf}

The following proposition  shows that
the quasi-T\"oplitz vector fields satisfy, slightly modulating
the parameters, the Lie algebra property. 

\begin{proposition}\label{festa2} {\bf (Lie bracket)}
Assume that $ X^{(1)}, X^{(2)} \in {\mathcal Q}^T_{\vec p} $
(see \eqref{priscilla}) and assume that
$ \vec p_1:= $ $ (s_1,r_1,\aaa_1, N_1, \theta_1,\mu_1,\l,\mathcal O) $
with $ N_1 \geq N_0 $, $ \mu_1 \leq \mu $, $  \teta_1
\geq  \teta $, $ s/2 \leq s_1 < s $, $ r/2 \leq r_1 < r $, $ \aaa_1 < \aaa $,  satisfy
\begin{eqnarray}
&&(\kappa+1) N_1^{b-L}  < \mu - \mu_1, \ \mu_1
N_1^{L-1} + (\kappa+1) N_1^{b-1} <  \teta_1  -  \teta\,, \label{Leviatan0}
\\
&&2 N_1 e^{- N_1^b \min\{\aaa- \aaa_1 , s -s_1\} /2} < 1 \, , \
b \min\{\aaa- \aaa_1 , s -s_1\} N_1^b > 2 \,  .
\label{Leviatan}
\end{eqnarray}
Then $[X^{(1)}, X^{(2)} ]\in {\mathcal Q}^T_{ \vec p_1} $
and, for some  $ C(n) \geq 1 $,
\begin{equation}\label{cippalippa}
 \| [X^{(1)}, X^{(2)} ] \|^T_{\vec p_1} \leq  C(n) \delta^{-1}  \| X^{(1)} \|^T_{\vec p}  \| X^{(2)} \|^T_{\vec p} \, ,
 \quad
\d := 
\min \Big\{ 1 - \frac{s_1}{s}, 1 - \frac{r_1}{r}  \Big\}\,.
\end{equation}
\end{proposition}
 The main point in the proof of the above proposition is the
 following purely algebraic result.

\begin{lemma}\label{nn} {\bf (Splitting lemma)}
Let $ X^{(1)}, X^{(2)} \in {\mathcal V}_{s,r,\aaa} $
and \eqref{Leviatan0} hold.
Then, for all $ N \geq N_1 $,
\begin{eqnarray}
& &
\Pi_{N, \teta_1,\mu_1}[X^{(1)}, X^{(2)} ]  = \nonumber  \\
&  & \Pi_{N, \teta_1,\mu_1} \big(  \big [ \Pi_{N, \teta,\mu}
X^{(1)} ,\Pi_{N, \teta,\mu} X^{(2)} \big] +  \big[ \Pi _{N,
\teta,\mu} X^{(1)}, \Pi^L_{N,\mu} X^{(2)} \big]   +\big[ \Pi^L
_{N,\mu} X^{(1)}, \Pi_{N,\teta,\mu} X^{(2)} \big] \qquad \qquad \label{ultima2}  \\
&& \qquad\quad \, + \, \big[  \Pi_{|k| \geq N^b\, {\rm or}\, |\pi| \geq  N^b} X^{(1)},  X^{(2)} \big]
+ \big[ \Pi_{|k|, |\pi| < N^b} X^{(1)},   \Pi_{|k| \geq N^b\, {\rm or}\, |\pi| \geq  N^b} X^{(2)}
\big] \big)
 \, . \label{ultima3}
\end{eqnarray}
\end{lemma}

Then the proof of Proposition \ref{festa2}
  follows as in \cite{BBP1} (see Proposition 3.1).
The point is to find a T\"oplitz approximation and
a T\"oplitz defect of
$\Pi_{N, \teta_1,\mu_1}[X^{(1)}, X^{(2)} ]$, recall \eqref{defto}.
A T\"oplitz approximation is obtained by \eqref{ultima2}
substituting $\Pi_{N, \teta,\mu}
X^{(i)} $, $i=1,2$, with their T\"oplitz approximations, thus 
yielding a vector field which is T\"oplitz  by
Lemma \ref{poisb}.
The remaining terms in \eqref{ultima2} are T\"oplitz defects. They are small
because contain commutators with the  T\"oplitz
defects of $\Pi_{N, \teta,\mu}
X^{(i)} $. 
The last terms  \eqref{ultima3}
 are exponentially small 
 by  
 \eqref{Leviatan} and \eqref{smoothT}.
The momentum-norms of the commutators
are estimated by Proposition \ref{settimiosevero}.
\smallskip

\begin{pfn}{\sc{of Lemma \ref{nn}.}}
We have
\begin{eqnarray}
[ X^{(1)}, X^{(2)} ] & = & [ \Pi_{|k|, |\pi| < N^b} X^{(1)},  \Pi_{|k| ,  |\pi| < N^b} X^{(2)} ]  \label{pr1}\\
& + & [  \Pi_{|k| \geq N^b \,{\rm or}\, |\pi| \geq  N^b} X^{(1)},  X^{(2)} ] + [  \Pi_{|k| , |\pi| < N^b} X^{(1)}, \Pi_{|k| \geq N^b\, {\rm or}\, |\pi| \geq  N^b} X^{(2)} ] \nonumber \, .
\end{eqnarray}
The last two terms are
  \eqref{ultima3}. We now prove that the right hand side
  of \eqref{pr1} gives the three terms in \eqref{ultima2}.
It is sufficient to study the case where $X^{(h)} $, $ h = 1, 2 $, are monomial vector fields
\be\label{mon12}
\mathfrak m_h = m_{k^{(h)}, i^{(h)}, \alpha^{(h)}, \beta^{(h)}; \vgot^{(h)}}  \ {\rm (see \ \eqref{commodo})} \
\quad {\rm with} \quad
|k^{(h)}|, |\pi(\mathfrak m_h )| < N^b\,,\; h=1,2 \, ,
\ee
and  analyze under which conditions the projection
$\Pi_{N, \teta_1,\mu_1} [\mathfrak m_1,  \mathfrak m_2 ]  $
is not zero.

By the formula of the commutator \eqref{commco} and
the definition of the projection $ \Pi_{N, \teta_1, \mu_1} $ (see Definition \ref{def:vfm}, in particular \eqref{perla})
we have to compute  $(D_\vgot \mathfrak m_1^{\vgot'} )[\mathfrak m_2^\vgot]$ only for
 $ \vgot' = z_m^{\s} $ with $|m|>\theta_1 N $ and
$ \vgot \in \mathtt V$, see \eqref{setV}.
\\[1mm]
$\bullet)$
{\sc case 1}: $ \vgot = x_i$ or $ \vgot = y_i $. By \eqref{perla}, in order to have a
non trivial projection  $ \Pi_{N, \teta_1, \mu_1}  (D_\vgot \mathfrak m_1^{z_m^{\s}} )[\mathfrak m_2^\vgot] $ it must be
\be\label{nlarge}
\alpha^{(1)}+ \beta^{(1)}+\alpha^{(2)}+ \beta^{(2)} =
e_n + \gamma\,,\quad |n|>\theta_1 N\,, \ \
{\mathop\sum}_{l\in \Z\setminus \mathcal I} |l| \gamma_l< \mu_1 N^L\,.
\ee
We claim that
\begin{equation}\label{cipolla0}
\alpha^{(1)}+ \beta^{(1)}
=e_n+\gamma^{(1)}\,,\quad
\alpha^{(2)}+ \beta^{(2)}
=\gamma^{(2)}
\,,\quad
{\mathop\sum}_{l\in \Z\setminus \mathcal I} |l| \gamma_l^{(h)}< \mu_1 N^L\, , \
h=1,2\,,
\end{equation}
which implies that $ \mathfrak m_1$ is $ (N,\theta_1,\mu_1)$--linear (see \eqref{perla}), hence
$ (N,\theta,\mu)$--linear, and $\mathfrak m_2 $ is $(N,\mu_1) $--low (see \eqref{esimia}), hence
$(N,\mu)$--low. Thus
$ \Pi_{N,\theta,\mu} \mathfrak m_1=\mathfrak m_1$ and
$\Pi^L_{N,\mu}\mathfrak m_2=\mathfrak m_2 $ and we obtain
the second  (and third by commuting indices)  term in the right hand side of \eqref{ultima2}.
By \eqref{nlarge},  the other possibility instead of  \eqref{cipolla0} is
\be\label{assur1}
\alpha^{(1)}+ \beta^{(1)}
=\tilde\gamma^{(1)}\,,\quad
\alpha^{(2)}+ \beta^{(2)}
= e_n +  \tilde\gamma^{(2)} \,,\quad
{\mathop\sum}_{l\in \Z\setminus \mathcal I} |l| \tilde\gamma_l^{(h)}<
\mu_1 N^L\, , \ h=1,2\,.
\ee
In such a case, since  $ |\pi(\mathfrak m_2)| < N^b $ we get (recall $ \mathfrak m_2 = \mathfrak m_2^{\vgot} $ with $ \vgot = x, y $),
$$
N^b > |\pi (k^{(2)},\alpha^{(2)},\beta^{(2)})|
\stackrel{\eqref{momp}, \eqref{assur1}} \geq
|n|- \sum_{l\in \Z\setminus \mathcal I} |l| \tilde\gamma_l^{(2)}
-\kappa |k^{(2)}| \stackrel{\eqref{nlarge}, \eqref{assur1}, \eqref{mon12}} \geq \theta_1 N- \mu_1 N^L -\kappa N^b
$$
which contradicts \eqref{caracalla}.
\\[1mm]
$\bullet)$
{\sc case $ 2$:} $ \vgot = z^{\s_1}_j$, $ j \in\Z\setminus \mathcal I.$
Only for this case we use  \eqref{Leviatan0}.
In order to have a non trivial projection  $ \Pi_{N, \teta_1, \mu_1}  (D_\vgot \mathfrak m_1^{z_m^{\s}} )[\mathfrak m_2^\vgot] $
it must be
\be\label{split2}
\alpha^{(1)}+ \beta^{(1)}+\alpha^{(2)}+ \beta^{(2)}-e_j
=e_n+\gamma\,,\quad |n|>\theta_1 N\,, \ \ \
{\mathop\sum}_{l\in \Z\setminus \mathcal I} |l| \gamma_l< \mu_1 N^L\,.
\ee
We have the two following  possible cases:
\begin{align}
&
\alpha^{(1)}+ \beta^{(1)}
=e_j+e_n+\gamma^{(1)}\,,\quad
\alpha^{(2)}+ \beta^{(2)}
=\gamma^{(2)}
\,,\quad
{\mathop\sum}_{l\in \Z\setminus \mathcal I} |l| \gamma_l^{(h)}< \mu_1 N^L\, ,  \
h=1,2
\label{ernesto}
\\
&
\alpha^{(1)}+ \beta^{(1)}
= e_j + \tilde \gamma^{(1)}\,,\quad
\alpha^{(2)}+ \beta^{(2)} = e_n+\tilde\gamma^{(2)}
\,,\quad
{\mathop\sum}_{l\in \Z\setminus \mathcal I} |l| \tilde\gamma_l^{(h)}< \mu_1 N^L\, ,  \
h=1,2
\label{evaristo}
\end{align}
where
$ \gamma^{(1)}+\gamma^{(2)}=\tilde\gamma^{(1)}+\tilde\gamma^{(2)} = \gamma $.
Note that, since we  differentiate  $\mathfrak m_1 $ with respect to $ \vgot =z^{\s_1}_j $
the monomial $ \mathfrak m_1 $ must depend on $  z^{\s_1}_j $ and so the following case does not arise:
$$
\alpha^{(1)}+ \beta^{(1)}
=  \tilde \gamma^{(1)}\,,\quad
\alpha^{(2)}+ \beta^{(2)} = e_j + e_n+\tilde\gamma^{(2)}
\,,\quad
{\mathop\sum}_{l\in \Z\setminus \mathcal I} |l| \tilde\gamma_l^{(h)}< \mu_1 N^L\, ,  \
h=1,2 \, .
$$
In the case \eqref{ernesto}, the monomial $\mathfrak m_2$ is $(N,\mu)$--low and
we claim that  $\mathfrak m_1$ is $(N,\theta,\mu)$--linear.
Indeed, since
\be\label{pm2b}
|\pi (\mathfrak m_2)| \stackrel{\eqref{momp}} = | \pi (k^{(2)},\alpha^{(2)},\beta^{(2)}) - \s_1 j | <N^b
\ee
we get $| j|\leq | \pi (k^{(2)},\alpha^{(2)},\beta^{(2)})|+ N^b $. Hence
\begin{eqnarray*}
|j|+ {\mathop\sum}_l \gamma_l^{(1)} |l| & \leq &
|\pi (k^{(2)},\alpha^{(2)},\beta^{(2)})|+ N^b+
{\mathop\sum}_l \gamma_l^{(1)} |l|
\leq \kappa |k^{(2)}|+ {\mathop\sum}_l\gamma_l |l| +N^b \\
& \stackrel{\eqref{mon12}, \eqref{split2}} \leq &  (\kappa +1) N^b+ \mu_1 N^L
\stackrel{\eqref{Leviatan0}}\leq \mu N^L
\end{eqnarray*}
namely $ \mathfrak m_1 $ is $(N,\theta,\mu)$--linear (see \eqref{perla} with $ \g = e_j  + \g^{(1)} $).
Hence
$ \Pi_{N,\theta,\mu} \mathfrak m_1=\mathfrak m_1$ and
$\Pi^L_{N,\mu}\mathfrak m_2=\mathfrak m_2 $ and we obtain
the second term (and third by commuting indices) in the right hand side of \eqref{ultima2}.

In the case \eqref{evaristo} we claim that  both $\mathfrak m_1,\mathfrak m_2 $ are $(N,\theta,\mu)$--linear
so we obtain the first term in the right hand side of \eqref{ultima2}. Since,  by \eqref{split2},
$ | n | > \teta_1 N > \teta N $ we already know that $ \mathfrak m_2 $ is $(N,\theta,\mu)$--linear.
Finally,  
$\mathfrak m_1  $ is $(N,\theta,\mu)$--linear because
\begin{eqnarray*}
|j| \stackrel{\eqref{pm2b}} > |\pi (k^{(2)},\alpha^{(2)},\beta^{(2)})| -N^b \! \! \! \! & \stackrel{\eqref{momp}, \eqref{evaristo}} \geq & \! \! \! \!
|n|- {\mathop\sum}_{l\in \Z\setminus \mathcal I} |l| \tilde\gamma_l^{(2)}
-\kappa |k^{(2)}| - N^b \\
\! \!  \! \! & \stackrel{ \eqref{split2}, \eqref{evaristo} , \eqref{mon12} } > & \! \! \! \theta_1 N- \mu_1 N^L -(\kappa+1) N^b
\stackrel{\eqref{Leviatan0}}>\theta N
\end{eqnarray*}
concluding the proof.
\end{pfn}

The quasi-T\"oplitz character of a vector field is preserved under the flow
of a quasi-T\"oplitz vector field. As the corresponding Proposition 3.2 of \cite{BBP1}, the proof is  an iteration of
Proposition \ref{festa2}.
\begin{proposition}\label{main2} {\bf (Lie series)}
Let $ X,Y \in {\mathcal Q}^T_{\vec p} $
(see \eqref{priscilla}).
Assume $
\vec p  \, ':=(s',r',\aaa', N_0', \theta',\mu',\l,\mathcal O) $ satisfies
$ s/2 \leq s' < s $, $ r/2\leq r' < r $, $\aaa'<\aaa$ , $ \mu' < \mu $, $  \teta' >  \teta $, and
\begin{equation}\label{stoppa}
    N_0' \geq \max \{N_0, \bar N \} \, , \quad
    \bar N := \exp \big(
\max\big\{ 2 / b, (L-b)^{-1}, (1-L)^{-1},8
\big\} \big)\, ,
\end{equation}
\begin{equation}\label{Giobbe}
   ( \kappa +1) (N_0')^{b-L} \ln N_0'  \leq \mu - \mu' \,, \ \
(7  + \kappa) (N_0')^{L-1} \ln N_0' \leq   \teta'  -  \teta \, ,
\end{equation}
\begin{equation}\label{Giobbe1}
 2(N_0')^{-b}\ln^2 N_0' \leq b \min \{ s-s', \aaa - \aaa' \} \,.
\end{equation}
There is $  c(n) >  0  $ such that, if the smallness condition
\be\label{piccof2}
\| X \|^T_{\vec p}  \leq c(n) \, \d
\ee
holds (with $ \d $ defined in \eqref{diffusivumsui}),
then $ e^{ {\rm ad}_X } Y \in {\mathcal Q}^T_{\vec p \, '} $ and
\be\label{gPhifbis}
\| e^{{\rm ad}_X} Y \|^T_{\vec p  \, '} \leq 2 \| Y \|^T_{\vec p} \, .
\ee
Moreover, for  $ h  = 0,1, 2 $, and  coefficients $ 0 \leq b_j\leq 1/j!$, $j\in\mathbb{N},$
\be\label{gPhif12bis}
\Big\| {\mathop\sum}_{j \geq h}  b_j \,{\rm ad}_X^j (Y)
\Big\|^T_{{\vec p}  \, '} \leq 2 (C\d^{-1} \| X\|^T_{\vec p})^{h} \| Y \|^T_{\vec p} \, .
\ee
\end{proposition}

\section{An abstract KAM theorem}\label{sec:4}\setcounter{equation}{0}

We consider a family of linear integrable vector fields with constant coefficients
\begin{equation}\label{90210}
{\mathcal N}  (\xi) :=   \oo (\xi) \partial_x +
 \ii \OO (\xi) z \partial_z  -  \ii   \OO  (\xi) {\bar z} \partial_{\bar z}
\end{equation}
defined on the phase space
$ \mathbb{T}^n_s  \times \mathbb{C}^n
\times \ell^{a,p}_{{\cal I}} \times \ell^{a,p}_{{\cal I}}  $,
where the tangential sites $ {\cal I}
 \subset \Z $ are symmetric as in \eqref{Isym},
the  space $\ell^{a,p}_{\cal I} $
is defined in \eqref{ellcal},  the tangential  frequencies
$ \oo  \in \R^n $
and the normal frequencies $   \OO \in \R^{\Z \setminus {\cal I}} $
depend on real parameters
$ \xi \in {\cal O} \subset \mathbb{R}^{n/2}$ 
(where $ n / 2 = $ cardinality of $ {\cal I}^+ $, see \eqref{Isym}),
and satisfy
\be\label{sio}
{\oo}_j (\xi ) =
{ \oo}_{-j} (\xi) \, , \ \forall j \in {\cal I } \, ,  \quad
 \OO_j (\xi ) =   \OO_{-j} (\xi) \, , \ \forall j \in \Z \setminus {\cal I} \, .
\ee
For each $ \xi $ there is an invariant $ n$-torus ${\cal T}_0 = \mathbb{T}^n \times \{0 \} \times \{0 \} \times \{0 \} $
with frequency $ \oo ( \xi ) $. In its normal space, the origin $ (z,\bar z ) = 0  $ is an elliptic fixed point with
proper frequencies $   \OO (\xi ) $.
The aim is to prove the persistence of a large portion of this family  of linearly stable
tori under small perturbations
\be\label{perturbVF}
\PP(x,y,z, \bar z; \xi) = \PP^{(x)} \partial_x + \PP^{(y)} \partial_y +
 \PP^{(z)} \partial_z + \PP^{(\bar z)} \partial_{\bar z} \, .
\ee

\begin{description}
\item $ {\bf (A1)} $ {\sc Parameter dependence}.
The map $  \oo : {\cal O} \to {\mathbb R}^n $,
$ \xi \mapsto  \omega ( \xi )$,  is Lipschitz continuous.
\end{description}
With in mind the application to DNLW we assume
\begin{description}
\item $ {\bf (A2)} $
{\sc Frequency asymptotics}.
\begin{equation}\label{carbonara}
\OO_j (\xi) =   |j| + a(\xi) +  b(\xi) |j|^{-1} + O( j^{-2} )  \ \ \ {\rm as} \  \ \  |j| \to + \infty \, .
\end{equation}
Moreover the map
$  ( \OO_j - |j|)_{j \in \Z \setminus {\cal I}} : {\cal O} \to \ell_\infty  $ is Lipschitz continuous.
\end{description}

By (A$1$) and (A$2$),  the Lipschitz semi-norms
of the frequency maps satisfy, for some $1\leq M_0 < \infty $,
\begin{equation}\label{M}
|\oo|^{\rm lip} + |   \OO |^{\rm lip}_{\infty} \leq M_0 \qquad
{\rm where} \qquad
|   \OO |^{\rm lip}_{\infty} := \sup_{\xi \neq \eta \in \cal O}
\frac{ |\OO (\xi) - \OO (\eta )|_\infty }{| \eta - \xi |} \,
\end{equation}
and $ | z |_\infty := \sup_{j \in \Z \setminus {\cal I}} | z_j | < + \infty $.
\begin{description}
\item $ {\bf (A3)} $  {\sc Regularity}.
The  vector field $ \PP $ in  \eqref{perturbVF} maps
$ \PP : D (s,r) \times {\cal O} \to $ $
\C^n \times \C^n \times \ell^{a,p}_{{\cal I}}   \times \ell^{a,p}_{{\cal I}} $
 for some $ s,r>0 $.
Moreover $ \PP $ is {\sc reversible}
(Definition \ref{reversVF}), {\sc real-coefficients},
{\sc real-on-real}  {\sc Even} (Definition \ref{rean}).
\end{description}

Finally, in order to obtain the asymptotic expansion for the perturbed frequencies we also assume
\begin{description}
\item $ {\bf (A4)} $  {\sc Quasi-T\"oplitz}.
The perturbation vector field $ \PP $ is
quasi-T\"oplitz, see Definition \ref{topbis_aa}. 
\end{description}

Recalling  \eqref{perturbVF} and
the notations in \eqref{baronerosso}, \eqref{phom},
 we define
\be\label{primopezzo}
\PP^y(x)\partial_y :=\Pi^{(-1)} \PP^{(y)} \partial_y\,,
\qquad  \PP_*:=\PP-\PP^y (x) \partial_y
\ee
and we  denote by $\PP_*^{(-1)},\PP_*^{(0)}$  the terms of degree $ -1 $ and $ 0 $
respectively of $ \PP_* $, see \eqref{phom}.
Let
\begin{equation}\label{pingu}
\ot (\xi)  := (\oo_j (\xi)  )_{j \in {\cal I}^+ }
\in\mathbb R^{n/2}\,,\qquad
{\rm then}\ \ \   \oo=(\ot,\ot)\ \
{\rm by}\ \  \eqref{sio}\,.
\end{equation}
\begin{theorem} {\bf (KAM theorem)} \label{thm:IBKAM}
Fix
 $ s, r, \mathtt a >  0 $, $  \cc < \teta,\mu < \CC   $,
 $ N_0 \geq \bar N $ (defined in \eqref{stoppa}).
Let $ \g \in (0,\g_*),$ where $\g_*=\g_*(n,s,\mathtt a)<1$
is a (small) constant.
Let
$ \l := \g / M_0 $ (see \eqref{M})
and $\vec p := (s,r,\aaa,N_0,\theta,\mu,\l,\mathcal O)$. 
Suppose that the vector field $ \XX = {\cal N} + \PP  $ satisfies ${\rm (A1)}$-${\rm (A4)}$.
If
\be\label{KAMcondition}
 \g^{-1} \| \PP_* \|_{\vec p}^T\leq 1
\quad \mbox{and}\quad
 \varepsilon := \max \Big\{ \gamma^{- 2/3}
  \| \PP^y (x) \partial_y \|_{s,r, \mathtt a,   {\cal O}}^\l,
 \g^{-1} \| \PP_*^{(-1)}\|_{\vec p}^T,
 \g^{-1} \| \PP_*^{(0)}\|_{\vec p}^T  \Big\}
\ee
is small enough, then 
\\[1mm] $ \bullet $ {\bf (Frequencies)} There exist
Lipschitz functions $ \o^\infty :  \R^{n/2} \to {\mathbb R}^n $, $ \O^\infty : \R^{n/2} \to  \ell_\infty $,
$ a^{\infty}: \R^{n/2} \to {\mathbb R} $ 
(recall that $ {\cal O} \subset \mathbb{R}^{n/2} $) 
such that 
$ \o^\infty=(\ot^\infty,\ot^\infty) $, $ \ot^\infty := ( \o^\infty_j )_{j \in {\cal I}^+} \in\R^{n/2} $, 
and 
\begin{equation}\label{muflone}
| \o^\infty -  \oo |+\l | \o^\infty -  \oo |^{\rm lip} \, ,
\ \  | \O^\infty -  \OO  |_{\infty}+
\l| \O^\infty -   \OO  |_{\infty}^{\rm lip}\
 \leq \ C \g \varepsilon \, , \  | a^{\infty} |\leq C \g \e \, ,
\end{equation}
$$
\o^\infty_j (\xi ) = \o^\infty_{-j} (\xi) \, , \ \forall j \in {\cal I } \, , \quad
\Om_j^\infty (\xi ) = \Om_{-j}^\infty (\xi) \, , \ \forall j \in
\Z \setminus {\cal I} \, ,
$$
\be\label{freco} \sup_{\xi \in
\R^{n/2} } | \O^\infty_j (\xi) -  \OO_j(\xi) -
 a^\infty (\xi) |
\leq   \g^{2/3} \e \, \frac{C}{|j|} \, , \quad \forall |j| \geq C_\star \g^{-1/3} \, .
\ee
$ \bullet $  {\bf (KAM normal form)}
for every $\xi$ belonging to
\begin{eqnarray}\label{Cantorinf}
{\cal O}_\infty & := & \Big\{ \xi \in {\cal O} \ : \
 \forall h \in \Z^{n/2}, \, i, j
 \in \Z \setminus {\cal I}, \, p \in \Z \, ,
 \nonumber \\
  & &
\ \, | \ot^\infty(\xi) \cdot h  +  \O_j^\infty |
\geq  2 \g \langle h \rangle^{-\t} \, ,  \ \
   |\ot^\infty(\xi) \cdot h  + \O_i^\infty(\xi)  + \O_j^\infty(\xi)  |
   \geq 2  \g \langle h \rangle^{-\t} \, ,
   \nonumber \\ 
  & &
\ \, | \ot^\infty(\xi) \cdot h  -  \O_i^\infty(\xi)  +  \O_j^\infty(\xi) |
  \geq 2 \g \langle h \rangle^{-\t}
  \quad    {\rm if}\quad  h \neq 0 \   {\rm or}\
   i\neq \pm j  \, , \nonumber \\
& &
\ \,  |\ot^\infty (\xi)\cdot h +p|
\geq 2\g^{2/3} \langle h \rangle^{-\t} , \
{\rm if} \ (h,p) \neq (0,0)  \nonumber  \\
& &
\ \, | \ot (\xi ) \cdot h |\geq 2 \g^{2/3} \langle h \rangle^{- n /2}\,,\ \forall\, 0 < |h| < \g^{-1 / (7n)}
\label{clinea} \Big\}
\end{eqnarray}
there exists an even, analytic, close to the identity diffeomorphism
\be\label{eq:Psi}
\Phi(\cdot;\xi) : D(s/4, r/4)
\ni(x_\infty,y_\infty, z_\infty, {\bar z}_\infty)
\mapsto (x,y,z, \bar z)\in D(s, r)\,,
\ee
(Lipschitz in $\xi$)
such that the transformed vector field
\be\label{Hnew}
\XX_\infty= {\mathcal N} _{\infty} + {\cal P}_\infty
:=
\Phi_\star(\cdot;\xi) \XX = (D \Phi ( \ ; \xi ))^{-1}
\XX \circ \Phi ( \ ; \xi )
\quad {has} \quad \big(\PP_\infty^{\leq 0}\big)_{|E} = 0 \, ,
\ee
see  \eqref{Pleq2}, \eqref{pandistelle}.
Moreover
$ {\cal N}_{\infty} 
$ is a constant coefficients linear normal form vector field as \eqref{90210}
with frequencies $ \o^\infty (\xi) $,  $ \O^\infty (\xi) $, and $ \PP_\infty $ is reversible, real-coefficients,
real-on-real, even. Finally $ (\XX_\infty)_{|E} =  ({\cal S} \XX_\infty)_{|E} $.
\end{theorem}

As a consequence we derive

\begin{corollary}
For all $  \xi \in {\cal O}_\infty $, the  map
$  \T^{n/2} \ni \vec x_\infty \mapsto \Phi \big((\vec x_\infty,\vec x_\infty), 0,0,0;\xi\big) \in E $
is an $n/2$-dimensional
analytic invariant torus of the vector field $ \XX = {\cal N } + \PP $.
Such torus is  {\it linearly stable on $ E $}
and, in particular, it has {\it zero Lyapunov exponents on $E$}.
\end{corollary}

The set $ {\cal O}_\infty $   in \eqref{Cantorinf} could be empty. In the next theorem we bound
its measure.

\begin{theorem}{\bf (Measure estimate)} \label{thm:measure}
Let $ {\cal O} := {\cal O}_\rho := \big\{ \xi := ( \xi_j )_{j \in {\cal I}^+} \in \R^{n/2} \, : \,  0
< \rho \slash 2 \leq |\xi_j| \leq \rho   \big\} $.
Assume that the frequencies are affine functions of $ \xi $ 
\be\label{omegaxi}
\ot (\xi) = 
{\bar \om} + A \xi \, , \ {\bar \om}= (\l_j)_{j\in\mathcal I^+}
 \in \R^{n/2} \, ,
\quad \OO_j (\xi ) = \l_j + \l_j^{-1} \,\vec a \cdot \xi   \, , \
\forall j \notin {\cal I} \, ,
\ee
where   $A \in {\rm Mat}(n/2\times n/2)$, $ {\rm det} A \neq 0 $, and $\vec a \in \R^{n/2}$ 
are continuous functions in $\mm.$
Fix a compact interval of masses $[\mm_1,\mm_2]\subset (0,\infty)$
and take $ \mm  \in [\mm_1,\mm_2] $ such that 
\begin{equation}\label{pota}
 (\l_i^{-1} \pm \l_j^{-1})(A^T)^{-1}
  \vec  a \,,\ \
 \l_j^{-1} (A^T)^{-1}  \vec  a
  \notin \Z^{n/2} \setminus \{ 0 \} \, ,
 \ \forall i,j \in \Z \setminus {\cal I}\,,
\ |i|,|j|\leq C_0
 \, ,
\end{equation}
for a suitably large constant $C_0 := C_0 (\mm_1,\mm_2, A, \vec a, \bar \om) $.
Then  the Cantor like set  $ {\cal O}_\infty $ defined in \eqref{Cantorinf},
 with exponent
$ \t > \max\{n + 3, 1/ b \} $ ($ b $ is fixed in \eqref{figaro}),
satisfies, for $ \rho \in (0, \rho_0 (\mm)) $ small,
\be\label{consolatrixafflictorum}
| {\cal O} \setminus {\cal O}_\infty | \leq C(\t) \rho^{ \frac{n}{2} - 1} \g^{2/3} \, .
\ee
\end{theorem}

The proof of Theorem \ref{thm:measure} is similar to that of the
analogous Theorem 4.2 of \cite{BBP1}. 
The specific form  $ \OO_j (\xi) $ in  \eqref{omegaxi} is motivated 
by  application to the DNLW, see \eqref{attila}.
Clearly \eqref{omegaxi}  implies  \eqref{carbonara}.
The asymptotic estimate \eqref{freco} is
the  key point in order to prove \eqref{consolatrixafflictorum}
(in particular for the second order Melnikov conditions at the third line of \eqref{clinea}).
At the end of section \ref{sec:KAM} we explain how
the finitely many  conditions in \eqref{pota} are used  to
estimate the measure 
\begin{equation}\label{sublevels}
|\{ \xi \in {\cal O} \, : \, |\ot^\infty (\xi) \cdot h + \Om^\infty_i (\xi) - \Om^\infty_j (\xi)| < \g \langle h \rangle^{-\tau}  \} | 
\leq \gamma \rho^{\frac{n}{2} -1} \langle h \rangle^{-\tau} , \ h \neq 0 \, , i, j \in \Z \setminus {\cal I}  . 
\end{equation}
This is the main difference with respect to \cite{BBP1}-Lemma 6.1.

\section{Homological equations}\setcounter{equation}{0}\label{sec:hom}

The  integers $ k \in \Z^n $ have indexes in $\mathcal I $   (see \eqref{Isym}), namely $k = (k_h)_{h\in \mathcal I} $.

In the sequel by $ a \lessdot b $ we mean that there exists  $ c  >  0 $ depending
only on $n,\mm, \kappa$ such that $a\leq  c b$.

\begin{definition}{\bf (Normal form vector fields)}\label{nonno}
The normal form vector fields are
\be\label{paderborn}
\mathcal N
:= \partial_\o +\Ngot \ugot \partial_\ugot =
\partial_\o  + \ii  \O z \partial_{z} -\ii  \O {\bar z} \partial_{\bar z}
 =
\o(\xi)\cdot\partial_x +
\ii \sum_{j \in \Z \setminus {\cal I} }  \O_j(\xi) z_j \partial_{z_j}
-\ii
\sum_{j \in \Z \setminus {\cal I}}  \O_j(\xi)  {\bar z}_j \partial_{{\bar z}_j}
\ee
where  the frequencies $ \o_j(\xi), \O_j(\xi) \in \R ,$
$\forall\, \xi \in\mathcal O\subseteq \R^{n/2},$
 are real and symmetric Lipschitz functions
\be\label{Fresim}
\o_{-j} = \o_j  \, , \ \forall j \in {\cal I } \, ,  \quad \ \Om_{-j} = \Om_j  \, , \ \forall j \in \Z \setminus {\cal I} \, ,
\ee
the matrix $\Ngot$ is diagonal
\begin{equation}\label{westfalia}
\Ngot=
\left( \begin{array} {ccc} {\bf 0}_n &0 & 0\\ 0& \ii   \O & 0 \\
0 & 0 & -\ii  \O \end{array}\right)
\,, \ \
 \O:={\rm diag}_{j \in \Z \setminus {\cal I}}( \O_j ) \, ,
\end{equation}
and
there exists $j_* >0$ such that (recall \eqref{carbonara})
 \be\label{Omasin}
\sup_{\xi\in\mathcal O}
\Big| \O_j(\xi) -  \OO_j (\xi) -
 a(\xi) \Big|\lessdot  \frac{\g }{|j|} \, ,
\quad \forall\, |j| \geq j_* \, ,
\ee
(see \eqref{carbonara}) for some
Lipschitz function
$a:\mathcal O\to \R $, independent of $ j $.
 \end{definition}

Note that $ \mathcal N \in {\mathcal R}_{rev}^{\leq 0} $, see Definition \ref{pajata}.
The symmetry condition \eqref{Fresim} implies the resonance relations $ \Om_{-j} - \Om_j = 0 $ and
$ \om \cdot k = 0 $ for all $ k \in \Z^n_{\rm odd} $ defined in \eqref{Zodd}.

\subsection{Symmetrization}\label{sec:symme}

For a  vector field $ X $,
we define its ``symmetrized" ${\cal S} (X) $
by linearity on the monomial vector fields:
\begin{definition}\label{SVF}
The symmetrized monomial vector fields are defined by
\begin{eqnarray}\label{scoglio2}
&& \mathcal S(e^{\ii k\cdot x}\partial_{x_j}) :=
\partial_{x_j}\,,\ \  \mathcal S(e^{\ii k\cdot x}y^i\partial_{y_j}) :=
y^i\partial_{y_j}\,,\  \forall k \in \Z^n_{\rm odd} \, ,\,   |i| = 0,1, \, j \in \mathcal I \, ,  \\
&& \mathcal S(e^{\ii k\cdot x} z_{\pm j} \partial_{z_j}) := z_{j} \partial_{z_j}
\,,\quad
\mathcal S(e^{\ii k\cdot x}  \bar z_{\pm j} \partial_{\bar z_j}) :=  \bar z_{j} \partial_{\bar z_j}\,, \
\forall  k \in \Z^n_{\rm odd}\,,\  j\in\Z\setminus \mathcal I\, ,
\qquad
\label{scoglio3}
\end{eqnarray}
and $ {\cal S} $ is the identity on the other monomial vector fields.
\end{definition}

By \eqref{scoglio2}-\eqref{scoglio3} we write
$ \mathcal SX = X +  X' + X'' $
where
\begin{eqnarray}\label{scoglio4}
X' & := &
{\mathop\sum}_{k \in \Z^n_{\rm odd}, j\in\Z\setminus\mathcal I}
X^{(z_j)}_{k,0,e_{j},0}
(1-e^{\ii k\cdot x} )z_{j} \partial_{z_j}+
X^{(\bar z_j)}_{k,0,0,e_{j}}
(1-e^{\ii k\cdot x})\bar  z_{j} \partial_{\bar z_j}
\end{eqnarray}
and
\begin{eqnarray}\label{scoglio5}
X'' & := &
\sum_{k \in \Z^n_{\rm odd}, k \neq 0, j \in\mathcal I}
X^{(x_j)}_{k,0,0,0} (1-e^{\ii k\cdot x}) \partial_{x_j}+ \sum_{k\in \Z^n_{\rm odd}, k \neq 0, j \in\mathcal I, |i| = 0,1}
X^{(y_j)}_{k,i,0,0}
(1-e^{\ii k\cdot x})y^i \partial_{y_j}
\nonumber
\\
&&
\quad+
\sum_{k\in \Z^n_{\rm odd}, j\in\Z\setminus\mathcal I}
X^{(z_j)}_{k,0,e_{-j},0}
(z_j-e^{\ii k\cdot x} z_{-j}) \partial_{z_j}+
X^{(\bar z_j)}_{k,0,0,e_{-j}}
(\bar z_j-e^{\ii k\cdot x}\bar  z_{-j}) \partial_{\bar z_j}\,.
\end{eqnarray}
The  ``symmetric" subspace $ E $ defined in \eqref{pandistelle} is invariant under the flow evolution
generated by the  vector field $  X $,
because $ X : E \to E $. Moreover the vector fields $ X $ and $ {\cal S} (X) $ coincide on  $ E $: 
\begin{proposition}\label{carbonara1}
 $ X_{|E} = (\mathcal S X)_{|E} $.
\end{proposition}

As a consequence $  v(t) \in E $ is a solution of $\dot v=X(v)$ if and only if
it is a solution of  $\dot v=\mathcal ({\cal S} X)(v) $,
and we  may replace the vector field
$ X $ with its symmetrized $ {\cal S} (X) $ without changing the dynamics on the invariant subspace $ E $.
The following lemma shows that both the $\mathtt a $-momentum and T\"oplitz norms of the
 symmetrized vector field $ {\cal S} (X) $ are controlled by those of $ X $.

\begin{proposition}\label{topotip} For $N_1\geq N_0 $ (defined in \eqref{caracalla})
which satisfy
\begin{equation}\label{farlocco}
N_1 e^{-N_1^b \min\{ s,\aaa\}}\leq 1\,, \quad
b N_1^b \min\{ s,\aaa\}\geq 1\, ,
\end{equation}
the norms of the symmetrized vector field satisfy 
\begin{eqnarray}
i) \, \|\mathcal S X\|_{s,r,\aaa} \leq \|X\|_{s,r,\aaa} \, , \quad  ii) \,
 \|\mathcal S X\|_{s,r,\aaa}^{\rm lip}
\leq \|X\|_{s,r,\aaa}^{\rm lip} \, , \quad iii) \,
\|\mathcal S X\|_{s,r,\aaa,N_1,\theta,\mu}^T
\leq 9 \|X\|_{s,r,\aaa,N_1,\theta,\mu}^T \, .
\label{carbonara4}
\end{eqnarray}
Moreover, if $ X $ is reversible, or real-coefficients,  or real-on-real,
or even,  the same holds for $\mathcal S X $.
\end{proposition}

\begin{pf} In order to prove \eqref{carbonara4}-i)  we first note that
the symmetrized monomial vector fields
$ \partial_{x_h} $, $ y^i\partial_{x_h} $, $ z_{j} \partial_{z_j} $, $ {\bar z}_{j} \partial_{{\bar z}_j} $ in
\eqref{scoglio2}-\eqref{scoglio3}
have zero momentum and are independent of $ x $.
Hence  their contribution to the $ a$-momentum  norm \eqref{normadueA}
is smaller or equal than the contribution of the (not yet symmetrized) monomials
$ e^{\ii k \cdot x}\partial_{x_j} $, $ e^{\ii k \cdot x} y^i\partial_{x_j} $, $ e^{\ii k \cdot x} z_{\pm j} \partial_{z_j} $,
$ e^{\ii k \cdot x} {\bar z}_{\pm j} \partial_{{\bar z}_j} $ of $ X $.  This proves \eqref{carbonara4}-i).
\\[1mm]
{\sc Proof of the \eqref{carbonara4}}-$ii$). The estimate \eqref{carbonara4}  follows by
\be
i) \ \| X' \|_{s,r,\aaa}^T \leq 6  \|X\|^T_{s,r,\aaa} \, , \qquad   \label{resti12}
ii) \  \| X'' \|_{s,r,\aaa}^T \leq 2 \|X\|^T_{s,r,\aaa} \, . 
\ee
{\sc Proof of \eqref{resti12}-$i$)}.
We claim that,  for $ N \geq N_1 $, the projection
$ \Pi_{N,\theta, \mu} X' = \tilde X' + N^{-1}\hat X' $
with
\be\label{calippo}
\tilde X' \in \mathcal T_{s,r,\aaa} \, , \quad
\|\tilde X'\|_{s,r,\aaa}\leq 6 \|X\|^T_{s,r,\aaa} \,,  \quad
\|\hat X'\|_{s,r,\aaa}\leq 5 \|X\|^T_{s,r,\aaa} \, ,
\ee
implying \eqref{resti12} (also because $ \| X' \|_{s,r,\mathtt a} \leq 2\| X \|_{s,r,\mathtt a}  $).
In order to prove \eqref{calippo} we write  the $ (N,\teta, \mu)$-projection as
\be \label{arrosto}
\Pi_{N,\theta ,\mu} X' =
U+U^{-} +U_\bot  + U_\bot^-
\ee
where
\begin{eqnarray}
&&
U := \sum_{k\in \mathcal K_N, |j|>\theta N}
X^{(z_j)}_{k,0,e_{j},0} (1-e^{\ii k\cdot x}) z_{j} \partial_{z_j}\,, \quad
U^{-} :=\sum_{k\in \mathcal K_N, |j|>\theta N}
X^{(\bar z_j)}_{k,0,0,e_{j}} (1-e^{\ii k\cdot x})\bar  z_{j} \partial_{\bar z_j}\,,
\nonumber
\\
&&
U_\bot := \sum_{|j|>\theta N }
\Big(\sum_{k   \in \Z^n_{\rm odd} \setminus \mathcal K_N} X^{(z_j)}_{k,0,e_{j},0} \Big) z_{j} \partial_{z_j}\,,\quad
U_\bot^- := \sum_{ |j|>\theta N} \Big( \sum_{k  \in \Z^n_{\rm odd} \setminus \mathcal K_N }
X^{(\bar z_j)}_{k,0,0,e_{j}} \Big)
\bar  z_{j} \partial_{\bar z_j}\,,
\nonumber
\end{eqnarray}
and
$ \mathcal K_N :=  \big\{ k \in \Z^n_{\rm odd} \, , \ |\pi(k)|, |k| < N^b \big\}  $,
$ \pi(k):=\sum_{j\in\mathcal I} j k_j $.
Then \eqref{resti12} follows by Steps 1)-2) below.
\\[1mm]
{\sc Step $ 1 $)} {\it The projection $ \Pi_{N,\theta, \mu} (U + U^-) =
({\tilde U} + {\tilde U}^-) + N^{-1} ({\hat U} + {\hat U}^- ) $ with
\be\label{calippo1}
{\tilde U} \, , {\tilde U}^- \in \mathcal T_{s,r,\aaa}\,, \ \|\tilde U\|_{s,r,\aaa}, \,
\|{\tilde U}^-\|_{s,r,\aaa} \leq 6 \|X\|^T_{s,r,\aaa} \,,
\quad \| {\hat U} \|_{s,r,\aaa}, \| {\hat U}^- \|_{s,r,\aaa}\leq 6 \|X\|^T_{s,r,\aaa}  \, .
\ee}
Since $ X $ is quasi-T\"oplitz, Lemma \ref{proJ} implies that the projection
\begin{equation}\label{aleph}
\Pi_{\rm diag}\Pi^{(0)}X= \sum_{k\in\Z^n, j\in\Z\setminus\mathcal I}
X^{(z_j)}_{k,0,e_{j},0}
e^{\ii k\cdot x} z_{j} \partial_{z_j}+
\sum_{k\in\Z^n,j\in\Z\setminus\mathcal I}
X^{(\bar z_j)}_{k,0,0,e_{j}}
e^{\ii k\cdot x}\bar  z_{j} \partial_{\bar z_j}
=:W+W'
\end{equation}
is quasi-T\"oplitz as well and ($\|\cdot\|^T_{s,r,\aaa}$ is short for
$\|\cdot\|^T_{s,r,\aaa,N_1,\theta,\mu}$)
$$
\|W\|^T_{s,r,\aaa}, \|W'\|^T_{s,r,\aaa} \leq   \|\Pi_{\rm diag}\Pi^{(0)}X\|^T_{s,r,\aaa} \stackrel{\eqref{fhT}}\leq \| X\|^T_{s,r,\aaa}\,.
$$
By  \eqref{margherita} we have $\Pi_{\rm diag}\Pi^{(0)}\Ta_{s,r,\aaa}\subset \Ta_{s,r,\aaa} $, hence
 Lemma \ref{proJ1} applied to $ W $ implies that for every $N\geq N_1$ there exist ($N$-dependent)
\begin{equation}\label{benten}
\tilde W=
\sum_{|\pi(k)|,|k|<N^b, |j|>\theta N}
\tilde W_k
e^{\ii k\cdot x} z_{j} \partial_{z_j}
\,,\qquad
\hat W=
\sum_{|\pi(k)|,|k|<N^b, |j|>\theta N}
\hat W_{k,j}
e^{\ii k\cdot x} z_{j} \partial_{z_j}
\end{equation}
(note that $\tilde W $ is $(N, \teta, \mu)$-linear and T\"oplitz)
with
\begin{equation}\label{JLB}
\Pi_{N,\theta,\mu} W
=
\sum_{|\pi(k)|,|k|<N^b, |j|>\theta N}
X^{(z_j)}_{k,0,e_{j},0}
e^{\ii k\cdot x} z_{j} \partial_{z_j}
 =\tilde W+ N^{-1}\hat W
\end{equation}
and
$ \|\tilde W\|_{s,r,\aaa}\,,\
\|\hat W\|_{s,r,\aaa}
\leq
\frac32 \| W\|_{s,r,\aaa}^T
\leq
\frac32 \| X\|^T_{s,r,\aaa} $.  By \eqref{arrosto},\eqref{aleph},\eqref{benten} and \eqref{JLB}
we have
$$
U = \! \sum_{k\in \mathcal K_N, |j|>\theta N} \! \tilde W_k (1-e^{\ii k\cdot x}) z_{j} \partial_{z_j}
+
N^{-1} \!\! \! \sum_{k\in \mathcal K_N, |j|>\theta N} \!\!\!  \hat W_{k,j}
(1-e^{\ii k\cdot x}) z_{j} \partial_{z_j}
=:\tilde U+ N^{-1}\hat U\,.
$$
Note that $\tilde U$ is T\"oplitz.
Moreover
$$
\|\hat U\|_{s,r,\aaa}\stackrel{\eqref{normadueA}}\leq
\sup_{\|z\|_{a,p}<r}
\Big\| \Big(
\sum_{k\in \mathcal K_N}
2 e^{\aaa|\pi(k)|} e^{s|k|} |\hat W_{k,j}|
| z_{j}|
\Big)_{|j|>\theta N}\Big\|_{s,r}
\stackrel{\eqref{benten}}\leq
2\|\hat W\|_{s,r,\aaa}
\leq 3\|X\|^T_{s,r,\aaa}\,.
$$
An analogous estimate holds true for $\tilde U$.
A similar decomposition holds for $U^{-} $ in \eqref{arrosto}.
\\[1mm]
{\sc Step $ 2 $)}  $ N \| U_\bot \|_{s,r,\aaa}\,,\ N \| U_\bot^- \|_{s,r,\aaa}\leq  \|X\|_{s,r,\aaa}\,. $
\\[1mm]
We have
\begin{eqnarray*}
\|U_\bot \|_{s,r,\aaa}
&\stackrel{\eqref{normadueA}}=&\sup_{\|z\|_{a,p}<r} \Big\|
\Big( \Big|\sum_{k \in  \Z^n_{\rm odd} \setminus \mathcal K_N} X^{(z_j)}_{k,0,e_{j},0} \Big|
|z_{j}| \Big)_{|j|>\theta N} \Big\|_{s,r}
\\
& \leq &  \sup_{\|z\|_{a,p}<r} \Big\| \Big( e^{-N^b \min\{ s,\aaa\}} \sum_{|\pi(k)| \ or\ |k|\geq N^b}
e^{\aaa|\pi(k)|+s|k|} |X^{(z_j)}_{k,0,e_{j},0}| |z_{j}|  \Big)_{|j|>\theta N} \Big\|_{s,r}
\\
& \leq & e^{-N^b \min\{ s,\aaa\}} \|X\|_{s,r,\aaa}
\stackrel{\eqref{farlocco}}\leq N^{-1}\|X\|_{s,r,\aaa}
\end{eqnarray*}
and similarly for $ U_\bot $.
\\[1mm]
{\sc Proof of \eqref{resti12}-$ii$).} The estimate \eqref{resti12}-ii follows by
\begin{equation}\label{cacioepepe}
\| \Pi_{N,\theta,\mu} X''\|_{s,r,\aaa}
\leq 2 N^{-1} \|X\|_{s,r,\aaa}\,, \quad \forall\, N\geq N_0\,.
\end{equation}
In order to prove \eqref{cacioepepe} we note that  the momentum of
$ e^{\ii k \cdot x} z_{-j} \partial_{z_j} $ with $|k|<N^b,$ $|j|>\theta N,$ $N\geq N_1\geq N_0$,  satisfies
\begin{equation}\label{papapig}
|\pi(k,e_{-j},0;z_j)|= \Big| {\mathop\sum}_{h\in\mathcal I} h k_h - 2 j \Big|
\geq 2|j|- \kappa |k|
\geq 2\theta N -\kappa N^b
\stackrel{\eqref{caracalla}}>
N > N^b
\end{equation}
(where $\kappa:=\max_{h\in\mathcal I}|h|$, recall \eqref{caracalla}).
Then by \eqref{scoglio5} and \eqref{perla} the projection
$ \Pi_{N,\theta,\mu}X''=V+V'$  with
$$
V:=
\sum_{|j|>\theta N}
\big(\sum_{k\in \mathcal K_N}
X^{(z_j)}_{k,0,e_{-j},0}\big)
 z_j \partial_{z_j}
\,,\qquad
V':=
\sum_{|j|>\theta N}
\big(\sum_{k\in \mathcal K_N}
X^{(\bar z_j)}_{k,0,0,e_{-j}}\big)
\bar z_j \partial_{\bar z_j}\,.
$$
We have
\begin{eqnarray}
\| V\|_{s,r,\aaa}
&\stackrel{\eqref{normadueA}}=&
\sup_{\|z\|_{a,p}<r} \Big\|
\Big(
 \big| \sum_{k\in \mathcal K_N}
X^{(z_j)}_{k,0,e_{-j},0} \big|
 |z_j|
\Big)_{|j|>\theta N}
\Big\|_{s,r}
\label{peppapig1}
\\
%
&\stackrel{\eqref{papapig}}\leq&
\sup_{\|z\|_{a,p}<r} \Big\|
\Big(
  \sum_{k\in \mathcal K_N}
  e^{-\aaa N}
e^{\aaa|\pi(k,e_{-j},0;z_j)|}
|X^{(z_j)}_{k,0,e_{-j},0} |
 |z_{-j}|
\Big)_{|j|>\theta N}
\Big\|_{s,r}
\nonumber
\\
&\leq&
e^{-\aaa N} \|X\|_{s,r,\aaa} \stackrel{\eqref{farlocco}}\leq N^{-1} \|X\|_{s,r,\aaa}
\nonumber
\end{eqnarray}
where in \eqref{peppapig1} we have used that
 the domain $\{ \|z\|_{a,p}<r\}$ is invariant
under the map $ z_j \mapsto z_{-j} $.
Since a similar estimate holds for $V'$,
\eqref{cacioepepe} follows.

Finally, the vector field $ \mathcal S X $ is even because $ {\cal S} X_{|E} = X_{|E}  $
(Proposition \ref{carbonara1}) and $ X $ is even. Since $ X $ is real-coefficients,
Definition \ref{SVF} immediately  implies that
$ \mathcal S X $ is real-coefficients.
Since $ X $ is reversible and real-on-real,  \eqref{tango} and \eqref{reality2} enable to
check that
$ X' , X'' $ in \eqref{scoglio4}-\eqref{scoglio5} are reversible  and real-on-real, and so
$ \mathcal S X $.
\end{pf}

\begin{remark}\label{simme}
The assumptions $X \in\mathcal R_{rev},$
$Y\in\mathcal R_{a-rev},$ $X=\mathcal S X,$
$Y=\mathcal S Y$ are not sufficient to imply
$[X,Y]=\mathcal S [X,Y]$, as the example
$X=\ii (z_{-1}\partial_{z_2}+z_{1}\partial_{z_{-2}}
-\bar z_{1}\partial_{\bar z_{-2}}- \bar z_{-1}\partial_{\bar z_2})$,
$Y=z_{2}\partial_{z_1}+z_{-2}\partial_{z_{-1}}
+\bar z_{-2}\partial_{\bar z_{-1}}+ \bar z_{2}\partial_{\bar z_1} $ shows.
\end{remark}

\subsection{Homological  equations and quasi-T\"oplitz property}

We consider the homological equation
\begin{equation}\label{tebe}
{\rm ad}_{\mathcal N} F = R - [R] \,
\end{equation}
where
\begin{equation}\label{zenone}
R\in  \mathcal{R}_{rev}^{\leq 0} \
({\rm see} \,
{\rm Definition} \, \ref{pajata}),
\ \ \ R = \mathcal S R \ ({\rm see} \, {\rm Definition} \, \ref{SVF})
\end{equation}
and
\be\label{[R]} [R] := \langle R^x\rangle\partial_x + 
{\mathop\sum}_{j \in \Z\setminus \mathcal I}
 \langle  R^{z_j z_j } \rangle z_j \partial_{z_j}+
\langle   R^{{\bar z}_j {\bar z}_j }
 \rangle \bar z_j \partial_{\bar z_j}\,,
\ee
where $ \langle \cdot \rangle $ denotes the average with respect to the angles $ x.$
By Lemmata \ref{lem:closed} and \ref{embeh2}
and since $ {\cal N} \in \mathcal R_{rev}^{\leq 0} $
(see Definition \ref{nonno}), the action
$ {\rm ad}_{\mathcal N} :\mathcal{R}_{a\mbox{-}rev}^{\leq 0} \to
\mathcal{R}_{rev}^{\leq 0} $.
The commutator
\be\label{lorella}
{\rm ad}_{\mathcal N} F = {[F,{\mathcal N} ]} =
\begin{cases}
\big( \partial_\o F^\ugot -\Ngot F^\ugot \big)\partial_\ugot \qquad \qquad \qquad \qquad \qquad \ {\rm if} \ F = F^{(-1)}
\cr
\partial_\o F^x \partial_x + \big(
\partial_\o F^{\ugot,\ugot} + [F^{\ugot,\ugot} ,\Ngot]
\big) \ugot \partial_\ugot   \qquad \quad {\rm if} \ F = F^{(0)}
\end{cases}
\ee
(recall the notations in  \eqref{ndg}-\eqref{baronerosso})
where $[F^{\ugot,\ugot} ,\Ngot]=F^{\ugot,\ugot} \Ngot -\Ngot
F^{\ugot,\ugot} $ is the usual commutator between matrices (and $
\Ngot $ is defined in \eqref{westfalia}). We solve \eqref{tebe} when
\be\label{3PRO}
R=R^{(h)}_K:= \Pi_{|k|<K} \Pi_{|\pi|<K} R^{(h)} \, , \quad h = 0,-1 \, , \quad K \in \N
\ee
(recall the projections \eqref{minas}, \eqref{tirith} and \eqref{phom}).

\begin{definition}\label{cipolla}{\bf (Melnikov conditions)}
Let $\g>0.$  The frequencies
$\o(\xi)=(\ot(\xi),\ot(\xi)),$ $\ot\in\R^{n/2},$ $
\O(\xi)$ satisfy the Melnikov conditions (up to $ K > 0 $)  at 
$\xi\in\R^{n/2},$
if: $\forall\, h \in \Z^{n/2}$, $ | h |  < K $, $ i, j \in \Z \setminus {\cal I} $,
\begin{eqnarray}
\label{delfi1q}
     |\ot(\xi)\cdot h|
     & \geq &
     \g  \langle h \rangle^{-\t}  \quad
     {\rm if} \quad  h\neq 0 \,,
     \\
\label{delfi2q} | \ot(\xi)\cdot h  +  \O_j |
&\geq&
\g \langle h \rangle^{-\t}
\,,
\\
\label{delfi4q}
|\ot(\xi)\cdot h  + \O_i(\xi)  +
\O_j(\xi)  |
&\geq&
 \g \langle h \rangle^{-\t}
 \,, \\
\label{delfi3q} | \ot(\xi)\cdot h  -  \O_i(\xi)  +  \O_j(\xi) |
& \geq &
\g
\langle h \rangle^{-\t}
  \quad    {\rm if}\quad  h\neq 0 \ \ \
     {\rm or}\ \  i\neq \pm j\,,
\end{eqnarray}
where $\langle h \rangle:=\max\{ |h|,1\} $ and
$ \t > 1 / b $.
\end{definition}
For $k\in \Z^n$ we set $k_\pm:=(k_j)_{j\in\mathcal I^\pm}\in\Z^{n/2},$
namely $k=(k_+,k_-).$
Then
\begin{equation}\label{pesca}
    \o\cdot k=\ot\cdot h\,,\ \ {\rm with}\ \
    h:=k_+ + k_-\in\Z^{n/2} \ \ {\rm and}\ \
    k\notin \Z^n_{\rm odd}\
    \stackrel{\eqref{Zodd}}\Longrightarrow \ h\neq 0\,.
\end{equation}
Note that $|h|\leq |k_+|+|k_-|=|k|.$

\begin{lemma} {\bf (Solution of homological equations)} \label{lem:HOMO}
Let
 $ s, r, \mathtt a>0,$  $ K >  0 .$
 Let $\mathcal O\subset \R^{n/2}$ and assume that
 the Melnikov conditions \eqref{delfi1q}-\eqref{delfi3q}
 are satisfied $\forall\,\xi\in\mathcal O.$
 Then, $\forall\,\xi\in\mathcal O,$
 the homological equation \eqref{tebe}
 with $ R =R(\cdot;\xi)$ as in \eqref{zenone},\eqref{3PRO}
has a unique  solution
$F=F(\cdot;\xi)$
$$
F\in \mathcal R_{a\mbox{-}rev}^{\leq 0} \, , \qquad F=\mathcal S F \, , \qquad
F=\Pi_{|k|<K}\Pi_{|\pi|<K} F
$$
with $ \langle F^y \rangle=0,$ $\langle F^{y,y} \rangle =0$, $
\langle F^{z_i^\pm, z_i^\pm} \rangle=0$. It satisfies
\be\label{la1}
\| F\|_{s,r,\aaa,\mathcal O}
 \leq
 \g^{-1}K^\tau \| R \|_{s,r,\aaa,\mathcal O}
\ee
\be\label{corviale}
\| F\|_{s,r,\aaa,\mathcal O}^{\rm lip}
 \lessdot
 \g^{-1}K^\tau \| R \|_{s,r,\aaa,\mathcal O}^{\rm lip}
 +
\g^{-2} K^{2\t+1}\big(
|\o|^{\rm lip}_{\mathcal O}+
|\O|^{\rm lip}_{\mathcal O}
\big) \| R \|_{s,r,\aaa,\mathcal O}
\,.
\ee
\end{lemma}

\begin{pf}
By \eqref{lorella} the homological equation
\eqref{tebe} splits into
\begin{equation}\label{tebebis}
\partial_\o F^\ugot -\Ngot F^\ugot = R^\ugot  \,,
\quad
\partial_\o F^x = R^x - \langle R^x\rangle\,,\quad
\partial_\o F^{\ugot,\ugot} + [F^{\ugot,\ugot} ,\Ngot] = R^{\ugot,\ugot} - [R]^{\ugot,\ugot}\,.
\end{equation}
Since $ R = \mathcal S R $ (recall \eqref{zenone}), by
\eqref{scoglio2} we get
\begin{equation}\label{parmenide}
R^x(x)= \langle R^x \rangle + {\mathop\sum}_{k\notin\Z^n_{\rm odd} }R^x_k
e^{\ii k\cdot x} \, , \quad {\rm similarly\ for}\quad R^y(x) \,,\
R^{y,y}(x) \, .
\end{equation}
Since $ R $ is reversible and even the average
\be\label{venerdi1}
\langle R^y \rangle = 0 \, ,\quad \langle R^{y,y} \rangle = 0
\ee
By \eqref{westfalia}  the
first equation in \eqref{tebebis} amounts to
$\partial_\o F^y  = R^y $,  $ \partial_\o F^z - \ii {\O} F^z = R^z $,
$ \partial_\o F^{\bar z} + \ii { \O}  F^{\bar z} = R^{\bar z} $.
By \eqref{westfalia},  the third equation in \eqref{tebebis} splits into
$    \partial_\o F^{y,y}  = R^{y,y} $,
$ \partial_\o F^{y,z}+\ii F^{y,z}{ \O} = R^{y,z} $
(and the analogous equations   for $ F^{y,\bar z}, $ $ F^{z,y}
$, $ F^{\bar z,y} $),
$ \partial_\o F^{z,\bar z}-\ii F^{z,\bar z} { \O}
-\ii { \O} F^{z,\bar z} =  R^{z,\bar z} $
(analoguosly  for $ F^{\bar z, z} $),
\begin{equation}\label{WFRbis}
   \partial_\o F^{z,z} +\ii F^{z, z} { \O} -
 \ii { \O} F^{z, z} 
= R^{z,z} - [R]^{z,z}
\end{equation}
(analogously for $ F^{\bar z,\bar z} $).
By \eqref{delfi1q}, \eqref{parmenide}, \eqref{venerdi1} and \eqref{pesca}
the equations for $ F^x, F^y, F^{yy}  $ are uniquely (having zero average)
solved, i.e., $F^x (x) =
{\mathop \sum}_{k\notin\Z^n_{\rm odd}} F^x_k e^{\ii k \cdot x}$ with
$F^x_k :=  -\ii{R^x_k}/{ \om \cdot k}$.
Similarly the equations for
$F^{z^\s},$ $ F^{y, z^\s}, $ $ F^{z^\s,y}
$, $\s=\pm$  and $F^{z,\bar z},$ $F^{\bar z,z}$
are solved by \eqref{delfi2q} and \eqref{delfi4q} respectively.

For $i,j\in\Z\setminus\mathcal I,$ developing in Fourier series $ F^{z_i z_j}(x)= {\mathop\sum}_{k \in
\Z^n}
 F^{z_i z_j}_k e^{\ii k\cdot x} $, equation \eqref{WFRbis} becomes
\be\label{WFRter}
\ii(\o\cdot k  +  \O_j -   \O_i) F^{z_i z_j}_k = R^{z_i z_j}_k - [R]^{z_i z_j}_k \, .
\ee
If $ i \neq \pm j $ then
\eqref{WFRter} is easily solved by
\eqref{delfi3q}. Otherwise, since
$R=\mathcal S R$ and by
\eqref{scoglio3},
\be\label{rimmel}
{\rm if} \ i = j \  \Longrightarrow \
R^{z_i z_i}_k=0\,,\, \forall\, k\in \Z^n_{\rm odd}\setminus\{ 0\}; \quad
{\rm if} \ i = - j \ , \, (i \neq 0)\  \Longrightarrow \
R^{z_i z_{-i}}_k=0\,,\  \forall\, k \in \Z^n_{\rm odd}\, .
\ee
 Then
 \eqref{WFRter} is solved by \eqref{delfi3q}  and \eqref{pesca}.

 The properties of anti-reversibility,  anti-real-coefficients, real-on-real, and parity for the vector field  solution
 $ F $ are easily verified.
  The estimates \eqref{la1}-\eqref{corviale} directly follow by
  bounds on the small divisors in the Melnikov
 conditions  \eqref{delfi1q}-\eqref{delfi3q}
 (and \eqref{pesca}) and the expression of $ F $.
\end{pf}

%



The  solution of the homological equation  is quasi-T\"oplitz.

\begin{proposition}\label{FP}
{\bf (Quasi-T\"oplitz)}
Let the normal form $ {\mathcal N} $ be as in Definition \ref{nonno} 
 and assume that $ R \in {\mathcal Q}^T_{s,r,\aaa} (N_0,  \teta, \mu ) $.
Let $ F $ be the (unique) solution of the homological equation
\eqref{tebe} found in Lemma  \ref{lem:HOMO}, for all  $ \xi \in\mathcal O $ satisfying
the Melnikov conditions \eqref{delfi1q}-\eqref{delfi3q}.
If, 
in addition, 
\be\label{s+s-}
|\ot (\xi)
\cdot h + p | \geq \g^{2/3} \langle h \rangle^{-\tau}  \, , \ \ \forall  |h|\leq K, \ p \in \Z \, ,
\ \ (h,p)\neq(0,0)\,,
\ee
then  $ F=F(\cdot;\xi) \in {\mathcal Q}^T_{s,r,\aaa} (N_0^* , \teta, \mu) $ with
\be\label{N0large}
N_0^*:= \max\big\{N_0 \, , j_*, \ \hat{c}\g^{-1/3} K^{\tau +1} \big\}
\ee
for a (suitably large) constant $\hat{c} := \hat{c} (\mm , \kappa ) \geq 1 $.
Moreover
\be\label{Fijstimeh}
\| F(\cdot;\xi) \|^T_{s,r,\aaa, N_0^*,  \teta, \mu} \leq 4\hat{c} \g^{-1} K^{2 \t}
\| R (\cdot;\xi)\|_{s,r,\aaa, N_0,  \teta, \mu}^T \, .
\ee
\end{proposition}

\begin{pf} The proof follows step by step the one of the analogous Proposition 5.1 of \cite{BBP1}. \end{pf}

\section{Proof of  Theorem \ref{thm:IBKAM} }\label{sec:KAM}
\setcounter{equation}{0}
\subsection{First step}\label{volkswagen}

We perform a preliminary change of variables in order to improve
the smallness conditions of the perturbation. In particular we
want to average out the term $\PP^y (x)\partial_y$ defined in
\eqref{primopezzo}. We  introduce the symmetrized vector fields
(see Definition \ref{SVF})
\begin{equation}\label{aurelia}
R^y (x)\partial_y:=\mathcal S \PP^y (x)\partial_y\,, \qquad
R:=\mathcal S \PP\,, \qquad X:=\mathcal S \mathcal X=\mathcal N+R
\end{equation}
(since $\mathcal S \mathcal N=\mathcal N$). 
By assumption (A3) and the last statement of Proposition
\ref{topotip},  $ R \in {\cal R}_{rev}$ (see Definition
\ref{pajata}). Moreover Proposition \ref{carbonara1} implies that
$ X_{|E}=\mathcal X_{|E} $.

Next we  study the homological equation \be\label{F0P0} - {\rm
ad}_{\mathcal N}  F + \Pi_{|k|<\g^{-1/(7n)}} R^y \partial_y =
\langle R^y \rangle\partial_y \stackrel{\eqref{venerdi1}} =0 \ee
because $ R $  is reversible and even.

\begin{lemma}\label{FST}
For all $ \xi $ in
$ {\cal O}_* := \big\{   \xi \in {\cal O} \ : \  | \ot
 (\xi) \cdot h | \geq \g^{2/3} \langle h \rangle^{- n/2} \, , \, \forall 0 \, < |h| < \g^{-1/(7n)} \big\} $
the homological equation \eqref{F0P0} admits a unique solution
with
$ \langle F \rangle = 0 $ which  satisfies \be\label{sporta} \|
F\|^T_{3s/4,r,\aaa,N_0,\theta,\mu,\lambda,\mathcal O_*} =
 \| F  \|_{3s/4,r,\aaa,\mathcal O_*}^\l \leq C( s) \e\,.
\ee
Moreover   $ F \in {\cal R}_{a-rev}^{\leq 0} $
and $\mathcal S F = F$.
\end{lemma}

We now apply Proposition \ref{main2}  with $\vec p \rightsquigarrow
(3s/4, r, \aaa, N_0,  \teta ,\mu,\l , \mathcal O_*) $ and $\vec
p\,' \rightsquigarrow \vec p_0$ with
$ \vec p_0 := $ $ (s/2, r/2,\aaa/2, N_0^{(0)},
4\teta/3, 3\mu/4,\l,\mathcal O_*) $
where
 $ N_0^{(0)} \geq \max\{ N_0,\bar N\} $
 (recall \eqref{stoppa})
 is chosen large enough
 so that
 \eqref{stoppa}-\eqref{Giobbe}-\eqref{Giobbe1}  are satisfied and \eqref{sporta} imply
 condition \eqref{piccof2} for $ \e $ sufficiently small.
 Let $ \bar \Phi $ be the time $ 1$-flow of  $ F $ (so that $e^{{\rm ad}_{F }}= \bar \Phi_{\star}$).
 Since the quasi-T\"oplitz norm is non-increasing with $ N_0 $ (see \eqref{dino})
 we may also take $ N_0 \geq \bar N $ large enough so that
\eqref{farlocco} (with $ N_0  \rightsquigarrow N_1  $) 
holds. Hence
\begin{eqnarray}\label{zione}
\| e^{{\rm ad}_{F} }
 (R- R^y  \partial_y) \|^T_{\vec p_0}
&  \stackrel{\eqref{gPhifbis} } \leq &
 2  \| R -  R^y  \partial_y \|_{s, r,\aaa, N_0, \teta,\mu,\l,\mathcal O_*}^{T} \nonumber \\
 &  \stackrel{\eqref{aurelia}, \eqref{carbonara4}}  \leq &
 18  \| {\cal P} - {\cal P}^y(x) \partial_y  \|_{s, r,\aaa, N_0, \teta,\mu,\l,\mathcal O_*}^{T}
 \stackrel{\eqref{primopezzo}, \eqref{KAMcondition}} < 18 \g  \, .
\end{eqnarray}
 Similarly \eqref{gPhif12bis} (with $h\rightsquigarrow 1$,
$b_j\rightsquigarrow 1/j!$) implies,
for $h=-1,0,$
\begin{equation}\label{zionebis}
\Big\|\Big( e^{{\rm ad}_{F} }
 (R- R^y  \partial_y) - (R- R^y  \partial_y)\Big)^{(h)}\Big\|^T_{\vec p_0}
 \lessdot
 \|
 \PP_*  \|_{s, r,\aaa, N_0, \teta,\mu,\l,\mathcal O_*}^{T}
    \| F\|^T_{3s/4,r,\aaa,N_0,\theta,\mu,\lambda,\mathcal O_*}\!\!
    \stackrel{\eqref{sporta},\eqref{KAMcondition}}\leq
    \!\!
    C(s) \g\e \,.
\end{equation}
Since the commutator $ [ F, R^y(x) \partial_y ] = [ F^y(x)
\partial_y, R^y(x) \partial_y ] = 0 $ we deduce $ e^{{\rm ad}_{F
}} (R^y \partial_y) = R^y  \partial_y $, and, using also
\eqref{F0P0}, we get $ e^{{\rm ad}_{F}} {\cal N} = {\cal N} + {\rm
ad}_{F}  {\cal N} $. Hence, using \eqref{F0P0},
\be
 e^{{\rm ad}_{F }} X =
 {\cal N} + \Pi_{|k|\geq \g^{-1/(7n)}} R^y \partial_y+
e^{{\rm ad}_{F }}  (R-R^y \partial_y) =: {\cal N}_0 + P_0
\label{accazero}
\ee
where $ \mathcal N_0:=\mathcal N $. Then we  consider the
symmetrized vector field
\begin{equation}\label{flaminia}
X_0  :=  \mathcal S (e^{{\rm ad}_{F }} X) ={\cal N}_0 +
R_0\,,\qquad R_0:=\mathcal S P_0 \,.
\end{equation}
Since $R^y (x) \partial_y $ depends on the variable $ x $ only we
have
\be\label{RCA} \|\mathcal S  \Pi_{|k|\geq \g^{-1/(7n)}} R^y(x)
\partial_y \|^T_{\vec p_0} = \| \mathcal S  \Pi_{|k|\geq
\g^{-1/(7n)}} R^y(x) \partial_y \|^\l_{ s/2, r, \mathtt a, {\cal
O}} \lessdot \g \e \, , \ee
arguing as for \eqref{sporta},  using  \eqref{carbonara4}, \eqref{smoothl}, 
and for $\g<\g_*$ small (depending on $ s $ and $ n $).
Recollecting  \eqref{flaminia}, \eqref{accazero}, \eqref{zione},
\eqref{RCA} and \eqref{zionebis} we get

\begin{lemma} The constants
$ \bar\varepsilon_0 :=  \varepsilon_0^{(-1)}+ \varepsilon_0^{(0)} $,
$ \varepsilon_0^{(h)} := \g^{-1} \| R_0^{(h)} \|^T_{\vec p_0} $, $ h=-1,0 $,
$ \Theta_0 := \g^{-1} \| R_0 \|^T_{\vec p_0} $
satisfy
$ \e_0^{(h)}\leq C(s,n) \e $, $ h =-1,0 $, $ \Theta_0\leq 2^8$,
where $ \e $ is defined in \eqref{KAMcondition}.
\end{lemma}

The vector fields $ P_0, R_0 \in {\mathcal R}_{rev}$ because $ F
\in {\cal R}_{a-rev}$ (Lemma \ref{FST}), $ R \in {\cal R}_{rev} $,
and 
using Proposition \ref{topotip}.
Similarly, since $ \XX \in {\cal R}_{rev} $ (by the hypothesis of
Theorem \ref{thm:IBKAM}) the vector field
\begin{equation}\label{salaria}
\mathcal X_0:= e^{{\rm ad}_F} \mathcal X= {\bar \Phi}_\star \XX
\in {\cal R}_{rev}\,.
\end{equation}
Proposition \ref{carbonara1} implies that $ X_{|E}= ({\cal S}
{\mathcal X})_{|E} = \mathcal X_{|E} $ (see \eqref{aurelia}) and $
{X_0}_{|E} = (e^{\rm ad_F} X)_{|E} $ (see \eqref{flaminia}).
Moreover, since $ F $ is even,
Lemma \ref{montalto} (applied with $Y \rightsquigarrow  F $) and \eqref{salaria}  
imply
\begin{equation}\label{tiburtina}
{\mathcal {X}_{0}}_{|E}={X_{0}}_{|E}\,.
\end{equation}

\subsection{ The KAM step}

We now describe the iterative scheme which produces a sequence of
quasi-T\"oplitz vector fields $X_\nu$  with parameters $\vec
p_\nu= (s_\nu,r_\nu,\aaa_\nu,N^{(\nu)}_0,\theta_\nu,\mu_\nu,\l
,\mathcal O_\nu) $, $\lambda = \gamma/M_0$, and such that
${X_\nu^{\leq 0 }}_{\vert E}$ tends to zero as $ \nu \to + \infty
$. For compactness of notation we drop the index $\nu$ and write
"$+$" for $\nu+1$.
 \\[1mm]
 \noindent{\bf Iterative hypotheses.}
Suppose
 $  \cc < \teta,\mu < \CC   $, $ N_0 \geq \bar N $
 (defined in \eqref{stoppa}),
 $\mathcal O\subseteq\R^{n/2}$.
Let $X=\mathcal N+R$, where
 $\mathcal N$ is  a normal form vector field
 (see Definition \ref{nonno})
 with  Lipschitz frequencies $\o(\xi),\O(\xi)$,
 $\xi\in\R^{n/2}$
and \eqref{Omasin} holds  with some $a(\xi)$, $\forall\, |j|\geq
\CC N_0$ (namely $j_*=\CC N_0$). Moreover $|\o|^{\rm
lip}_{\R^{n/2}}, |\O|^{\rm lip}_{\R^{n/2}}\leq M\leq 2 M_0.$
 The perturbation $R$ satisfies
 $\|R\|^T_{\vec p}<\infty$,
    $R\in \mathcal R_{rev} $,   $\mathcal S R= R$.
 We finally fix some $K$ and we assume that
 $\CC N_0\geq \hat{c}\g^{-1/3}K^{\t +1}$
 (where $\hat{c}$ is the constant introduced in
 \eqref{N0large}).
\smallskip

We now describe a {\em KAM step}, namely a change of variables
generated by the time-1 flow of a vector field $F$ and such that
\begin{equation}\label{cecio}
    X_+:=\mathcal S e^{{\rm ad}_F}X=:\mathcal S \Phi_\star X=
    \mathcal N_++ R_+
\end{equation}
 still
satisfies the iterative hypotheses, with  slightly different
parameters, and a much smaller new perturbation $ R_+ $, see
\eqref{ultimax+}.

\smallskip

\noindent {\bf The new normal form $ \mathcal N_+ $.} Set (recall
\eqref{ndg})
\begin{equation}\label{montecalvo}
R^{\leq 0}_K:= \Pi_{|k|<K}\Pi_{|\pi|<K} R^{\leq 0} =
\Pi_{|k|<K}\Pi_{|\pi|<K} R^{(-1)} +\Pi_{|k|<K}\Pi_{|\pi|<K}
R^{(0)} =:  R_{K}^{(-1)}+ R_{K}^{(0)}\,.
\end{equation}
Since $ R \in {\cal R}_{rev}$ then $ R^{\leq 0}_K \in \mathcal
R_{rev}^{\leq 0} $ and $ \mathcal S R^{\leq 0}_K= R^{\leq 0}_K. $
The new normal form is defined for $ \xi \in \mathcal O $ as
\be\label{newNF}
{\cal N}^+  := {\cal N} + \hat {\cal N} \, ,
\ee
\be\label{ennepiu}
\hat {\cal N}   \stackrel{\eqref{[R]}} {:=}   [R^{\leq 0}_K] =
\langle R^x\rangle\partial_x + \sum_{j \in \Z\setminus \mathcal I}
 \langle  R^{z_j z_j } \rangle z_j \partial_{z_j}+
\langle   R^{{\bar z}_j {\bar z}_j }
 \rangle \bar z_j \partial_{\bar z_j} =
 \hat \o \cdot \partial_x + \ii \sum_{j\in \Z\setminus \mathcal I} \hat \O_j z_j (\partial_{z_j} -  \bar z_j\partial_{\bar z_j})
\ee because, since $ R^{\leq 0}_K $ is real-coefficients
and real-on-real (Definition \ref{rean})
\begin{equation}\label{triciclo}
\langle R^{z_jz_j} \rangle = \ii {\hat \Omega}_j \, , \  {\hat
\Omega}_j \in \R\,,
 \ \
   \langle R^{\bar z_j\bar z_j} \rangle \stackrel{\eqref{reality2}} =  - \ii {\hat \Omega}_j
\,, \ \forall\,j\in\Z\setminus\mathcal I\,, \ \ \hat \o_j:=
\langle R^{x_j}\rangle \in\R \, , \forall\,j\in\mathcal I \, .
 \end{equation}
Moreover, since  $R $ is even,
 $\hat \o\,, \hat \O$ satisfy \eqref{Fresim},
 namely
$ \hat \o_j \stackrel{\eqref{parity}}= \hat \o_{-j} $,
$ \hat \O_j \stackrel{\eqref{parity}}= \hat \O_{-j} $.
Note that $\hat {\cal N}$ only depends on $R^{(0)}$ and that $\hat {\cal N}-\langle R^x\rangle\partial_x= \Pi_{\rm diag} R$.

The following lemma on the asymptotic of the frequencies is based on the projection Lemma \ref{proJ} for
$ \Pi_{\rm diag} $ similarly to Lemma 5.2 of \cite{BBP1}.
\begin{lemma}\label{frodo} It results
$ \sup_{\xi\in \mathcal O}|\hat \o|, | \hat
\O|_\infty \leq 2 \| R^{ (0)} \|_{s,r,\aaa}$,
$ |\hat \o|^{\rm lip}_{\mathcal O}, | \hat \O|_{\infty, \mathcal O}^{\rm
lip} \leq 2 \| R^{ (0)} \|_{s,r,\aaa}^{\rm lip}  $
and there
exist $ \hat a :\mathcal O\to \R $ satisfying
$ \sup_{\xi\in\mathcal O}|\hat a(\xi) |  \leq 2 \| R^{ (0)}
\|^T_{s,r,\aaa,N_0,\theta,\mu} $
such that
$$
\sup_{\xi\in\mathcal O} |\hat\O_j(\xi)-\hat a(\xi)| \leq
\frac{40}{|j|} \| R^{ (0)} \|^T_{s,r,\aaa,N_0,\theta,\mu}  \, ,
\quad
\forall\, |j| \geq \CC(N_0+1) \, .  
$$
\end{lemma}

\noindent {\bf The new vector field $ X_+ $.} We decompose
$$
X = \mathcal N + R = \mathcal N + R^{\leq 0}_K + (R - R^{\leq 0}_K )
$$
where $ R^{\leq 0}_K $  is defined in \eqref{montecalvo}. We apply
Lemma \ref{lem:HOMO} and Proposition \ref{FP} with $\mathcal O
\rightsquigarrow \mathcal O_+ := \big\{ \xi\in\mathcal O\; | \;  \rm{
\eqref{delfi1q}-\eqref{delfi3q}\;
 and \; \eqref{s+s-} \; hold} \big\} $.
 Let $ F =F^{\leq 0}_K = F^{(-1)}_K+F^{(0)}_K$ $ \in \mathcal
R_{a-rev}^{\leq 0} $ be the unique solution of the homological
equation
\begin{equation}\label{capocotta}
{\rm ad}_{\mathcal N}F = R^{\leq 0}_K-[R^{\leq 0}_K]\,.
\end{equation}
The bounds  \eqref{corviale},  $ | \o |^{\rm lip}, |\O |^{\rm lip}
\leq M \leq 2 M_0 $, and \eqref{Fijstimeh} (with $R
\rightsquigarrow R^{(h)}_K$, $h=-1,0$) imply \be\label{effe}
\|F^{(h)}\|^T_{\vec p_\star} \lessdot \gamma^{-1}
K^{2\tau+1}\|R^{(h)}\|^T_{\vec p}\,,\quad h=-1,0\,, \quad {\rm
where}\quad \vec p_\star:= (s,r,\aaa,6 N_0,\theta,\mu, \l,\mathcal
O_+)\, . \ee
Note that
in \eqref{N0large}-\eqref{Fijstimeh}
$ N_0^* = 6 N_0 $
because, by the iterative hypothesis,  $ j_* = \CC N_0\geq
\hat{c}\g^{-1/3}K^{\t +1}$.

We introduce the new parameters
\begin{equation}\label{magliana}
\vec p_+:=(s_+,r_+,\aaa_+,N_0^+,\theta_+,\mu_+, \l,\mathcal
O_+)\,,
\end{equation}
where $ s/2 \leq s_+ < s $, $ r/2\leq r_+ < r $,  $0<\aaa_+<\aaa,$
$N_0^+\geq 7 N_0,$ $\theta_+ > \theta, $ $ \mu_+ < \mu $,
such that 
\begin{equation}\label{Giobbe+}
 (\kappa+1)  (N_0^+)^{b-L} \ln N_0^+  \leq \mu - \mu_+, \quad
 (7  + \kappa) (N_0^+)^{L-1} \ln N_0^+ \leq   \teta_+  -  \teta \, ,
 \end{equation}
 \be\label{lap1}
 2(N_0^+)^{-b}\ln^2 N_0^+ \leq b\min \{ s-s_+,
 \aaa - \aaa_+ \} \, ,
\ee and note that $ N_0^+ \geq \bar N $ defined in \eqref{stoppa}
(by the iterative hypothesis $N_0\geq \bar N $). If, moreover, the
smallness condition \be\label{piccoF} \| F \|^T_{\vec p_\star}\leq
c(n) \, \d_+ \, , \quad \d_+ := \min \Big\{  1 - \frac{s_+}{s}, 1
- \frac{r_+}{r}  \Big\} \ee holds (see \eqref{piccof2}), then
Proposition \ref{main2} (with $\vec p\rightsquigarrow \vec p_\star,$
$\vec p \, '\rightsquigarrow \vec p_+$, $\d \rightsquigarrow
\d_+$) implies that the time $ 1 $-flow generated by $ F $  maps $
D(s_+,r_+) $ into $ D(s,r) $. The transformed and symmetrized
vector field is
\begin{equation}\label{ixspiu}
X^+ := {\cal S} e^{{\rm ad}_F} X
 \stackrel{\eqref{trasfoXY}} =
{\cal S} \Big( X + {\rm ad}_F (X) + {\mathop\sum}_{j \geq 2}
\frac{1}{j!}{\rm ad}_F^j (X) \Big) = {\cal N}^+ + R^+
\end{equation}
with the new normal form $ {\cal N}^+ $ defined in \eqref{ennepiu}
and, by  \eqref{capocotta}, the  new perturbation \be\label{newNP}
R^+ := {\cal S } \Big( R - R^{\leq 0}_K + {\rm ad}_{F} (R^{\leq
0}) + {\rm ad}_{F} (R^{\geq 1})  + {\mathop\sum}_{j \geq 2}
\frac{1}{j!}{\rm ad}_F^j (X)\Big) \ee
where 
$  R^{\geq 1} := {\mathop\sum}_{j \geq 1} R^{(j)} $, see \eqref{phom}, so
that $ R = R^{\leq 0} + R^{\geq 1} $.

We set \be\label{xhx} \varepsilon^{(h)} := \g^{-1} \| R^{(h)}
\|^T_{\vec p} \, , \; h=-1,0 \, ,\quad \ \bar\varepsilon :=
 \varepsilon^{(-1)}+ \varepsilon^{(0)} \, , \quad \
 \Theta := \g^{-1} \| R \|^T_{\vec p}
\ee and the corresponding quantities $ \e^{(h)}_+, \bar \e_+,
\Theta_+$ for $ R^+ $ with parameters $\vec p_+$ defined in
\eqref{magliana}.

\begin{proposition} {\bf (KAM step)} \label{KAMstep}
Assume that the parameters $ \vec p, \vec p_+$ (see \eqref{magliana})  satisfy \eqref{Giobbe+},
\eqref{lap1},  and that
\be\label{piccoloKAM} \d_+^{-1}
K^{2\t+1}\bar\e \quad is \ small \ enough  \, , \quad \Theta \leq 2^9 \, ,
\ee
where $ \d_+ $ is defined in \eqref{piccoF}.
Then, by \eqref{effe},   the solution  $ F \in {\cal R}_{rev}^{\leq 0} $ of the homological equation
\eqref{capocotta} satisfies  \eqref{piccoF}
 and the transformed vector field $  X^+$ in \eqref{ixspiu} is well
defined. The new normal form is \eqref{newNF}-\eqref{ennepiu} with frequencies satisfying
 Lemma \ref{frodo}.
The new perturbation $ R^+ \in {\cal R}_{rev} $ in
\eqref{newNP} satisfies $ R^+ = {\cal S} R^+ $ and (see \eqref{xhx})
\begin{eqnarray}
  \qquad \e^{(-1)}_+
  & \lessdot &
  \d_+^{- 2} K^{4\t+2} \bar\e^2 +   \e^{(-1)}  \,
  e^{- K \min \{s - s_+,\,
\aaa - \aaa_+\}}
 \nonumber \\
\e^{(0)}_+ & \lessdot & \d_+^{- 2 } K^{4\t+2} \big(  \e^{(-1)} +
\bar\e^2 \big)  +
   \e^{(0)} \, e^{- K \min \{s - s_+,\,
\aaa - \aaa_+\}} \label{kamst}
\\
 \label{ultimax+}
\Theta_+ \, &\leq&
  \Theta
( 1 + C \d_+^{- 2 } K^{4\t+2} \bar\e )  \, .
 \end{eqnarray}
\end{proposition}
\begin{pf}
We analyze each term of $R^+ $ in \eqref{newNP}.
We first claim that
\begin{equation}
\label{2pez} \big\|  {\rm ad}_F ( R^{\leq 0} ) \big\|_{\vec p_+}^T
+ \big\| {\mathop\sum}_{j \geq 2} \frac{1}{j!}{\rm ad}_F^j (X) \big\|_{\vec
p_+}^T \lessdot \d_+^{-2} \gamma  K^{2 (2\t+1)} \bar\varepsilon^2
\, .
\end{equation}
We have
\begin{eqnarray*}
{\mathop \sum}_{j \geq 2} \frac{1}{j!}{\rm ad}_F^j (X) \!\! \! & = & \!\! \! {\mathop \sum}_{j \geq 2}
\frac{1}{j!}{\rm ad}_F^j ({\cal N} + R) =
{\mathop \sum}_{j \geq 2} \frac{1}{j!}{\rm ad}_F^{j-1} ( {\rm ad}_F {\cal N}) +  {\mathop \sum}_{j \geq 2} \frac{1}{j!}{\rm ad}_F^j (R)  \\
& \stackrel{\eqref{capocotta} } = & {\mathop\sum}_{j \geq 2}
\frac{1}{j!}{\rm ad}_F^{j-1} ( [R^{\leq 0}_K] - R^{\leq 0}_K) +
{\mathop\sum}_{j \geq 2} \frac{1}{j!}{\rm ad}_F^j (R) \, .
\end{eqnarray*}
As we have already noticed, by \eqref{Giobbe+}, \eqref{lap1},
\eqref{piccoF} we can apply Proposition \ref{main2} (with $\vec p
\rightsquigarrow \vec p_\star,$
 $\vec p {\,}{'} \rightsquigarrow \vec p_+$,
 $\d \rightsquigarrow \d_+ $,  $h \rightsquigarrow 2$)
obtaining
\begin{equation}\label{uno}
\Big\| {\mathop\sum}_{j \geq 2} \frac{1}{j!}{\rm ad}_F^j (R) \Big\|_{\vec
p_+}^T
 \stackrel{\eqref{gPhif12bis}}\lessdot
 \Big( \d_+^{-1} \| F \|^T_{\vec p_\star} \Big)^2 \| R \|^T_{\vec p_\star}
 \stackrel{\eqref{effe}, \eqref{xhx} } \lessdot
  \d_+^{-2} K^{2 (2\t+1)} {\bar \e}^2 \g  \, \Theta\,.
\end{equation}
In the same way we get (with $h \rightsquigarrow 1$)
\begin{eqnarray}
&&\Big\| {\mathop\sum}_{j \geq 2} \frac{1}{j!}{\rm ad}_F^{j-1}
\big([R^{\leq 0}_K]-R^{\leq 0}_K\big) \Big\|_{\vec p_+}^T = \Big\|
{\mathop\sum}_{j \geq 1} \frac{1}{(j+1)!}{\rm ad}_F^{j} \big([R^{\leq
0}_K]-R^{\leq 0}_K\big) \Big\|_{\vec p_+}^T
\nonumber \\
&& \stackrel{\eqref{gPhif12bis}} \lessdot
   \d_+^{-1} \| F \|^T_{\vec p_\star}  \|
   [R^{\leq 0}_K]-R^{\leq 0}_K \|^T_{\vec p_\star}
   \leq
     \d_+^{-1} \| F \|^T_{\vec p_\star}  \| R^{\leq 0}_K \|^T_{\vec p_\star}
     \stackrel{\eqref{effe}, \eqref{xhx}}  \lessdot
   \d_+^{-1}  K^{2\t+1} \g \bar\varepsilon^2  \, .  \label{due}
\end{eqnarray}
Finally, by Proposition \ref{festa2}, applied with $\vec
p\rightsquigarrow\vec p_\star$, $ \vec p_1\rightsquigarrow \vec
p_+,$ $\d\rightsquigarrow  \d_+$ (note that conditions
\eqref{Leviatan0}-\eqref{Leviatan} follow by
\eqref{Giobbe+}-\eqref{lap1}), we get
\begin{equation}
\big\|  {\rm ad}_F ( R^{\leq 0} ) \big\|_{\vec p_+}^T
 \stackrel{\eqref{cippalippa}} \lessdot
\d_+^{-1} \| F \|^T_{\vec p_\star} \| R^{\leq 0} \|^T_{\vec
p_\star}
 \stackrel{\eqref{effe}, \eqref{xhx}}\lessdot
  \d_+^{-1} K^{2\t+1} \gamma \, \bar\varepsilon^2  \, . \label{tre}
\end{equation}
The bounds \eqref{uno}, \eqref{due},  \eqref{tre},  and $ \Theta
\leq 2^9 $ (see \eqref{piccoloKAM}), prove \eqref{2pez}.

\smallskip

We now prove \eqref{ultimax+}. Again by Proposition \ref{festa2} we
get
\begin{equation}
\big\|  {\rm ad}_F ( R^{\geq 1} ) \big\|_{\vec p_+}^T  \lessdot
 \d_+^{-1} \| F \|^T_{\vec p_\star}\| R^{\geq 1} \|^T_{\vec p}
 \stackrel{  \eqref{effe}, \eqref{xhx}}\lessdot
 \d_+^{-1} K^{2\t+1} \gamma \,  \bar\varepsilon \, \Theta  \label{tre0}
\end{equation}
and \eqref{ultimax+} follows by \eqref{newNP},
\eqref{carbonara4},  \eqref{fhT},  \eqref{xhx} \eqref{tre0},
\eqref{2pez} and $ \bar\varepsilon \leq  3 \Theta $ (which follows
by \eqref{xhx} and \eqref{fhT}).

\smallskip

We now consider $  R^{(h)}_+ $, $ h = 0,-1 $. Recalling the degree
decomposition  $ F = F^{(-1)} + F^{(0)} $,  formula
 \eqref{sumd}  implies that
the term $ {\rm ad}_F R^{\geq 1}  $ in \eqref{newNP} does not
contribute to  $ R^{(-1)}_+ $. On the other hand,  its
contribution to $ R^{(0)}_+ $ is  $ [R^{(1)}, F^{(-1)}] $.
Again by  \eqref{cippalippa}, \eqref{effe},  
\eqref{xhx} and \eqref{fhT}, we get
\begin{equation}
\label{contrib12} \|  [   R^{(1)}, F^{(-1)}]\|_{\vec p_+}^T
\lessdot \d_+^{-1} \g K^{2\t+1} \varepsilon^{(-1)} \Theta\,.
\end{equation}
The contribution of $ R - R^{\leq 0}_K  $ in  \eqref{newNP} to  $
R^{(h)}_+ $, $ h = 0,-1 $, is
$ \Pi_{|k|<K}\Pi_{|\pi|\geq K} R^{(h)}+ \Pi_{|k|\geq K} R^{(h)} $.
By \eqref{smoothT} (recall $s s_+^{-1}<2$),
\eqref{fhT},  and \eqref{xhx}, we get
\begin{equation}\label{lem:ultra}
 \big\| \Pi_{|k|<K}\Pi_{|\pi|\geq K} R^{(h)}+ \Pi_{|k|\geq K} R^{(h)} \big\|^T_{\vec p_+}
 \leq
3 e^{- K \min \{s - s_+,\, \aaa - \aaa_+\}} \g \varepsilon^{(h)}
\, .
 \end{equation}
In conclusion,  \eqref{kamst} follows by \eqref{newNP},
\eqref{carbonara4},  \eqref{2pez}, \eqref{contrib12},
\eqref{lem:ultra} and $ \Theta \leq 2^9 $.
\end{pf}

\noindent
{\bf KAM iteration.} Once the KAM step has been proved, the proof of Theorem \ref{thm:IBKAM}
is concluded by an usual KAM iteration.
The scheme is very similar to that in
 \cite{BBP1} (and \cite{BB10}) and we skip it. We only focus on the main difference,
 which is the symmetrization procedure.

For every $i\in\mathbb N$
we construct a
close-to-the-identity, analytic, even (Definition
\ref{rean}) change of variables $ \Phi^i $
(obtained as the time-1 flow of the
the solution $F_i$ of the homological equation
\eqref{capocotta} at the $i^{\rm th}$ step)
such that  (recall \eqref{cecio} and \eqref{ixspiu})
\begin{equation}\label{garbatella}
X_i:=\mathcal S \Phi_\star^i X_{i-1}=:\mathcal N_i+R_i
\,, \quad R_i\in\mathcal R_{rev}\,, \
R_i=\mathcal S R_i
\end{equation}
($\Phi^i_\star $ is the lift to the tangent
space (recall \eqref{Hnew}).
Since the algorithm is ``quadratic'' (recall \eqref{kamst}),
 the quasi-T\"{o}plitz (with suitable $i$-dependent parameters)
norm of the $-1$ and $0$ degree
terms of $R_i$ converges  \textsl{super-exponentially} to zero.
 Let
\begin{equation}\label{tormarancia}
X_\infty:= \lim_{i\to\infty} X_i=
\lim_{i\to\infty} \mathcal S \Phi_\star^i X_{i-1}=
\mathcal N_\infty +R_\infty \quad
\mbox{where}\ \ \ \mathcal N_\infty:= \lim_{i\to\infty} {\mathcal
N}_i \,, \quad R_\infty:= \lim_{i\to\infty} R_i  \,  .
\end{equation}
By \eqref{garbatella} and the convergence of
the $-1$ and $0$ degree
terms of $R_i$ we get
\begin{equation}\label{aventino}
R_\infty = {\cal S} R_\infty \, , \qquad  R_\infty^{\leq 0} = 0 \,
.
\end{equation}
The  transformation $\Phi$ in \eqref{eq:Psi} is defined by
$ \Phi:=\lim_{\nu\to\infty} {\bar \Phi} \circ
\Phi^0\circ \Phi^1 \circ \cdots \circ \Phi^\nu $
where $ \bar \Phi $ is defined in section \ref{volkswagen} as the time $
1$-flow of  $ F $ defined in Lemma \ref{FST}. The map $ \Phi $ is even
because $ \Phi^i $, $ i \geq  0 $, and $ \bar \Phi $ are even.
Let us show  the proof of  \eqref{Hnew}.
We have  that
\begin{equation}\label{labaro}
\XX_\infty=\Phi_\star \XX= \lim_{i\to\infty} \XX_i \quad \mbox{where}\quad
\XX_{i} := \Phi^i_\star \XX_{i-1}\,, \  i\geq 1\,, \ \XX_0 \
\mbox{defined in }\eqref{salaria}\,.
\end{equation}
The vector field $ \XX_\infty  \in {\cal R}_{rev} $ 
because $ \XX_0   \in {\cal R}_{rev}  $ (see \eqref{salaria}) and
each $ \XX_i \in {\cal R}_{rev} $ because $ \Phi^i_\star = e^{{\rm
ad}_{F_i}} $ with $F_i \in\mathcal R_{a-rev} $ (then use  Lemma
\ref{embeh2}).
The relation between the
``auxiliary'' vector field $X_\infty$ and the ``true'' vector field
$\XX_\infty$ is given by the following
\begin{lemma}\label{testaccio}
$ (\XX_\infty)_{|E}=(X_\infty)_{|E} $.
\end{lemma}
\begin{pf}
The lemma follows by proving
$(\XX_i)_{|E}=(X_i)_{|E},$ $\forall\, i\geq 0.$ The inductive
basis for $i=0$ is \eqref{tiburtina}. Let us assume that
$(\XX_{i-1})_{|E}=(X_{i-1})_{|E}.$ Then
$$
(\XX_i)_{|E}-(X_i)_{|E}
\stackrel{\eqref{garbatella},\eqref{labaro}}{=}
(\Phi^i_\star\XX_{i-1})_{|E}-(\mathcal S\Phi^i_\star X_{i-1})_{|E}
= \big(\Phi^i_\star(\XX_{i-1}-X_{i-1}) \big)_{|E} \equiv 0
$$
 by Proposition \ref{carbonara1} and   Lemma \ref{montalto}
(used with $X\rightsquigarrow \XX_{i-1}-X_{i-1}$, $Y
\rightsquigarrow F_i$  with
$e^{{\rm ad}_{F_i}}=\Phi^i_\star$).
\end{pf}

We have already chosen $\mathcal N_\infty$ in \eqref{tormarancia},
then $\mathcal P_\infty$ in \eqref{Hnew} is $\mathcal
P_\infty=\XX_\infty - \mathcal N_\infty.$ It is now simple to show
that $ \big(\PP_\infty^{\leq 0}\big)_{|E} = 0.$ Indeed
$$
\big(\PP_\infty^{\leq 0}\big)_{|E} = \big((\XX_\infty-\mathcal
N_\infty)^{\leq 0}\big)_{|E} \stackrel{\eqref{tormarancia}}{=}
\big((\XX_\infty-X_\infty + R_\infty)^{\leq 0}\big)_{|E}
\stackrel{\eqref{aventino}}{=} \big((\XX_\infty-X_\infty)^{\leq
0}\big)_{|E} \stackrel{\eqref{dicastro}}{\equiv} 0\,.
$$
by Lemma \ref{testaccio}. Finally $\PP_\infty\in\mathcal R_{rev}$
because $\mathcal N_\infty $ and
$\XX_\infty \in\mathcal R_{rev}$. This concludes the proof of
\eqref{Hnew}.
\\[1mm]
{\bf Proof of  \eqref{sublevels}}. 
By $ {\det} A \neq 0 $ and \eqref{muflone}  the action-to-frequency map $ \ot^\infty $ 
 is invertible. Introducing $ \zeta =  \ot^\infty (\xi) $ as parameters, we obtain
 $\xi = (\ot^\infty)^{-1} (\zeta) = A^{-1} ( \zeta - {\bar \om} ) + O(\e\gamma) $ and,
 using also  \eqref{omegaxi}, 
$$
\ot^\infty (\xi) \cdot h + \Om^\infty_i (\xi) - \Om^\infty_j (\xi)
= f_{h,i,j}(\zeta)+ r_{i,j}(\zeta)\,, \quad  | r_{i,j}|=O(\gamma \e), \  | r_{i,j}|^{\rm lip}=O(\e) \, , 
$$
 $$
f_{h,i,j}(\zeta )
:=c_{h,i,j}\cdot \zeta + d_{i,j}\,, \ 
c_{h,i,j}:=  h  +  (\l_i^{-1} - \l_j^{-1}) A^{-T} \vec a  \,,
 \ 
d_{i,j}:=
\l_i - \l_j  
-  (\l_i^{-1} - \l_j^{-1}) {\bar \om} \cdot  A^{-T} \vec a \, .   
$$
Then  \eqref{sublevels}
follows immediately 
if  $ | c_{h,i,j} | >  \bar c >  0 $ because $ | c_{h,i,j} \cdot \partial_{\zeta} f_{h,i,j} | \geq \bar c^2 >  0 $
and  $ |  r_{i,j}|^{\rm lip}=O(\e) $. 
Now, since $ h \in \Z^{n/2} \setminus \{0 \} $ and  $ |  (\l_i^{-1} - \l_j^{-1}) A^{-T}  \vec  a| = O ( \l_i^{-1} + \l_j^{-1}) $,
the coefficient  $ | c_{h,i,j} | > 1 / 2 $,  for  $\min\{|i|,$ $ |j|\} \geq C $ large. 
On the other hand, if $ | i | \leq C $ and  $ |j| \geq C_0 $ with $C_0$ large  enough (or permuting the role of $i$ and $j$) 
the coefficient  $ | d_{i,j}|\geq 1 $. In this case
 $  | f_{h,i,j}   + r_{i,j}  | > 1 /8 $ for all $ \zeta \in 
\ot^\infty ( {\cal O} ) $ unless  $ |c_{h,i,j} \cdot \bar \om | \geq 1 /4 $  (for $\e,\rho$ small). Hence
$ | \bar \om  \cdot \partial_{\zeta} f_{h,i,j} | =  |  \bar \om \cdot c_{h,i,j} | > 1 / 4 $ 
and, again, \eqref{sublevels} follows. 
Finally, 
the first condition in \eqref{pota} and $ h \neq 0 $ imply 
$ \min \{ | c_{h,i,j} |\ {\rm for}\  |i|, |j| \leq C_0 \}  > 0 $ and so \eqref{sublevels}.

\section{Proof of Theorem \ref{thm:DNLW}}\label{sec:BNF} \setcounter{equation}{0}

By hypothesis, the analytic nonlinearity  $ g $ has the convergent Taylor expansion
$$
g(\mathtt x,\mathtt y,\mathtt y_{\mathtt x},\mathtt v) = \kappa_1
{\mathtt y}^3 +\kappa_2 {\mathtt y} {\mathtt y}_{\mathtt x}^2
+\kappa_3 {\mathtt y} {\mathtt y}_t^2 + {\mathop \sum}_{k+h+l \geq
5} \, g^{(k,h,l)}( \mathtt x ) \mathtt y^k \mathtt y^h_{\mathtt x}
\mathtt v^l
 $$
where  $ k, h, l \in \N $ and
\begin{equation}\label{boudicca}
\| g^{(k,h,l)}\|_{a_0,p} < C^{k+h+l} \quad {\rm for\  some}\quad
a_0>0 \, , \ p>  1 / 2 \, , \ C > 0  \,,
\end{equation}
having identified each function $g^{(k,h,l)}( \mathtt x )$ with
the Fourier series $\{g^{(k,h,l)}_{j_0}\}_{j_0 \in \Z}\in
\ell^{a_0,p}$, recall \eqref{ellcal}. As phase space we consider $ u , \bar u  \in\ell^{a,p} $
with  $ a := a_0 / 2  $.  The coefficients $\mathtt g_j^+ $ in \eqref{piggione} are
\begin{eqnarray}\label{patatinefritte}
 \mathtt g_j^+ :=  \mathtt g_j
 &=&
-\sum_{d=3, d\geq  5}\sum_{\atop{   j_0+ \sum_{i=1}^d \s_i j_i=j}
}(\sqrt{2})^{-d-1} \mathtt g_{\vec \s,\vec \jmath,j_0} u^{\vec
\s}_{\vec \jmath} =: \mathtt g_j^{(=3)}+ \mathtt g_j^{(\geq 5)}
\end{eqnarray}
where $\vec \jmath=(j_1,\ldots,j_d) \in \Z^d $, $\vec
\s=(\s_1,\ldots,\s_d) \in \{\pm \}^d $ and
 $u^{\vec \s}_{\vec \jmath}= {\mathop \prod}_{i=1}^d u_{j_i}^{\s_i}$.
The coefficients $g_{\vec \s,\vec \jmath,j_0}$ are 
$$
g_{\vec \s,\vec \jmath,j_0} = \sum_{ h+k+l= d} (-1)^l \, \ii^{h+l}
\s_{k+1}\cdots \s_{k+h+l} \frac{ j_{k+1}\cdots  j_{k+h}
}{\l_{j_1}\cdots \l_{j_{k+h}}} g^{(k,h,l)}_{j_0}\,.
$$
We consider \eqref{tuc} as the equations of motion of the vector
field $ \mathcal N_0 + G $ where (recall \eqref{piggione})
\begin{equation}\label{lupin}
\mathcal N_0  := \sum_{\s=\pm,\, j\in \Z}  \s \ii\l_j u_j^\s
\partial_{u_j^\s} \, , \quad
G = \sum_{\s=\pm, j\in\Z}G^{(u_j^\s)}\partial_{u_j^\s} \, ,
\quad 
    G^{(u_j^\s)} := \ii \s \mathtt g_{\s j}\, ,
    \quad     G=G^{(=3)}+G^{(\geq 5)} \, .
\ee
Note that
$ G^{(u^+_{-j})}=- G^{(u^-_j)} $
and that $G^{(=3)}$ has zero momentum by \eqref{patatinefritte},
\eqref{lupin}.
Moreover 
$G$ is  reversible (w.r.t. the involution $S$ in \eqref{Seven}),
real-coefficients, real-on-real, even, namely $G\in \mathcal
R_{rev}$ (Definition \ref{pajata} in absence of $x,y$-variables).

\begin{lemma}\label{aragorn}
Set $ \aaa:= a_0 \slash 2  $ (where $ a_0 $ is defined in \eqref{boudicca}). 
Then,  for  $ R := R ( C ) > 0 $ small enough (where $ C $ is defined in \eqref{boudicca}),  it results
\begin{equation}\label{boundG}
\|G\|_{R,\mathtt a}\,, \, \|G^{(=3)}\|_{R,\mathtt a}
  \lessdot R^2\,,   \quad
  \|G^{(\geq 5)}\|_{R,\mathtt a}
  \lessdot R^4\,.
\end{equation}
Moreover  $G, G^{(=3)}, G^{(\geq 5)} \in \mathcal Q^T_{R,\mathtt
a}(N_0,\CCquarti,\quattrocc)$, for $N_0$ satisfying
\eqref{caracalla}, and
\begin{equation}\label{radiatore}
  \|G\|_{R,\mathtt a,N_0,\CCquarti,\quattrocc}^T\,,\,
    \|G^{(=3)}\|_{R,\mathtt a,N_0,\CCquarti,\quattrocc}^T
 \lessdot R^2 \, ,
 \quad
  \|G^{(\geq 5)}\|_{R,\mathtt a,N_0,\CCquarti,\quattrocc}^T
 \lessdot R^4
 \,.
\end{equation}
\end{lemma}

\begin{pf}
We first note that (recall also that  $ a :=  a_0 \slash 2  $)
\be\label{chiave} \Big\| \Big( \sum_{  j_0+\sum_{i=1}^d \s_i j_i=j
} e^{\mathtt  a |j_0|}
 |g^{(k,h,l)}_{j_0}| |u^{\vec \s}_{\vec \jmath}|
 \Big)_{j\in\Z}
 \Big\|_{a,p}
   \lessdot
    \|g^{(k,h,l)}\|_{a_0,p}(\|u\|_{a,p}+\| \bar u\|_{a,p})^d \, .
\ee Indeed
$$
\sum_{  j_0+\sum_{i=1}^d \s_i j_i=j  } e^{\mathtt  a |j_0|}
|g^{(k,h,l)}_{j_0}| |u^{\vec \s}_{\vec \jmath}|
 \leq
\big(f^{(k,h,l)}* \tilde u* \tilde u*  \dots *\tilde u\big)_{j}\,,
\qquad \forall\, j\in \Z\,,
$$
where   $f^{(k,h,l)} :=(e^{\mathtt  a |j_0|}
 g^{(k,h,l)}_{j_0})_{j_0\in \Z}$,  $ \tilde u
 :=
  (\tilde u_n)_{n\in\Z} $,
  $ \tilde
u_n:=|u_n|+|\bar u_n|$,   $*$ denotes the convolution of sequences
and
$$
 \|f^{(k,h,l)} *\tilde u* \tilde u*  \dots *\tilde u \|_{a,p}\lessdot
 \|f^{(k,h,l)} \|_{a,p} \|\tilde u\|_{a,p}^d\lessdot
\| g^{(k,h,l)} \|_{a_0,p}  (\|u\|_{a,p}+\|\bar u\|_{a,p})^d
$$
by the Hilbert algebra property of $\ell^{a,p}$ and since $\mathtt a = a_0-a $.

Now we  rewrite the sum in \eqref{patatinefritte} as $ \mathtt
g_j= {\mathop\sum}_{|\a|+|\b|\geq 3} (\mathtt g_j)_{\a,\b} u^\a \bar u^\b $
where $(\mathtt g_j)_{\a,\b}$ can be explicitly computed from
\eqref{patatinefritte} but has a complicated combinatorics. In
order to compute the norm $ \|  G \|_{R, \mathtt a} $ we note that
${1}/{|\l_l|}\lessdot 1$, ${|l|}/{|\l_l|} \leq 1 $ and
\be\label{patatinefritte4} u^{\vec \s}_{\vec \jmath}= u^\a \bar
u^\b\, \quad\Longrightarrow\quad
 \pi(\a,\b; u^\s_j)= {\mathop\sum}_{1\leq i\leq d} \s_i j_i -\s j\,.
\end{equation}
We have (recall \eqref{lupin})
\begin{eqnarray*}
  \|G\|_{R,\mathtt a}
  &\stackrel{\eqref{normaEsr}}=&
\sup_{\|u\|_{a,p},\|\bar u \|_{a,p} < R } R^{-1} \sum_{\s=\pm}
\Big\|\Big(\sum_{|\a|+|\b|\geq 3} e^{\mathtt  a |\pi(\a,\b;
u^\s_j)|}|(g_{\s\! j})_{\a,\b}| |u^\a| |\bar u^\b| \Big)_{j\in\Z}
\Big\|_{a,p}
  \\
  &\stackrel{\eqref{patatinefritte4}}\lessdot&
  R^{-1} \!\!\! \sup_{\|u\|_{a,p},\|\bar u \|_{a,p} < R }
  \Big\| \Big( \sum_{d\geq  3}(\sqrt{2})^{-d-1}\sum_{ j_0+\sum_{i=1}^d \s_i j_i=j } 
  \sum_{ h+k+l= d}e^{\mathtt  a |j_0|} |g^{(k,h,l)}_{j_0}| |u^{\vec \s}_{\vec \jmath}|\Big)_{j\in\Z}
\Big\|_{a,p}
\\
& \stackrel{\eqref{chiave}} \lessdot& R^{-1}  \!\!\!
\sup_{\|u\|_{a,p},\|\bar u \|_{a,p} < R } \sum_{d\geq
3}\sum_{h+k+l=d}(\sqrt{2})^{-d-1}
 \|g^{(k,h,l)}\|_{a_0,p}(\|u\|_{a,p}+\| \bar u\|_{a,p})^d
 \stackrel{\eqref{boudicca}}\lessdot R^2
\end{eqnarray*}
proving \eqref{boundG} for $ R $ small enough with respect to the constant $ C $ in \eqref{boudicca}.

Let us now prove the estimate \eqref{radiatore} for  the
quasi-T\"oplitz norm of $ G$ (the estimates for $G^{(=3)} $ and $
G^{(\geq 5)} $ are analogous). For $N\geq N_0$, by
\eqref{patatinefritte} and \eqref{lupin}
 we deduce that
 $$
\Pi_{N,\CCquarti,\quattrocc}  G = \sum_{ |m|,|n|>(\CCquarti) N}
G^{\s,m}_{\s', n }\  u_n^{\s'}
\partial_{u_m^\s } = \tilde G + N^{-1} \hat G
$$
where (recall \eqref{perla}, \eqref{bibi}, \eqref{giglio})
\begin{eqnarray*}
G^{\s,m}_{\s', n }  & := & -\ii \s \sum_{d\geq  2}\sum_{
\sum_{i=1}^d |j_i|< 4 N^L\,, |j_0|< N^b\atop
 j_0+\sum_{i=1}^d \s_i j_i= \s m -\s' n }
 (\sqrt{2})^{-d-2}
\!\!\!\!\!\!\!\sum_{ h+k+l= d+1} (\ii)^{h+l}(-1)^l
g^{(k,h,l)}_{j_0} c_{\vec \s,\vec \jmath,\s ',n}^{(k,h,l)} u^{\vec
\s}_{\vec \jmath}  \label{riga1}
\\c_{\vec \s,\vec \jmath,\s ',n}^{(k,h,l)}&:=&\frac{  j_{k+1}\dots j_{k+h-1}\s_{k+1}\dots \s_{k+h+l-1}}
 {\l_{j_1}\dots \l_{j_{k+h-1}} } \, \bigg(k\frac{ \s_k j_k}{\l_n}
 +
h\frac{ \s ' n}{\l_n}+  l \frac{\s '}{\l_{j_{k+h}}}\bigg)
 \, .  \label{riga3}
\end{eqnarray*}
The  T\"oplitz approximation $\tilde G $
is obtained 
substituting the coefficients $c_{\vec \s,\vec \jmath,\s
',n}^{(k,h,l)} $
 with their T\"oplitz approximation $\tilde c_{\vec \s,\vec \jmath,\s ',n}^{(k,h,l)}$ defined by
replacing  $1/ \l_n $ by  $ 0 $ and  $ n / \l_n $  by the sign $
{\mathtt s}(n) $.

Since $ 0 \leq \l_n-|n|\leq \sqrt \mathtt m$, $ \forall n\in\Z$,
and $\l_n\geq |n|> (\CCquarti) N $, the Taylor coefficients of
$\tilde G$ and of the corresponding defect  $\hat G$ are uniformly
bounded. Then, arguing as in the proof of \eqref{boundG},
we deduce that $ \|\tilde G\|_{R,\mathtt a} $, $ \|\hat
G\|_{R,\mathtt a} \lessdot
 R^2 $.
  Note that
$\tilde{c}_{\vec \s,\vec \jmath,\s ',n}^{(k,h,l)}$ depends on $n$
only through $\mathtt s(n).$ Since by \eqref{giglio} $\mathtt
s(n)=\s \s' \mathtt s(m)$ we have that
  $\tilde G\in \mathcal T_{R,\mathtt a}(N, \CCquarti,\quattrocc)$
  (Definition \ref{matteo_aa}).
  By Definition \ref{topbis_aa}
 we get \eqref{radiatore}.
\end{pf}

For the Birkhoff normal form step we need
the following lemma 
proved in \cite{BBP1} (Lemma 7.2 and formula (7.21), see also
\cite{Po3}).

\begin{lemma}\label{vaccinara}
There exists an absolute constant $c_* >  0 $, such that,
for every $\mm\in(0,\infty) $ and  $j_i\in\Z$, $\s_i=\pm$,
$i=1,2,3,4$ satisfying $ \s_1 j_1+\s_2 j_2 + \s_3 j_3 +\s_4 j_4 =
0 $ but not satisfying \be\label{zerisecchi} j_1=j_2 \, ,  \
j_3=j_4 \, , \  \s_1 = -\s_2 \, , \ \s_3=-\s_4 \ {\rm (or  \
permutations \ of\ the\ indexes)} \, , \ee we have
\begin{equation}\label{pecorino}
|\s_1 \l_{j_1}+\s_2 \l_{j_2}+\s_3 \l_{j_3}+\s_4 \l_{j_4}| \geq
\frac{c_* \mm}{(n_0^2+\mm)^{3/2}} > 0 \quad \mbox{where} \quad
n_0:= \min \{ \langle j_1 \rangle,  \langle j_2  \rangle,
\langle j_3 \rangle,   \langle j_4  \rangle \} \, .
\end{equation}
\end{lemma}

Then we define the projections $G_1$ and $G_2$ of $G^{(=3)}$ as follows:
the vector field $-\ii G_1^{(u_j^+)}$ is the projection of
${\mathtt g}_j^{(=3)}$ (recall \eqref{patatinefritte}) onto the
indexes $(\s_1,\s_2,\s_3,+)$, $(j_1,j_2,j_3,j)$ which satisfy
\eqref{zerisecchi} with $j_1\in\mathcal I.$ Let $-\ii
G_2^{(u_j^+)}$ be the projection of ${\mathtt g}_j^{(=3)}$
 onto the indexes
 $j_1,j_2,j_3\not\in\mathcal I$
if $j\not\in\mathcal I$ and zero otherwise. We have that
 \begin{equation}\label{heidi}
    \|G_1\|_{R,\aaa} =\|G_1\|_{R,0}\,,
    \ \|G_2\|_{R,\aaa}=\|G_2\|_{R,0} \lessdot R^2 \, ,
\end{equation}
and, for $N_0'$ large enough,
\begin{equation}\label{radiatoreter}
    \|G_1\|_{R,\aaa,N_0,\CCquarti,\quattrocc}^T
 \lessdot R^2\,,
 \qquad
\|G_2\|_{R,\aaa,N_0,\CCquarti,\quattrocc}^T
 \lessdot R^2
 \,.
\end{equation}
The estimates \eqref{heidi} and \eqref{radiatoreter} follows by
\eqref{lessore} and the analogous estimates \eqref{boundG} and
\eqref{radiatore} for $G$, since $G_1,G_2$ are  projections
(recall \eqref{PROJE}) of $G$, satisfying \eqref{tototo}.

\begin{proposition}\label{BNF}
{\bf (Birkhoff normal form)} For any $ {\cal I} $ as in
\eqref{Isym},  and $ \mm > 0 $, there exists $ R_0 >  0 $ and
a real analytic  change of variables
$ \Gamma 
: B_{R/2}\times B_{R/2} \subset \ell^{a,p}\times \ell^{a,p} \to $
$ B_{R}\times B_{R} \subset \ell^{a,p}\times \ell^{a,p} $, $ 0 < R < R_0  $, that takes the vector field $ \mathcal N_0 + G $  into
\be\label{Birkhoff} \big( D \Gamma^{-1}[\mathcal N_0+G]\big) \circ
\Gamma =
 \mathcal N_0+G_1+G_2 + G_3 \ee
where $ G_1 $, $ G_2$ satisfy \eqref{radiatoreter},
$G_3$ satisfy $ G_3^{(u^+_{-j})}=- G_3^{(u^-_j)} $ 
 and for $N_0'$ large enough
\begin{equation}\label{radiatorebis}
\|G_3\|_{R/2,\aaa/2,N_0',7/4,3}^T
 \lessdot R^4
 \,.
\end{equation}
Finally $\mathcal N_0 +G_1+G_2+G_3\in \mathcal R_{rev}$ (recall
Definition \ref{pajata} in absence of $x,y$-variables).
\end{proposition}

\begin{pf}
Let us define the generating function $ F := {\mathop \sum}_{j\in \Z,\s=\pm}
F^{(u_j^\s)}  \partial_{u_j^\s} $ with \be\label{laF} F^{(u_j^\s)}
:= \sum_{ {\s_1 \l_{j_1} + \s_2 \l_{j_2} + \s_3 \l_{j_3} - \s
\l_j\neq 0 \atop (j_1,j_2,j_3,j)\notin (\mathcal I^c)^4, \vec
\s\cdot \vec \jmath = \s j },   } -\frac14 \frac{\s   }{  \s_1
\l_{j_1} + \s_2 \l_{j_2} + \s_3 \l_{j_3} - \s \l_j } \
 {\mathtt g}_{\vec \s, \vec \jmath,0}\
u_{\vec \jmath}^{\vec \s}\,. \ee 
By Lemma \ref{vaccinara} and
arguing as in Lemma \ref{aragorn} we get  $ \|F\|_{R,\mathtt a}   =
\|F\|_{R,0}\lessdot R^2 $. Moreover we claim that
\begin{equation}\label{radiatoreF}
    \|F\|_{R,\mathtt a,N_0,\CCquarti,\quattrocc}^T
 \lessdot R^2\,.
\end{equation}
For $ N \geq N_0$, by \eqref{laF}
 we wish to write
$\Pi_{N,\CCquarti,\quattrocc} F = \tilde F + N^{-1} \hat F $,
where (recall \eqref{perla}, \eqref{bibi} and \eqref{giglio})
$$
\tilde F =
{\mathop \sum}_{ |m|,|n|>(\CCquarti) N} \tilde
F^{\s,m}_{\s', n }\  z_n^{\s'} \partial_{z_m^\s }
$$
is T\"{o}plitz. We define $\tilde F $ by using \eqref{laF} with
$j,j_3,\s_3\rightsquigarrow m,n,\s'$  and substituting as follows:
${\mathtt g}_{\vec \s, \vec \jmath,0}$ by its
 T\"{o}plitz approximation (given in Lemma \ref{aragorn})
 and
$d:=\s_1\l_{j_1}+\s_2\l_{j_2}+\s'\l_n-\s\l_m$ by $\tilde
d:=\s_1\l_{j_1}+\s_2\l_{j_2}+\s' |n|-\s |m|$. To estimate the
T\"{o}plitz defect $\hat F$ we consider first the case
$\s=\s'$. We have
$$
|d-\tilde d|= \big|\l_n-\l_m -|n|+|m|\big| \lessdot ( |n|^{-1} +
|m|^{-1} ) \lessdot |n|^{-1}\,,
$$
noting that  $1/2 \leq |n|/|m|\leq 2$ by $\s_1 j_1 + \s_2 j_2 = \s
m -\s' n$ and $|j_1|+|j_2|<4 N^L$, for $N\geq N_0$ large enough.
Then, since by \eqref{pecorino}, $1\lessdot |d|,$ for
$|n|\geq(\CCquarti) N$ and $N_0$ large enough, $ 1\lessdot
|d|-|\tilde d-d|\leq |\tilde d| .$ In particular $|\tilde d|\geq
const.>0 > 0 $ and $\tilde F, \hat F$ are well defined. Moreover
$$
\left|\frac{d}{\tilde d}-1\right|=\left|\frac{1}{\tilde
d}(d-\tilde  d)\right| \lessdot \frac{1}{|n|} \qquad {\rm and}
\qquad \big| \l_n-|n|\big|\lessdot \frac{1}{|n|}
$$
and, therefore, $ | |n|- d \tilde  d^{-1}\l_n \big|\lessdot 1 $.
In the case
 $\s=-\s'$, since $|j_1|+|j_2|<4 N^L$ and
$\l_m\geq |m|\geq N$,
  we get $|d|\geq |n|$.
Recollecting we have that, both in the case $\s=\s'$ and
$\s=-\s'$, the Taylor coefficients of $\tilde F, \hat F$ are
uniformly bounded and, arguing as in the proof of Lemma
\ref{aragorn}, we get $ \|\tilde F\|_{R,\mathtt a}\,, \ \|\hat
F\|_{R,\mathtt a} \lessdot
 R^2 $.
 We note that $\tilde F\in \mathcal T_{R,\mathtt a}(N, \CCquarti,\quattrocc)$;
 indeed $\s=\s'$ and by \eqref{giglio}
 $\mathtt s(m)=\mathtt s(n)$, so that
 $\tilde  d:=\s_1\l_{j_1}+\s_2\l_{j_2}+
 \mathtt s(m)(\s' n-\s m)$.
Then by Definition \ref{topbis_aa}
 we deduce \eqref{radiatoreF}.

With $\mathcal N_0$ defined in \eqref{lupin} we have
$$
\big[ \mathcal N_0 , u_{j_1}^{ \s_1} u_{j_2}^{ \s_2} u_{j_3}^{
\s_3} \,
\partial_{u_j^\s} \big] = \ii ( \s_1 \l_{j_1} + \s_2 \l_{j_2} +
\s_3 \l_{j_3} - \s \l_j ) u_{j_1}^{ \s_1} u_{j_2}^{ \s_2}
u_{j_3}^{ \s_3} \, \partial_{u_j^\s} \, .
$$
Then $ F $ in \eqref{laF} solves the homological equation $
[\mathcal N_0,F] + G^{(=3)}  + G^{(=3)} = G_1+ G_2 $. Then we
define $\Gamma$ as the time-$1$ flow generated by the vector field
$F$ and  \eqref{radiatorebis} follows by Proposition \ref{main2}
taking $R<R_0$ small enough and $N_0'$ large enough.

We claim that $F\in \mathcal R_{a-rev}$. Indeed $F$ is
real-on-real (recall Definition \ref{rean})  by \eqref{laF}. $F$
is anti-real-coefficients since the Taylor coefficients in
\eqref{laF} are real. $ F $ is anti-reversible (recall Definition
\ref{reversVF}) with respect to
 the involution $S$ in \eqref{Seven} since by \eqref{laF} we have
$F^{(u_j^\s)}\circ S=F^{(u_{-j}^{-\s})}$. Finally $ F $ is even
(recall Definition \ref{rean}) since, again by \eqref{laF}
$F^{(u_j^\s)}_{|E}=F^{(u_{-j}^\s)}_{|E}$ (with $E$ defined in
\eqref{parit}).

Then $\mathcal N_0 +G_1+G_2+G_3=e^{{\rm ad}F} (\mathcal N_0+G)\in
\mathcal R_{rev}$ by Lemma \ref{embeh2}.
\end{pf}

\subsection{Action-angle variables and conclusion of Proof of Theorem \ref{thm:DNLW}}\label{nutella}

Let us denote by  $  (u^+,u^-) =\Phi(x,y,z^+, z^-;\xi)  $ the
change of variable introduced in \eqref{AAV}. For $\rro>0$, let
(recall \eqref{Isym}) \be\label{paraO} \mathcal O_\rro := \big\{
\xi \in \R^{n/2}
 \, : \,  \rro / 2\leq
\xi_j \leq \rro \, , \ j\in\mathcal I^+ \big\} \, . \ee A vector
field $X=(X^{(u^+)},X^{(u^-)})$ is transformed by the change of
variable $\Phi$ in
\begin{eqnarray}
&&Y:=\Phi_\star X=\big( D\Phi^{-1}[X] \big)\circ \Phi\,,
\ \ {\rm with} \ \  Y^{(z^\s_j)}= X^{(u^\s_j)}\circ
\Phi\,,\quad \s=\pm\,,\quad j\in\Z\setminus \mathcal I\, ,
\label{ACAN}
\nonumber\\
&& Y^{(x_j)}= -\frac{\ii}{2} \Big( \frac{1}{u^+_j} X^{(u^+_j)} -
\frac{1}{u^-_j} X^{(u^-_j)}\Big)\circ \Phi\,, \quad Y^{(y_j)}=
\Big( u^-_j X^{(u^+_j)} + u^+_j X^{(u^-_j)}\Big)\circ \Phi\,,\quad
j\in\mathcal I\,, \nonumber
\end{eqnarray}

\begin{lemma} {\bf (Lemma 7.6 of \cite{BBP1})}
Let us take
\begin{equation}\label{condro}
0<16 r^2 < \rro\,,\quad  \
 \rro  = C_* R^2  \quad {\rm with} \quad C_*^{-1} := 48 n \kappa^{2p} e^{2 ( s+ a\kappa)} \, .
\end{equation}
where $ a=a_0/2 $, $ p>1/2 $ and 
$ \kappa $ is defined  in \eqref{caracalla}.
Then, for all $\xi \in  \mathcal O_\rro \cup \mathcal  O_{2 \rro}
$, the map $ \Phi(\,\cdot\,;\xi) : D(s,2r) \to B_{R/2} \times
B_{R/2}  \subset \ell^{a,p} \times \ell^{a,p}  $ is well defined
and analytic  $($$D(s,2r)$ is defined in \eqref{Dsr}$)$.
\end{lemma}
Given a vector field $ X : \, B_{R/2} \times B_{R/2} \to
\ell^{a,p} \times \ell^{a,p} $, the previous Lemma and
\eqref{ACAN} show that the transformed vector field $ Y:=\Phi_\star
X : D(s, 2r) \to \ell^{a,p} \times \ell^{a,p} $.
It results  that,  
if $ X $ is quasi-T\"oplitz in the variables $ (u, \bar u) $ then
$ Y $ is quasi-T\"oplitz in the variables $ (x,y,z, \bar z ) $
(see Definition \ref{topbis_aa}). We define the space of vector fields 
$$
{\cal V}^d_{R,\aaa} := \Big\{  X := X(u, \bar u ) \ : \ \|X\|_{R,\aaa}<\infty\ \ {\rm and}\  \
X^{(u_j^\s)} = \sum_{|\a^{(2)}+\b^{(2)}| \geq d}
X^{(u_j^\s)}_{\a,\b}  u^\a\bar u^\b \Big\}\, .
$$
\begin{proposition}{\bf (Quasi--T\"oplitz)} \label{qtop0}
Let $ N_0,\teta,\mu,\mu' $ satisfy \eqref{caracalla} and
\begin{equation}\label{fuso}
(\mu'-\mu)N_0^L> N_0^b\,,\qquad  N_0
2^{-\frac{N_0^b}{2\kappa}+1}<1 \, .
\end{equation}
If $ X \in \mathcal Q_{R/2,\aaa}^T (N_0,\theta,\mu') \cap {\cal
V}_{R/2,\aaa}^{d} $ with $ d = 0, 1 $, then $ Y := \Phi_\star X\in
\mathcal Q_{s,r,\aaa}^T(N_0,\theta,\mu)$ and
\begin{equation}\label{topstim}
 \|Y\|_{s,r,\aaa, N_0,\theta,\mu,\mathcal O_\rro}^T
 \lessdot (  8r/R )^{d-2}\|X\|_{R/2,\aaa, N_0,\theta,\mu'}^T \, .
 \end{equation}
\end{proposition}

The proof of Proposition \ref{qtop0} follows closely the analogous
Proposition 7.2 in \cite{BBP1} (replacing the Hamiltonians with
the vector fields). The following  lemma holds (see Lemma 7.11 in
\cite{BBP1}).

\begin{lemma}\label{artasersebis}
 Let $ X \in {\cal V}_{R/2,\aaa} $, $ Y:=\Phi_\star X$
   and $Y_0(x,y):=Y(x,y,0,0)-Y^{(y)}(x,0,0,0)\partial_y.$ Then,
assuming \eqref{condro}, $ \|Y_0\|_{s,2r,\aaa,\mathcal O_\rro\cup
\mathcal O_{2\rro}} \lessdot (R/r) \|X\|_{R/2,\aaa} $.
\end{lemma}

Recalling \eqref{ACAN} the vector field $ \mathcal N_0+G_1+G_2 +
G_3 $ in \eqref{Birkhoff} is transformed by the change of variable
\eqref{AAV} into
\begin{equation}\label{argiletum}
    \Phi_\star (\mathcal N_0+G_1+G_2 + G_3)=\mathcal N+\mathcal P=
  \mathcal N+\mathcal P_1+\PP_2+\PP_3
\end{equation}
where the normal form $\mathcal N$ is as in \eqref{90210} with
frequencies (satisfying \eqref{sio})
 as in \eqref{omegaxi}
\begin{eqnarray}\label{attila}
&&\o_j(\xi)=\l_j+ \l_j^{-1} \Big(-\frac14 \vec a_{|j|}\xi_{|j|} +
\vec a\cdot\xi\Big)
  \, ,
\forall j \in {\cal I}, \qquad \OO_j (\xi ) = \l_j + \l_j^{-1}
\,\vec a \cdot \xi   \, , \forall j \notin {\cal I} \, ,
\\
&& \vec a:= {\mathop\sum}_{1\leq l\leq 3} \kappa_l \vec
a^{(l)}\in\R^{n/2}\,,\ \ \vec a^{(1)}_i:=-\l_i^{-2}\,,\  \vec
a^{(2)}_i:=-i^2\l_i^{-2}\,,\  \vec a^{(3)}_i:=-1\,,\ 
\forall\, i\in\mathcal I^+\,.\  \ \label{attilabis}
\end{eqnarray}
Moreover the three terms of the perturbation are
\begin{eqnarray}
    \mathcal P^{(x_j)}_1 \!\!\!\!
    &:=& \!\!\!\!
    \frac{1}{\l_j} \left(-\frac14 \vec a_{|j|} \cdot y_j +
\frac12(\vec a,\vec a) y\right),\, \mathcal P^{(y_j)}_1=0,  \, j
\in\mathcal I, \,
    \mathcal P^{(z_j^\s)}_1 :=
    -\frac{\s \ii}{2\l_j}
   (\vec a,\vec a)\cdot y\  z_j^\s
   \,, \  \s=\pm\,,\  j\not\in\mathcal I,
\nonumber
\\
\PP_2 \!\! &:=& \!\! \Phi_\star G_2 \  \ {\rm (note\  that\
}\PP_2^{(x)}=\PP_2^{(y)}=0\,,\ \PP_2^{(z^\pm_j)}=G_2^{(u^\pm_j)},
j\notin \mathcal I
 {\rm )}, \quad
\PP_3 := \Phi_\star G_3
    \,.
\label{atrocita}
\end{eqnarray}

As in \eqref{primopezzo} we decompose the perturbation 
\begin{equation}\label{serse}
\PP=\PP^y (x;\xi) \partial_y+\PP_*\,,\qquad \PP^y(x;\xi)\partial_y
:=\Pi^{(-1)} \PP^{(y)} \partial_y =\Pi^{(-1)} \PP_3^{(y)}
\partial_y =\PP_3^{(y)}(x,0,0,0;\xi) \partial_y
 \, .
\end{equation}
\begin{lemma}\label{dario}
Let $s,r>0$ as in \eqref{condro} and $N$ large enough (w.r.t.
$\mm,{\cal I},L,b$). Then
\begin{equation}\label{milziadebis}
 \| \PP^y  \partial_y \|_{s,r, \aaa/2,
    {\cal O}}^\l
 \lessdot (1+\l/\rro) R^6 r^{-2} \,,
 \qquad
\| \PP_* \|_{\vec p}^T
 \lessdot (1+\l/\rro) (r^2 + R^5 r^{-1})\,,
\end{equation}
where
\begin{equation}\label{temistocle}
    \mathcal O=\mathcal{O}(\rro):=
   \big\{ \xi\in \R^n\  :\  2 \rro /3  \leq \xi_l\leq 3 \rro / 4\,,
   \ \  l=1,\ldots , n   \big\} \subset {\cal O}_\rro
\end{equation}
(the set ${\cal O}_\rro $ was defined in \eqref{paraO}) and $ \vec p
:= (s,r,\aaa/2,N,2,2,\l,\mathcal O) $.
\end{lemma}

\begin{pf}
By the definition \eqref{serse} we have
\begin{eqnarray}
\| \PP^y  \partial_y \|_{s,r, \aaa/2,   {\cal O}_\rro} &=&
\|\Pi^{(-1)} \PP_3^{(y)} \partial_y \|_{s,r, \aaa/2,   {\cal
O}_\rro}
 \stackrel{{\rm Lemma}\  \ref{lem:pro}} \leq
 \| \PP_3^{(y)} \partial_y
\|_{s,r, \aaa/2,   {\cal O}_\rro} \nonumber
\\
&\stackrel{\eqref{topstim},\eqref{atrocita}}  \lessdot&
  \Big(\frac{r}{R}
\Big)^{-2}\|  G_3  \|_{R/2,\aaa/2, N,
 7/4,\trecc}^T
\stackrel{\eqref{radiatorebis}}\lessdot \frac{R^6}{r^2}
\label{p00}
\end{eqnarray}
(applying \eqref{topstim} with $ d \rightsquigarrow 0 $, $ N_0
\rightsquigarrow N$,
 $ \teta  \rightsquigarrow 7/4 $, $ \mu
\rightsquigarrow  \duecc $,
 $ \mu'  \rightsquigarrow  \trecc $) and
taking $ N $ large enough so that \eqref{fuso} holds and $ N \geq
N_0' $ defined in Proposition \ref{BNF}.
 By \eqref{argiletum}, \eqref{atrocita} and \eqref{serse}
  we write
\be\label{barP} \PP_* = \PP_1+\PP_2+\PP_4+\PP_5 \qquad {\rm where}
\ee
$$
\PP_4 :=  \PP_3 (x,y,z, \bar z; \xi) - \PP_3 (x,y,0,0; \xi) \, ,
\quad \PP_5 :=  \PP_3 (x,y,0,0; \xi) -  \PP_3^{(y)} (x,0,0, 0;
\xi)\partial_y \, .
$$
We claim that \be\label{coma1}
\|\PP_1\|_{s,r,\aaa/2,N,\CCterzi,\duecc,\mathcal O_\rro}^T,\quad
\|\PP_2\|_{s,r,\aaa/2,N,\CCterzi,\duecc,\mathcal O_\rro}^T \quad
 \lessdot \quad r^2 
\, . \ee Indeed the estimate on $\PP_1$ follows since $\PP_1$ is
T\"oplitz and $\|\PP_1\|_{s,r,\aaa/2,\mathcal O_\rro} \lessdot
r^2$ by \eqref{atrocita}. On the other hand the estimate on
$\PP_2$ follows by \eqref{atrocita} and \eqref{radiatoreter} with
$ N \geq N_0 $ large enough to fulfill \eqref{caracalla}.

 By \eqref{atrocita} and
  \eqref{topstim} (with $ d \rightsquigarrow 1 $, $ N_0
\rightsquigarrow N $, $ \mu  \rightsquigarrow  \duecc$,
 $ \mu'  \rightsquigarrow  \trecc $),
for $ N $ large enough,  we get
\be\label{coma3} \| \PP_4
\|_{s,r,\aaa/2,N,\CCterzi,\duecc,\mathcal O_\rro}^T \lessdot
\Big(\frac{r}{R} \Big)^{-1} \|G_3\|^T_{R/2,\aaa/2,N_0',7/4,\trecc}
\stackrel{\eqref{radiatorebis}}\lessdot \Big(\frac{r}{R}
\Big)^{-1} R^4  = \frac{R^5}{r} \, . \ee Since $\PP_5$ does not
depend on the variables $z^\pm$ we get \be\label{coma2} \| \PP_5
\|_{s,r,\aaa/2,N,\CCterzi,\duecc,\mathcal O_\rro}^T = \| \PP_5
\|_{s,r,\aaa/2,\mathcal O_\rro} \stackrel{\text{Lemma}
\ref{artasersebis}}\lessdot \Big(\frac{r}{R} \Big)^{-1}\| G_3
\|_{R/2,\aaa/2} \stackrel{\eqref{radiatorebis}}\lessdot
\Big(\frac{r}{R} \Big)^{-1} R^4  = \frac{R^5}{r} \, . \ee In
conclusion, by \eqref{barP}, \eqref{coma1}, \eqref{coma3},
\eqref{coma2} we get $ \| \PP_*
\|_{s,r,\aaa/2,N,\CCterzi,\duecc,\mathcal O_\rro}^T \lessdot r^2+
R^5 r^{-1} $. In order to prove the
 estimates  \eqref{milziadebis} we have to prove  Lipschitz
 estimates
(see \eqref{filadelfia}, \eqref{fanfulla}). We first note that
the vector fields $\PP^y\partial_y  $ and $ \PP_* $ are
\textsl{analytic} in the parameters $\xi\in\mathcal{O}_\rro $.
Then we apply Cauchy estimates in the subdomain $\mathcal
O=\mathcal{O}(\rro)\subset \mathcal{O}_\rro $ (see
\eqref{temistocle}), noting that $\rro\lessdot {\rm dist}(\mathcal
O, \partial \mathcal{O}_\rro).$ Then
$$
\| \PP_* \|^{{\rm lip}}_{s,r,\aaa/2, {\cal O}} \lessdot \rro^{-1}
\| \PP_* \|_{s,r,\aaa/2, {\cal O}_\rro} \qquad {\rm and} \qquad \|
\PP^y\partial_y \|^{{\rm lip}}_{s,r,\aaa/2, {\cal O}} \lessdot
\rro^{-1} \| \PP^y\partial_y \|_{s,r,\aaa/2, {\cal O}_\rro}
$$
and \eqref{milziadebis} is proved.
\end{pf}

\smallskip

We now verify that the assumptions of Theorems
 \ref{thm:IBKAM}-\ref{thm:measure} are fulfilled by
 $ \mathcal N+\PP $ in \eqref{argiletum} with parameters
 $ \xi \in {\cal O}(\rro) $ defined in \eqref{temistocle}.
 Note that the sets ${\cal O}=[\rho/2,\rho]^n $
 defined in Theorem \ref{thm:measure}
 and ${\cal O}(\rro)$ defined in \eqref{temistocle}
 are diffeomorphic through $\xi_i\mapsto (7\rho +2\xi_i)/12$.
The frequency $\vec \o$ (recall \eqref{pingu}) defined in
\eqref{attila} has the form $\vec\o=\bar\o+A\xi$ in
\eqref{omegaxi}
  with
\begin{equation}\label{platea}
    \ A :=
    {\mathop\sum}_{1\leq l \leq 3} \kappa_l A^{(l)}\,,\qquad
     A^{(l)} := \Big({\rm diag}_{i\in\mathcal I^+} \l_i^{-1} \Big)
     \Big( -\frac14 {\rm Id}_{n/2}+ {\bf 1}_{n/2}\Big)
     \Big(  {\rm diag}_{i\in\mathcal I^+} \vec a_i^{(l)} \Big) \,,
\end{equation}
denoting by ${\bf 1}_{n/2}$ the $ (n/2) \times (n/2) $ matrix with
all entries equal to $1.$ Then $\vec \oo$ and $\OO_j$, defined in
\eqref{attila} satisfy \eqref{omegaxi} and
 hypotheses (A1)-(A2) follow.
Moreover (A3)-(A4) and the quantitative bound \eqref{KAMcondition}
follow by \eqref{milziadebis},
choosing 
\begin{equation}\label{leonida}
s = 1,\  r = R^{1 + \frac34}\,,
 \ \rro = C_* R^2 \  {\rm as \ in \
}\eqref{condro},  
\ \ N \ {\rm as \ in \ Lemma \ }\ref{dario},\
  \theta = \CCterzi, \
   \mu = \duecc,
   \
\gamma = R^{3 + \frac15}
\end{equation}
and taking $ R $ \textsl{small enough}. Hence Theorem
\ref{thm:IBKAM} applies.

Let us verify that also the assumptions  of Theorem
\ref{thm:measure} are fulfilled. Since $ {\bf 1}_{n/2}^2= (n/2)
{\bf 1}_{n/2}$
by \eqref{platea} we get that 
the matrix $ A $ is invertible with
\begin{equation}\label{plateabis}
 A^{-1}:= \Big({\rm diag}_{i\in\mathcal I^+} 1 / \vec a_i 
 \Big)
     \Big( -4 {\rm Id}_{n/2}+\frac{16}{2n-1} {\bf 1}_{n/2}\Big)
     \Big({\rm diag}_{i \in \mathcal I^+} \l_i  \Big) \,,
    \end{equation}
 for 
all $ \kappa_1, \kappa_2, \kappa_3 $ such that 
$\vec a_i\stackrel{\eqref{attilabis}} {:=}{ \mathop \sum}_{1\leq l\leq 3}\kappa_l \vec a_i^{(l)} = - ( \kappa_1 + \kappa_2 i^2 + \kappa_3 \lambda_i^2  ) \lambda_i^{-2} 
\neq 0 $, $\forall i \in \mathcal I^+$, see \eqref{NRF}. Moreover, by
\eqref{attilabis} and \eqref{plateabis} we have$  \big(A^T\big)^{-1}  \vec  a =  4 {\bar \om}   \slash (2n-1) $  and 
then condition \eqref{pota} is equivalent to 
\eqref{potapota} (note that $  \big(A^T\big)^{-1}  \vec  a $ does {\it not} depend on  $ \kappa_1, \kappa_2, \kappa_3 $).

Finally we deduce that the Cantor  set of parameters $
\mathcal{O}_\infty \subset \mathcal{O} $ in \eqref{Cantorinf} has
asymptotically full density because
$$
\frac{|\mathcal{O} \setminus\mathcal{O}_\infty|}{| \mathcal{O}|}
\stackrel{\eqref{consolatrixafflictorum}} \lessdot \rho^{-1}
\g^{2/3}\stackrel{\eqref{leonida}}\lessdot R^{-2}
R^{\frac23 (3+\frac15)}=R^{\frac{2}{15}} \to 0 \, . 
$$
The proof of Theorem \ref{thm:DNLW} is now completed.

\begin{remark}\label{carmelo}
If  $g=g^{(=3)}+g^{(\geq 4)} $ (unlike \eqref{gleading})
then $\|G_3\|_{R/2,\aaa/2,N_0',7/4,3}^T
 \lessdot R^3 $ which does not fit the smallness condition
 of
 Theorem \ref{thm:IBKAM}. The  term of order four should be removed by a further step of
Birkhoff normal form.
If the term $g^{(=3)}$ depends on the space variable $\mathtt x $ nothing changes except
to check
the twist condition, see \eqref{attila}, \eqref{plateabis}. 
For simplicity,  we did not pursue these  points.
\end{remark}

\noindent
\small
{\it Massimiliano Berti}, Dipartimento di Matematica e
Applicazioni ``R. Caccioppoli", Universit\`a degli Studi di Napoli
Federico II,  Via Cintia, Monte S. Angelo,
I-80126, Napoli, Italy,  {\tt m.berti@unina.it}. 
\\[1mm]
{\it Luca Biasco}, Dipartimento di Matematica e Fisica, Universit\`a
``Roma Tre'', Largo S. L. Murialdo 1, I-00146, Roma, Italy, {\tt
biasco@mat.uniroma3.it}.
\\[1mm]
{\it Michela Procesi}, Dipartimento di Matematica ``G.
Castelnuovo", Universit\`a degli Studi di Roma la Sapienza,  P.le
A. Moro 5, 00185 ROMA, Italy,  {\tt
michela.procesi@mat.uniroma1.it}.
\\[1mm]
This research was supported by the European Research Council under
FP7
and partially by the PRIN2009
grant ``Critical point theory and perturbative methods for
nonlinear differential equations".

\end{document}